\documentclass{article}

\usepackage {amsmath}
\usepackage {hyperref}
\usepackage {amsfonts}
\usepackage {amsthm}
\usepackage[flushmargin]{footmisc}
\usepackage[parfill]{parskip}
\usepackage {amssymb}
\usepackage {fullpage}
\usepackage {mathrsfs}
\usepackage {amscd}
\usepackage{bbm}
\usepackage[all,cmtip]{xy}
\DeclareMathAlphabet{\mathpzc}{OT1}{pzc}{m}{it}

\title{Character sheaves on neutrally solvable groups}
\author{Tanmay Deshpande}
\date{}

\newtheorem {thm} {Theorem} [section]
\newtheorem {prop} [thm] {Proposition}
\newtheorem {conj} [thm] {Conjecture}
\newtheorem {lem} [thm] {Lemma}
\newtheorem {cor} [thm] {Corollary}

\theoremstyle{definition}
\newtheorem {defn} [thm] {Definition}

\newtheorem {prob} [thm]  {Problem}

\newtheorem {rk} [thm]  {Remark}
\newtheorem {ex} [thm] {Example}

\newcommand{\beq}{\begin{equation}}
\newcommand{\eeq}{\end{equation}}
\newcommand{\bthm}{\begin {thm}}
\newcommand{\ethm}{\end {thm}}
\newcommand{\bprop}{\begin {prop}}
\newcommand{\eprop}{\end {prop}}
\newcommand{\bprob}{\begin {prob}}
\newcommand{\eprob}{\end {prob}}
\newcommand{\bcor}{\begin {cor}}
\newcommand{\ecor}{\end {cor}}
\newcommand{\blem}{\begin{lem}}
\newcommand{\elem}{\end{lem}}
\newcommand{\bdefn}{\begin{defn}}
\newcommand{\edefn}{\end{defn}}
\newcommand{\brk}{\begin{rk}}
\newcommand{\erk}{\end{rk}}

\renewcommand{\subset}{\subseteq}
\newcommand{\xto}{\xrightarrow}
\newcommand{\onto}{\twoheadrightarrow}

\newcommand{\Perv}{\hbox{Perv}}

\newcommand{\av}{{\hbox{av}}}

\newcommand{\tMg}{\tM_{Ug,e}}

\newcommand{\tMbeel}{\tM_{G,\beel}}
\newcommand{\Mbeel}{\M_{G,\beel}}
\newcommand{\uast}{\underline{\ast}}

\renewcommand {\bar} {\overline}
\newcommand{\bpf}{\begin{proof}}
\newcommand{\epf}{\end{proof}}
\newcommand{\bex}{\begin{ex}}
\newcommand{\eex}{\end{ex}}
\newcommand{\rar}[1]{\stackrel{#1}{\longrightarrow}}
\newcommand{\f}{\mathbb}

\newcommand{\<}{\langle}
\renewcommand{\>}{\rangle}

\newcommand{\h}{\operatorname}
\newcommand{\e}{\operatorname{ev}}
\renewcommand{\c}{\operatorname{coev}}
\newcommand{\FPdim}{\operatorname{FPdim}}
\newcommand{\ind}{\operatorname{ind}}

\newcommand{\Sh}{\operatorname{Sh}}
\newcommand{\ad}{\operatorname{ad}}

\newcommand{\indg}{\h{ind}_{G'}^G}

\renewcommand{\L}{\mathcal{L}}
\newcommand{\R}{\mathscr{R}}
\renewcommand{\P}{\mathscr{P}}

\newcommand{\Padm}{{\mathscr{P}_{\h{adm}}(G)}}
\newcommand{\N}{\mathcal{N}}

\newcommand{\cM}{\mathcal{M}}

\newcommand{\Fq} {\mathbb{F}_q}

\newcommand{\Fpcl} {\overline{\mathbb{F}}_p}
\newcommand{\Fqcl} {\overline{\mathbb{F}}_q}
\newcommand{\Q} {\mathbb{Q}}

\renewcommand{\r} {\h{r}}
\newcommand{\F} {\mathbb{F}}
\newcommand{\un} {\mathbbm{1}}

\newcommand{\Dmon} {D^{\leq 0}_{mon}\Vec}

\newcommand{\normal} {\triangleleft}
\renewcommand{\t}{\widetilde}
\renewcommand{\phi} {\varphi}

\renewcommand{\k} {\mathtt{k}}
\newcommand{\K} {\mathcal{K}}
\newcommand{\hG} {\widehat{G}}
\newcommand{\lG} {\f{L}({G})}
\newcommand{\leG} {\f{L}_e({G})}
\newcommand{\lL} {\f{L}}
\newcommand{\M} {\mathscr{M}}

\newcommand{\D} {\mathscr{D}}

\newcommand{\C} {\mathcal{C}}

\newcommand{\Hom} {\h{Hom}}

\renewcommand{\Vec}{\h{Vec}}

\newcommand{\Qlcl} {\overline{\mathbb{Q}}_{\ell}}

\newcommand{\T}{\h{Tr}}

\newcommand{\pt}{\hbox{pt}}

\newcommand{\eel} {e\boxtimes e_\L}
\newcommand{\beel} {\bar{e_\L}}

\renewcommand{\tilde}{\widetilde}
\newcommand{\noin}{\noindent}
\newcommand{\A} {\mathcal{A}}
\newcommand{\fL} {\mathbb{L}}

\newcommand{\ta} {\tilde{a}}
\newcommand{\tb} {\tilde{b}}

\newcommand{\g} {{\gamma}}
\renewcommand{\l} {{\lambda}}

\newcommand{\id}{\operatorname{id}}

\newcommand{\Mfe}{\M_{G,f}}
\newcommand{\tM}{\widetilde{\M}}

\newcommand{\bit}{\begin{itemize}}
\newcommand{\eit}{\end{itemize}}
\newcommand{\loc}{\h{loc}}

\newcommand{\Irrep}{\h{Irrep}}
\newcommand{\Irrepe}{\h{Irrep}_{e}(G,F)}

\newcommand{\avg}{\av_{G/G'}}
\newcommand{\DG}{\D_G(G)}

\newcommand{\FunG}{\Fun([G],F)}
\newcommand{\DGp}{\D_{G'}(G')}

\newcommand{\DGn}{\D_{G_0}(G_0)}

\newcommand{\can}{\mathbb{K}}

\newcommand{\B}{\mathcal{B}}
\renewcommand{\O}{\mathcal{O}}

\newcommand{\Z}{\mathbb{Z}}

\newcommand{\tr}{\h{tr}}
\newcommand{\Fun}{\h{Fun}}
\newcommand{\CS}{\h{CS}}
\newcommand{\Dmu}{{\D_{\mu_3U}(\mu_3U)}}
\newcommand{\muU}{{\mu_3U}}

\newcommand{\bconj}{\begin{conj}}
\newcommand{\econj}{\end{conj}}
\begin{document}
\maketitle

\begin{abstract}
\noin Let $G$ be an algebraic group over an algebraically closed field $\k$ of characteristic $p>0$. In this paper we develop the theory of character sheaves on groups $G$ such that their neutral connected components $G^\circ$ are solvable algebraic groups. For such algebraic groups $G$ (which we call neutrally solvable) we will define the set $\CS(G)$ of character sheaves on $G$ as certain special (isomorphism classes of) objects in the category $\DG$ of $G$-equivariant $\Qlcl$-complexes (where we fix a prime $\ell\neq p$) on $G$. We will describe a partition of the set $\CS(G)$ into finite sets known as $\f{L}$-packets and we will associate a modular category $\M_L$ with each $\f{L}$-packet $L$ of character sheaves using a truncated version of convolution of character sheaves. In the case where $\k=\Fqcl$ and  $G$ is equipped with an $\Fq$-Frobenius $F$ we will study the relationship between $F$-stable character sheaves on $G$ and the irreducible characters of (all pure inner forms of) $G^F$. In particular, we will prove that the notion of almost characters (introduced by T. Shoji using Shintani descent) is well defined for neutrally solvable groups and that these almost characters coincide with the ``trace of Frobenius'' functions associated with $F$-stable character sheaves. We will also prove that the matrix relating the irreducible characters and almost characters is block diagonal where the blocks on the diagonal are parametrized by $F$-stable $\f{L}$-packets. Moreover, we will prove that the block in this transition matrix corresponding to any $F$-stable $\f{L}$-packet $L$ can be described as the crossed S-matrix associated with the auto-equivalence of the modular category $\M_L$ induced by $F$.
\end{abstract}

\tableofcontents

\section{Introduction}\label{s:i}
In \cite{De2} we embarked upon the journey towards developing a theory of character sheaves on a general affine algebraic group $G$ defined over an algebraically closed field $\k$ of characteristic $p>0$. The study of character sheaves was initiated by Lusztig for reductive groups in his works \cite{L}. Later Boyarchenko and Drinfeld developed the theory of character sheaves on unipotent groups in \cite{BD}, \cite{B}. In this paper, we develop the theory of character sheaves on neutrally solvable groups, namely groups $G$ such that their neutral connected components $G^\circ$ are solvable algebraic groups.

Let us fix two distinct primes $p,\ell$. In this paper, $\k$ will always denote an algebraically closed field of characteristic $p$. One case of particular interest is when $\k=\Fpcl$. All algebraic groups and schemes will be assumed to be over $\k$ unless mentioned otherwise. It is often more convenient to pass to the perfectizations of algebraic groups and schemes. Hence continuing our convention from \cite{De2}, by algebraic group we actually mean perfect quasi-algebraic group over $\k$ (i.e. the perfectization of an algebraic group). We refer to \cite[\S1.9]{BD} for more about this convention. 

We shall continue to use all the standard notations and conventions from \cite{BD}, \cite{De2}. Hence for example if $G$ is a (perfect quasi-) algebraic group, then $\DG$ denotes the $\Qlcl$-linear triangulated ribbon $\r$-category of conjugation equivariant constructible $\Qlcl$-complexes. In particular $\DG$ is a braided monoidal category (under convolution with compact supports) with a weak form of  duality (denoted by $\f{D}^{-}$) and for each $C\in \DG$ we have a functorial twisting automorphism $\theta_C:C\rar{\cong}C$. We also recall from \cite[Appendix A]{BD} that the (weak) duality functor $\f{D}^{-}$ is defined as the composition $\iota^*\circ \f{D}=\f{D}\circ \iota^*$ where $\f{D}$ denotes Verdier duality and $\iota:G\rar{} G$ denotes the inversion map. We let $\can_G$ denote the dualizing complex on $G$. In this paper we will always assume that our perfect quasi-algebraic group $G$ is defined as the perfectization of some specific algebraic group over $\k$. With this assumption, we have a canonical identification $\can_G\cong \Qlcl[2\dim G](\dim G)$.

Now suppose that $G$ is neutrally solvable. In this paper we will define a set $\CS(G)$ of isomorphism classes of some special objects in $\DG$ known as character sheaves on $G$. If $C\in \CS(G)$, then we will see that $C$ is a simple object in $\DG$, in the sense that $\Hom_{\DG}(C,C)=\Qlcl$. In \cite{De2}, we studied minimal idempotents in the braided monoidal category $\DG$ and proved that each such minimal idempotent is in fact locally closed (cf. \cite[\S2.3.2]{De2}) and can be obtained from an admissible pair for $G$. Let $\hG$ denote the set of (isomorphism classes of) minimal idempotents in $\DG$. For each $e\in \hG$, we will first define the set $\CS_e(G)$ of character sheaves in the full subcategory $e\DG\subset \DG$. If $e,e'\in \hG$ are non-isomorphic minimal idempotents, then $e\DG\cap e'\DG=0$, and hence it will follow that $\CS_e(G)$ and $\CS_{e'}(G)$ are disjoint sets. We will then define $\CS(G)=\coprod\limits_{e\in \hG}\CS_e(G)$. Moreover, we will see that the set $\CS(G)$ can be defined purely in terms of the $\Qlcl$-linear triangulated braided monoidal category $\DG$. In particular, this means that any $\Qlcl$-linear triangulated braided auto-equivalence of $\DG$ preserves the set $\CS(G)$ of character sheaves. 

Here we remark that for a neutrally solvable group, the character sheaves may not be perverse sheaves (even up to shift). This is in contrast with the fact that character sheaves on reductive as well as unipotent groups are perverse sheaves (at least up to shift, depending upon convention).

To define the set $\CS_e(G)$ we will use the fact that the minimal idempotent $e$ can be obtained from an admissible pair for $G$. We will further define a partition of the set $\CS_e(G)=\coprod\limits_{f\in\f{L}_e(G)}{\CS_{e,f}}(G)$ into what we will call $\f{L}$-packets of character sheaves. The set $\f{L}_e(G)$ parametrizing $\f{L}$-packets associated  with $e$ is equal to the set of equivalence classes of minimal quasi-idempotents (as defined in Appendix \ref{a:qiltc}) in the $\Qlcl$-linear triangulated braided monoidal category $e\DG$. Moreover, each $\f{L}$-packet of character sheaves $\CS_{e,f}(G)$ (which will often be denoted simply as $\CS_f(G)$) is a finite set. Associated with each such $\f{L}$-packet, we will define a modular category $\M_{G,f}$ whose simple objects correspond to the finite set $\CS_f(G)=\CS_{e,f}(G)$. We will construct the tensor product in this modular category by defining a truncated version of convolution. We will prove that any $\Qlcl$-linear triangulated braided auto-equivalence $\Phi:\DG\rar{\cong}\DG$ induces a permutation of $\f{L}$-packets and of character sheaves. In fact we obtain an induced equivalence of braided fusion categories $\Phi:\M_{G,f}\rar{\cong}\M_{G,\Phi(f)}$ for each $\lL$-packet $\CS_{e,f}(G)$.

After defining character sheaves, $\f{L}$-packets and their associated modular categories, we will apply these results to the character theory of the associated finite groups. For this we take $\k=\Fpcl$ and consider an $\Fq$-Frobenius map $F:G\rar{}G$ where $q$ is some power of $p$. We would like to study the character theory of the finite group $G(\Fq)=G^F$ using our theory of character sheaves on $G$. However, as observed in \cite{B} and \cite{De3}, it is more natural to consider all the pure inner forms $G^{gF}$ at the same time. Here $g\in G$ and $gF:=\ad(g)\circ F:G\rar{}G$ is another $\Fq$-Frobenius map for $G$. The pure inner forms of $G^F$ are parametrized by the finite set $H^1(F,G)$ of $F$-twisted conjugacy classes in $G$. Note that if $G$ is connected, then by Lang's theorem $H^1(F,G)$ is singleton and we have only one pure inner form. We refer to \cite[\S2.4.1]{De3}, \cite[\S1.2]{De4} for more on these notions. 

We will also continue to use all the notation from \cite{De4}. In particular we have the set 
\beq
\Irrep(G,F):=\coprod\limits_{\<g\>\in H^1(F,G)}\Irrep(G^{gF})
\eeq
of irreducible characters (over $\Qlcl$) of all the pure inner forms and the commutative Frobenius $\Qlcl$-algebra
\beq
\Fun([G],F):=\prod\limits_{\<g\>\in H^1(F,G)}\Fun(G^{gF}/\sim)
\eeq
of class functions on all pure inner forms (under convolution of functions). This space has the standard Hermitian inner product with respect to which $\Irrep(G,F)\subset \FunG$ is an orthonormal basis. Note that we will often use the same symbol to denote an irreducible representation as well as its character.
\brk\label{r:complex}
We will fix an identification $\Qlcl\cong \f{C}$. Once we fix this, we also obtain an inclusion $\Q_p/\Z_p\hookrightarrow \Qlcl^\times$. We remark that this choice is not important since all the functions that will arise in this paper can be chosen so as to take values in the subfield of $\Qlcl$ where `complex conjugation' is defined unambiguously.
\erk
Recall that in \cite[\S3.2]{De4} we defined an analogue of Asai's twisting operator, namely a unitary operator $\Theta^*:\FunG\rar{\cong}\FunG$. This unitary operator is obtained using a certain permutation $\Theta$ of the rational conjugacy classes in all pure inner forms which preserves the geometric conjugacy classes. In general, the action of $\Theta^*$ on $\FunG$ can be rather complicated.

Also recall that in \cite{De4},  for each positive integer $m$, we have defined the $m$-th Shintani descent map
\beq
\Sh_m:\Irrep(G,F^m)^F\hookrightarrow \FunG
\eeq 
which is well defined only up to scaling by $m$-th roots of unity and that the image of this map is an orthonormal basis of $\FunG$ known as the $m$-th Shintani basis.

Suppose we have an object $C\in \DG$ and an isomorphism $\psi:F^*C\rar{\cong} C$. Then as in \cite[\S2.4.8]{De3}, we have its associated Frobenius trace function $\T_{C,\psi}\in\FunG$. Furthermore if $C$ is simple, in the sense that $\Hom_{\DG}(C,C)=\Qlcl$, then we proved in \cite[\S3.3]{De4} that $\Theta^*(\T_{C,\psi})=\theta_C\cdot \T_{C,\psi}$, where $\theta_C\in \Qlcl^\times$ is the twist automorphism of $C\in\DG$. 

We have seen that $\Irrep(G,F)$ is an orthonormal basis of $\FunG$. We will prove that the Frobenius trace functions associated with $F$-stable character sheaves on $G$ also form an orthonormal basis of $\FunG$. By the previous paragraph, class functions in this orthonormal basis are eigenvectors for the twisting operator $\Theta^*$. On the other hand, the action of $\Theta^*$ on the basis $\Irrep(G,F)$ can be rather complicated.

Let us summarize the preliminary versions of the main results that we will prove in this paper. We will state and prove more precise statements later.
\bthm\label{t:main1}
Let $G$ be a neutrally solvable group over any algebraically closed field $k$ of characteristic $p$. We will define a set $\CS(G)$ of character sheaves in $\DG$ such that we have the following:\\
(i) Each character sheaf $C\in \CS(G)$ is a simple object in $\DG$, in the sense that $\Hom_{\DG}(C,C)=\Qlcl$. We have a partition \beq\CS(G)=\coprod\limits_{e\in \hG}\CS_e(G)\eeq where $\CS_e(G)$ are character sheaves in the full subcategory $e\DG\subset \DG$. The set $\CS_e(G)$ will be defined purely in terms of the $\Qlcl$-linear triangulated braided monoidal structure of the category $e\DG$.\\
(ii) For each minimal idempotent $e\in \DG$, we have a partition
\beq
\CS_e(G)=\coprod\limits_{f\in\f{L}_e(G)}\CS_{e,f}(G)
\eeq
into $\f{L}$-packets of character sheaves, where $\lL_e(G)$ is the set of equivalence classes of minimal quasi-idempotents (cf. Appendix \ref{a:qiltc}) in $e\DG$. Each $\f{L}$-packet $\CS_{e,f}(G)=\CS_f(G)$ \footnote{We will sometimes drop $e$ from the notation, since $e$ is determined by $f$ considered as a minimal quasi-idempotent in $\DG$.} is a finite set.\\
(iii) Let $\lL(G)$ denote the set of equivalence classes of minimal quasi-idempotents in $\DG$. Then $$\f{L}(G)=\coprod\limits_{e\in\hG}\f{L}_e(G).$$\\
(iv) Associated with each $\f{L}$-packet  $\CS_{e,f}(G)$ of character sheaves is a modular category $\M_{G,f}$ whose simple objects are the character sheaves in that $\f{L}$-packet.\\
(v) Let $\Phi:\DG\rar{\cong}\DG$ be any $\Qlcl$-linear triangulated braided monoidal auto-equivalence. Then clearly we have an induced permutation (also denoted by $\Phi$) of the set $\hG$ of minimal idempotents. Then $\Phi$ preserves the set of character sheaves on $G$ along with its idempotent and $\f{L}$-packet decompositions. Namely we have $\Phi(\CS_e(G))=\CS_{\Phi(e)}(G)$ for each $e\in \hG$ and $\Phi(\CS_{e,f}(G))=\CS_{\Phi(e),\Phi(f)}$. In particular, we have an induced permutation $\Phi:\f{L}(G)\rar{}\f{L}(G)$ of the set of $\fL$-packets such that $\Phi(\f{L}_e(G))=\fL_{\Phi(e)}(G)$ for each $e\in \hG$. Moreover for each $\f{L}$-packet, we have an induced equivalence
\beq
\Phi:\M_{G,f}\rar{\cong} \M_{G,\Phi(f)}
\eeq
of braided fusion categories.
\ethm 

One of the main motivations behind the study of character sheaves on an algebraic group $G$ is their relationship with the character theory of the finite groups of the form $G(\Fq)$. Our next main result is in this direction.

\bthm\label{t:main2}
In addition to the assumptions in Theorem \ref{t:main1}, suppose now that $\k=\Fqcl$ and that we have an $\Fq$-Frobenius map $F:G\rar{}G$. This induces a $\Qlcl$-linear triangulated braided monoidal auto-equivalence $F:={F^{*}}^{-1}:\DG\rar{\cong}\DG$ and hence an induced permutation $F$ of the set of character sheaves and $\fL$-packets. Then we have:\\
(i) The set $\CS(G)^F=\coprod\limits_{f\in\fL(G)^F}\CS_{f}(G)^F$ of Frobenius stable character sheaves is finite. The set $\f{L}(G)^F$ of $F$-stable $\fL$-packets is also finite and for each $f\in \fL(G)^F$, the set $\CS_{f}(G)^F$ is nonempty.\\
(ii) We have a partition of the set $\Irrep(G,F)$ of irreducible characters of all pure inner forms of $G^F$ in terms of $\fL(G)^F$:
\beq
\Irrep(G,F)=\coprod\limits_{f\in \fL(G)^F}\Irrep_{f}(G,F).
\eeq
The sets $\Irrep_{f}(G)$ for $f\in \lL(G)^F$ are all non-empty, and are known as $\fL$-packets of irreducible characters. \\
(iii) For each $C\in \CS(G)^F$, let us choose $\psi_C:F^*C\rar{\cong} C$ such that $||\T_{C,\psi_C}||=1$. Then the set $\{\T_{C,\psi_C}|C\in \CS(G)^F\}$ is an orthonormal basis of $\FunG$ consisting of $\Theta^*$-eigenvectors. In particular, $|\Irrep(G,F)|=|\CS(G)^F|$. For each $f\in \fL(G)^F$, the $\Qlcl$-linear spans in $\FunG$ of the two sets $\{\T_{C,\psi_C}|C\in \CS_{f}(G)^F\}$ and $\Irrep_{f}(G,F)$ are equal. In particular, we have $|\Irrep_{f}(G,F)|=|\CS_{f}(G)^F|$.  We also deduce from this that the transition matrix between irreducible characters and Frobenius trace functions of $F$-stable character sheaves is block diagonal with blocks parametrized by the set $\fL(G)^F$.\\
(iv) For each $f\in \fL(G)^F$, we have a modular auto-equivalence $F:\M_{G,f}\rar{\cong}\M_{G,f}$ of the modular category $\M_{G,f}$. Then the unitary transition matrix between the two sets $\{\T_{C,\psi_C}|C\in \CS_{f}(G)^F\},\Irrep_{f}(G,F)\subset \FunG$ is given by the crossed S-matrix, as defined in \cite{De5}, (suitably normalized) associated with the modular auto-equivalence $F:\Mfe\rar{\cong}\Mfe$.
\ethm 
From this result, for each positive integer $m$, we obtain an $\fL$-packet decomposition of the set $\Irrep(G,F^m)$ with $\fL$-packets parametrized by $\fL(G)^{F^m}$. We also obtain that
\beq
\Irrep(G,F^m)^F=\coprod\limits_{f\in \fL(G)^F\subset \fL(G)^{F^m}}\Irrep_{f}(G,F^m)^F.
\eeq

Finally we relate the theory of character sheaves on neutrally solvable groups to Shintani descent:

\bthm\label{t:main3}
We continue in the setting of Theorem \ref{t:main2}. Then we have:\\
(i) Shintani descent respects the $\fL$-packet decomposition of irreducible characters, namely for each positive integer $m$ and each $f\in \fL(G)^F$, $\Sh_m(\Irrep_{f}(G,F^m)^F)$ is a basis of the subspace of $\FunG$ spanned by $\Irrep_{f}(G,F)$ (or equivalently of the subspace spanned by  $\{\T_{C,\psi_C}|C\in \CS_{f}(G)^F\}$). \\
(ii) There exists a positive integer $m_0$ such that for any positive integer $m$, the $m$-th Shintani basis of $\FunG$ only depends (up to scaling by roots of unity) on the residue of $m$ modulo $m_0$. The $m$-th Shintani basis (well defined up to scaling by roots of unity) where $m$ is any positive multiple of $m_0$ is known as the basis of almost characters of $\FunG$.\\
(iii) The orthonormal basis of almost characters of $\FunG$ agrees (up to scaling by roots of unity) with the orthonormal basis  $\{\T_{C,\psi_C}|C\in \CS(G)^F\}$ formed by Frobenius trace functions associated with $F$-stable character sheaves on $G$. In particular, almost characters are eigenvectors of the twisting operator $\Theta^*$.
\ethm

We conjecture that there should be an interesting theory of character sheaves on general algebraic groups:
\bconj\label{c:main}
Analogues (see Remark \ref{r:caveat} below) of Theorems  \ref{t:main1}, \ref{t:main2} and \ref{t:main3} hold for any algebraic group $G$ over $\k$.
\econj
As we have seen in \cite{De2} and as we will see from the arguments of this paper, we can essentially reduce this conjecture to the Heisenberg case.

\brk\label{r:caveat}
For a general algebraic group $G$, we do not expect the $\f{L}$-packets of character sheaves to be parametrized by equivalence classes of minimal quasi-idempotents in $\DG$. Nevertheless, we do expect that associated with every $\f{L}$-packet will be a modular category. For example, if $G$ is a reductive group, Lusztig has defined certain modular categories associated with $\f{L}$-packets using a notion of truncated convolution of character sheaves. In fact, many aspects of Conjecture \ref{c:main} are known to hold  for reductive groups (sometimes under some further restrictions) by the work of Lusztig, Shoji and others. 
\erk

Our approach to the theory of character sheaves on neutrally solvable groups is based on the notion of minimal quasi-idempotents in $\Qlcl$-linear braided triangulated categories. We describe this notion in the Appendix \ref{a:qiltc}. In \S\ref{s:cst}, we begin by describing the theory of character sheaves on a torus in terms of the notion of minimal quasi-idempotents.  We will then use this to study the case of Heisenberg idempotents in \S\ref{s:thc}, \S\ref{s:icahi}. In \S\ref{s:thc}, we will define character sheaves in the Heisenberg case and prove an analogue of Theorem \ref{t:main1} in this case (see Corollary \ref{c:modcat}, Definition \ref{d:csfg} and Remark \ref{r:csedgg}). We will prove in \S\ref{s:pcshc}, that in the Heisenberg case, the character sheaves are in fact perverse up to a certain shift. (This statement does not necessarily hold in the general case.) In \S\ref{s:icahi}, we study the relationship between irreducible characters and $F$-stable character sheaves in the Heisenberg case and prove the Heisenberg case versions of Theorem \ref{t:main2} (see \S\ref{s:rbiccs}, in particular Theorem \ref{t:relheis}) and Theorem \ref{t:main3}  (see \S\ref{s:sdthc}, in particular Corollary \ref{c:tsh}). In \S\ref{s:pmr}, we use the fact that all minimal idempotents in $\DG$ come from some admissible pair to prove our main results in the general case using the results already proved in the Heisenberg case. Theorem \ref{t:main1} is proved in \S\ref{s:dcs} (see Proposition \ref{p:lleg}, Definition \ref{d:csg} and Proposition \ref{p:main1iv}). In \S\ref{s:irap} and \S\ref{s:grpag} we study the set $\Irrep(G,F)$ using $F$-stable geometric conjugacy classes of admissible pairs for $G$. Theorem \ref{t:main2} is proved in \S\ref{s:ircs} and \S\ref{s:c} (see Theorems \ref{t:maingen} and \ref{t:fstap}). In \S\ref{s:sdgc}, we study Shintani descent in the general case and prove Theorem \ref{t:main3} (see Theorems \ref{t:sdgc} and \ref{t:fstap}). In \S\ref{s:example} we study in detail an instructive example, namely the character sheaves on the Borel subgroup of $SL_3$.

\section*{Acknowledgments}
I am grateful to V. Drinfeld for introducing me to the theory of character sheaves on algebraic groups and to T. Shoji for many useful discussions. I thank A. Beilinson for very helpful correspondence. This work was partially supported by World Premier Institute Research Center Initiative (WPI), MEXT, Japan.

\section{Character sheaves on a torus}\label{s:cst}
As a warm-up, we now describe a new way to look at character sheaves on a torus which is more suited for our approach to character sheaves on neutrally solvable groups. A character sheaf on a torus is just a multiplicative local system on $T$. However, in this section we will characterize the character sheaves on a torus $T$ using the notion of minimal quasi-idempotents (cf. Appendix \ref{a:qiltc}) in $\D(T)$ (or in $\D_T(T)$).  As a consequence, we will obtain a description of character sheaves on $T$ purely in terms on the $\Qlcl$-linear triangulated monoidal category $\D(T)$ (or $\D_T(T)$). We refer to Appendix \ref{a:qiltc} for the definition and properties of the abstract categorical notion of quasi-idempotents in $\Qlcl$-linear triangulated monoidal categories.

In this section, we will establish a canonical bijection between the set of character sheaves on a torus $T$ and the set of equivalence classes of minimal quasi-idempotents (cf. Defn. \ref{d:minqi}) in the categories $\D(T)$ and $\D_T(T)$. Let $\can_T\cong \Qlcl[2\dim T](\dim T)$ denote the dualizing complex on $T$. In this section we will continue to use the notations from \cite[\S7]{De2}. In particular $\C(T)$ denotes the $\Qlcl$-scheme whose $\Qlcl$-points parametrize the multiplicative local systems on $T$ and $\mathcal{M}_!:\D(T)\rar{}D^b_{coh}(\C(T))$ denotes the Mellin transform. We refer to \cite{GL} for a detailed discussion about these notions. 

\blem\label{l:mellin}
(i) Let $M\in \D(T)$ and let $\L$ be a multiplicative local system on $T$. Then we have a natural isomorphism $\L\ast M\rar{\cong}(\L\ast M)_1\otimes \L$ in $\D(T)$, where $(\L\ast M)_1\in D^b\Vec$ denotes the stalk of $\L\ast M$ at $1\in T$. (Here $\otimes$ is used in the sense of the previous section, or equivalently we may identify $D^b\Vec$ with constant sheaves in $\D(T)$ and then $\otimes$ denotes the usual (derived) tensor product in $\D(T)$.)\\
(ii) If $M\in \D(T)$ is nonzero, then there exists a multiplicative local system $\L$ on $T$ such that $\L\ast M\neq 0$.
\elem 
\bpf
Statement (i) follows readily from the definition of convolution and the fact that $\L$ is a multiplicative local system, i.e. $\mu^*\L\cong \L\boxtimes \L$ where $\mu:T\times T\to T$ is the multiplication in the torus. 

To prove (ii), we will use the Mellin transform and its properties that are studied in \cite{GL}. Suppose that $M$ is any nonzero object in $\D(T)$. Then by \cite[Prop. 3.4.5]{GL}, its Mellin transform $\cM_!(M)$ is also nonzero. In particular there exists a multiplicative local system $\L$ such that the (derived) pullback $i^*_{\L^{-1}}\cM_!(M)$ is nonzero where $i_{\L^{-1}}:\{\L^{-1}\}\hookrightarrow \C(T)$ denotes the inclusion of the closed point $\L^{-1}\in \C(T)$. Now by definition of the Mellin transform, $\cM_!(\L)$ is only supported at the closed point $\L^{-1}\in \C(T)$.
In particular this means that $\cM_!(\L)\otimes \cM_!(M)$ is nonzero. But since Mellin transform takes convolution to tensor product, we conclude that $\L\ast M$ must also be nonzero as desired.
\epf

\bthm\label{t:torcs}
Let $\L\in \D(T)$ be a multiplicative local system (which may also be naturally considered as an object of $\D_T(T)$). Then $e_\L:=\L \otimes \can_T$ is a minimal quasi-idempotent in $\D(T)$ (resp. $\D_T(T)$). Moreover, if $e\in \D(T)$ (resp. $e\in \D_T(T)$) is any minimal quasi-idempotent then there exists a unique (up to isomorphism) multiplicative local system $\L$ on $T$ such that $e_\L\ast e\neq 0$. In other words, $\L\mapsto e_\L$ defines a bijection between $\C(T)(\Qlcl)$ and the set of minimal quasi-idempotents in $\D(T)$ (resp. $\D_T(T)$) up to equivalence.
\ethm
\bpf
Using the isomorphism $\mu^*\L\cong \L\boxtimes \L$ we can define a natural isomorphism
\beq
\L\ast\L\rar{\cong}H^*_c(T,\Qlcl)\otimes \L
\eeq
and hence using the definition of $e_\L$, we obtain an isomorphism $e_\L\ast e_\L\rar\cong H^*_c(T,\can_T)\otimes e_\L$. Hence we see that $e_\L$ is a quasi-idempotent in $\D(T)$. By Lemma \ref{l:mellin}(i), for any $M\in \D(T)$ we have
$\L\ast M \rar{\cong} (\L\ast M)_1\otimes \L$ and hence 
\beq\label{e:elm}
e_\L\ast M\cong(\L\ast M)_1\otimes e_\L.
\eeq
In particular this holds for quasi-idempotents $M$. This proves that $e_\L$ is a minimal quasi-idempotent. 

On the other hand, let $e\in \D(T)$ be any minimal quasi-idempotent. In particular $e$ is nonzero and by Lemma \ref{l:mellin}(ii) there exists a multiplicative local system $\L$ such that $\L\ast e\neq 0$ and hence $e_\L\ast e\neq 0$. If $\L'$ is another multiplicative local system on $T$ which is not isomorphic to $\L$, then $\L\ast \L'=0$. The uniqueness of $\L$ (up to isomorphism) now follows. This completes the proof of the statements pertaining to the category $\D(T)$. 

The statements for $\D_T(T)$ can be proven in the same way by using a $T$-equivariant version of Lemma \ref{l:mellin}(i). It is clear that $e_\L\in \D_T(T)$ is a quasi-idempotent. Now let $e\in \D_T(T)$ be a $V$-quasi-idempotent such that $e_\L\ast e\neq 0$. Then we can define an isomorphism (cf. (\ref{e:elm})) 
\beq
e_\L\ast e\cong (\L\ast e)_1\otimes e_\L \mbox{ in $\D_T(T)$},
\eeq
where we consider the stalk $(\L\ast e)_1$ as an object of $\D_T(1)$. Now by Lemma \ref{l:wv'} (or rather by a slight generalization) we obtain that $(\L\ast e)_1\cong V\in D^b\Vec$. This proves that $e_\L\in \D_T(T)$ is a minimal quasi-idempotent. The remaining statements for $\D_T(T)$ now follow readily.
\epf

\brk\label{r:vt}
Let us denote the monic complex (see Definition \ref{d:monic}) $H^*_c(T,\can_T)\in \Dmon$ by $V_T$. Note that we have $V_T\cong H^*_c(T,\Qlcl)[2\dim T](\dim T)$. We have seen that for each multiplicative local system $\L$ on $T$, $e_\L$ is a $V_T$-quasi-idempotent in $\D(T)$.
\erk

\bcor\label{c:allminqi}
Each equivalence class of minimal quasi-idempotents in $\D(T)$ (or in $\D_T(T)$) contains a unique indecomposable object up to isomorphism, namely some $e_\L$ for a multiplicative local system $\L$ on $T$. Any minimal quasi-idempotent in $\D(T)$ (or in $\D_T(T)$) is of the form $W\otimes e_\L$ for some multiplicative local system $\L$ on $T$ and some $W\in \Dmon$.
\ecor
\bpf
Let $e\in \D(T)$ be any minimal quasi-idempotent, and say $e$ is a $V$-quasi-idempotent. Then by Theorem \ref{t:torcs} there exists a unique multiplicative local system $\L$ on $T$ such that $e\sim e_\L$ and clearly $e_\L$ is an indecomposable object of $\D(T)$ (as well as of $\D_T(T)$). Since $e,e_\L$ are equivalent minimal quasi-idempotents, we must have (cf. Lemma \ref{l:wv'})
\beq
e\ast e_\L\cong V_T\otimes e\cong V\otimes e_\L.
\eeq Now we see that the only indecomposable factors of the rightmost term in this isomorphism are $e_\L$ and its shifts. Hence the indecomposable factors of $e$ must also be of this form. In other words, we must have $e\cong W\otimes e_\L$ for some $W\in \Dmon$. This means that all minimal quasi-idempotents equivalent to $e_\L$ are of the form $W\otimes e_\L$ for some $W\in \Dmon$, in particular $e_\L$ is the only indecomposable minimal quasi-idempotent in its equivalence class.
\epf

\bdefn
A character sheaf in the $\Qlcl$-linear triangulated monoidal category $\D_T(T)$ (or in $\D(T)$) is an indecomposable minimal quasi-idempotent.
\edefn

Note that our convention is slightly different from the standard convention. In our convention, the quasi-idempotent $e_\L$ is said to be a character sheaf rather than the multiplicative local system $\L$. The results of this section prove that there is a canonical bijection $e_\L\longleftrightarrow \L$ between isomorphism classes of character sheaves on $T$ and of multiplicative local systems on $T$.

\bcor
Any triangulated monoidal auto-equivalence of $\D(T)$ (or of $\D_T(T)$) preserves the set of isomorphism classes of character sheaves on a torus $T$.
\ecor

\bprop\label{p:eldt}
For each $\L\in \C(T)(\Qlcl)$, the subcategory (see Definition \ref{d:qi}) ${}^{e_\L}\D(T)\subset \D(T)$ is the full subcategory formed by objects of the form $W\otimes e_\L$ for some $W\in D^b\Vec$. In other words, we have a triangulated functor $D^b\Vec\xto{}{}^{e_\L}\D(T), W\mapsto W\otimes e_\L$ which induces a bijection on the isomorphism classes of objects. Similarly, the full subcategory ${}^{e_\L}\D_T(T)\subset \D_T(T)$ is the full subcategory formed by objects of the form $W\otimes e_\L$ for some $W\in \D_T(1)$. 
\eprop 
\bpf
Let $M\in {}^{e_\L}\D(T)$. Hence $e_\L\ast M\cong V_T\otimes M.$ On the other hand, by (\ref{e:elm}) we have $e_\L\ast M\cong (\L\ast M)_1\otimes e_\L$ and $e_\L$ is indecomposable. The first part of the proposition now follows.

Now recall that we have an action of $\D_T(1)\subset \D_T(T)$ on $\D_T(T)$ which we will also denote by $\otimes$. If $M\in {}^{e_\L}\D_T(T)$, then as before we obtain that $e_\L\ast M\cong V_T\otimes M\cong (\L\ast M)_1\otimes e_\L$, where $(\L\ast M)_1\in \D_T(1)$. From this we see that $M\cong W\otimes e_\L$ for some $W\in \D_T(1)$.
\epf

\brk\label{r:tstr}
For a minimal quasi-idempotent $e_\L\in \D(T)$ (resp. $\D_T(T)$), let ${}^{e_\L}\D(T)^\Delta$ (resp. ${}^{e_\L}\D_T(T)^\Delta$) be the full triangulated subcategories of $\D(T)$ (resp. $\D_T(T)$) generated by ${}^{e_\L}\D(T)$ (resp. ${}^{e_\L}\D_T(T)$) (see Defn. \ref{d:qi}). Then we see that the categories ${}^{e_\L}\D(T)^\Delta$ and ${}^{e_\L}\D_T(T)^\Delta$ have unique $t$-structures such that $e_\L$ lies in their heart. The semisimple abelian subcategories of the hearts formed by the semisimple objects are both equivalent to the category $\Vec$ with $e_\L$ being the unique simple object. Moreover, it follows from Lemma \ref{l:sasha} that any non-degenerate bounded $t$-structure on either of these triangulated categories must necessarily be a shift of the above $t$-structures.
\erk

\section{Character sheaves in the Heisenberg case}\label{s:thc}
In this section we work with neutrally solvable groups $G$ and analyze the case of Heisenberg idempotents $e\in \DG$. We will define character sheaves (cf. Definition \ref{d:csfg}) and their $\fL$-packet decomposition inside the $\Qlcl$-linear triangulated braided monoidal category $e\DG$. We will state and prove a version of our main results (Theorems \ref{t:main1},\ref{t:main2} and \ref{t:main3}) for Heisenberg idempotents in this section and \S\ref{s:icahi}. The results from \cite{De2} would then allow us to prove (in \S\ref{s:pmr}) our main results for general minimal idempotents, since by {\it op. cit.} every minimal idempotent on a neutrally solvable group can essentially be obtained from a Heisenberg idempotent on a subgroup of $G$ by using the induction with compact supports functor. 

We begin by describing the setting of Heisenberg idempotents. Let $G$ be a neutrally solvable (perfect quasi-) algebraic group over $\k$. Let $G^\circ=UT$ where $U$ is the unipotent radical and $T$ is some maximal torus. Let $(H,\N)$ be a Heisenberg admissible (cf. \cite[\S2.6]{De2}) pair for $G$. Recall that this means that $H$ is a connected normal unipotent subgroup of $G$, $\N\in (H^*)^G$ a $G$-equivariant multiplicative local system on $H$, $U/H$ is commutative and the induced skewsymmetric biextension 
\beq\phi_\N:U/H\rar{}(U/H)^*\eeq is an isogeny. Associated with such a skewsymmetric isogeny we have the metric group $(K_\N,\theta)$ where $K_\N=\ker(\phi_\N)$. Let $e:=\N\otimes\can_H\in\DG$ be the Heisenberg idempotent associated with the admissible pair. Then from \cite{De2} we know that $e\in \DG$ is a minimal idempotent and that it is closed. We also know that the full subcategory $e\D_U(U)\subset e\D_U(G)$ is equivalent to the bounded derived category of the modular category $\M(K_\N,\theta)$ associated with the metric group $(K_\N,\theta)$. We refer to \cite{De1} for details.

\brk\label{r:fundim}
For a Heisenberg idempotent $e$ as above, we set $d_e:=\frac{\dim U-\dim H}{2}$ and call it the functional dimension of $e$.
\erk

Let $\Gamma:=G/U$. We have $\Gamma^\circ=T$. Let $\Pi_0:=\pi_0(G)=G/G^\circ = \pi_0(\Gamma)=\Gamma/T$. Note that conjugation by $\Gamma$ induces an action of $\Pi_0$ on $T$. We want to define character sheaves in the $\Qlcl$-linear triangulated braided monoidal category $e\DG$. We will first study the triangulated braided $\Pi_0$-crossed category 
\beq
e\D_{G^\circ}(G)=\bigoplus\limits_{G^\circ g\in \Pi_0}e\D_{G^\circ}(G^\circ g).
\eeq
Now by \cite[Prop. 8.10]{De2}, we have an equivalence $e\D_{G^\circ}(G^\circ)\cong e\D_U(U)\boxtimes \D_T(T)$ of braided triangulated categories and for each $\bar{g}=G^\circ g\in \Pi_0$ a triangulated equivalence
\beq
e\D_{G^\circ}(G^\circ g)\cong e\D_U(Ug)\boxtimes \D^{\bar{g}}_T(T)
\eeq
of $e\D_{G^\circ}(G^\circ)\cong e\D_U(U)\boxtimes \D_T(T)$-module categories, where $\D^{\bar{g}}_T(T)$ is the $T$-equivariant derived category of $T$ for the $\bar{g}$-conjugation action of $T$ on itself (see also \cite[\S8.5]{De2}).  We can obtain the braided monoidal category $e\DG$ as the $\Pi_0$-equivariantization of $e\D_{G^\circ}(G)$.

Finally, let us recall from \cite[\S2.3]{De1} that each of the monoidal categories $e\D_U(G), e\D_{G^\circ}(G)$ and $e\DG$ is a monoidal $\r$-category. The duality functor is defined as $M\mapsto M^\vee:=\f{D}^{-}M\otimes \can_H\cong \f{D}^-M[2\dim H](\dim H)$.

\subsection{The connected Heisenberg case}\label{s:chc}
In this section we suppose that $G=TU$ is a connected solvable group over $\k$ and that $(H,\N)$ is a Heisenberg admissible pair for $G$ with $e=\N\otimes\can_H$ being the corresponding Heisenberg idempotent. Let us begin by classifying the minimal quasi-idempotents in the triangulated monoidal category $e\D_U(G)$. By \cite[Thm. 2.26]{De2} we have a canonical monoidal equivalence $e\D_U(G)\cong e\D_U(U)\boxtimes \D(T)$. In particular  $e\D_U(G)$ is braided and we can talk of minimal quasi-idempotents in this category.
\bprop\label{p:mqics}
(i) Let $\L\in \D(T)$ be a multiplicative local system and let $e_\L=\L\otimes\can_T\in \D(T)$ be the corresponding  minimal quasi-idempotent. Then $e\boxtimes e_\L\in e\D_U(U)\boxtimes \D(T)\cong e\D_U(G)$ is a minimal quasi-idempotent. \\
(ii) The full subcategory (see Definition \ref{d:qi}) ${}^{(e\boxtimes e_\L)}e\D_U(U)\boxtimes \D(T)\subset e\D_U(U)\boxtimes \D(T)\cong e\D_U(G)$ associated with the quasi-idempotent $e\boxtimes e_\L$ is equal to $e\D_U(U)\boxtimes {}^{e_\L}\D(T)$. It is the full subcategory of objects of the form  $M\boxtimes e_\L$ for $M\in e\D_U(U)$.\\
(iii) The mapping $\L\mapsto e\boxtimes e_\L$ establishes a bijection between isomorphism classes of multiplicative local systems on $T$ and equivalence classes of minimal quasi-idempotents in $e\D_U(U)\boxtimes \D(T)$. Each equivalence class of minimal quasi-idempotents in $e\D_U(U)\boxtimes \D(T)$ contains a unique (up to isomorphism) indecomposable object, namely some $e\boxtimes e_\L$ for some $\L\in \C(T)(\Qlcl)$. Any minimal quasi-idempotent in $e\D_U(U)\boxtimes \D(T)$ is of the form $W\otimes (e\boxtimes e_\L)$ for some $\L\in \C(T)(\Qlcl)$ and some $W\in \Dmon$.
\eprop
\bpf
It is clear that $e\boxtimes e_\L$ is a  $V_T$-quasi-idempotent. Let us first prove (ii). Objects of $e\D_U(U)\boxtimes \D(T)$ are direct sums of objects of the form $A\boxtimes B$ with $A\in e\D_U(U)$ and $B\in \D(T)$. Note that by (\ref{e:elm}) we have  
\beq
(e\boxtimes e_\L)\ast (A\boxtimes B)\cong (\L\ast B)_1\otimes A\boxtimes e_\L.
\eeq
Hence for any $X\in e\D_U(U)\boxtimes \D(T)$, the convolution $(e\boxtimes e_\L)\ast X$ is an object of the form $Y\boxtimes e_\L$ for some $Y\in e\D_U(U)$.   Hence we see that the full subcategory ${}^{(e\boxtimes e_\L)}e\D_U(U)\boxtimes \D(T)\subset e\D_U(U)\boxtimes \D(T)$ associated with the quasi-idempotent $e\boxtimes e_\L$ is equal to $e\D_U(U)\boxtimes {}^{e_\L}\D(T)$. This is the full subcategory  of objects of the form $Y\boxtimes e_\L$ where $Y\in e\D_U(U)$. 

Now suppose that such an object is a nonzero quasi-idempotent, i.e. $(Y\boxtimes e_\L)\ast(Y\boxtimes e_\L)\cong V\otimes (Y\boxtimes e_\L)$ for some $V\in \Dmon$. Hence $(Y\ast Y)\boxtimes (V_T\otimes e_\L)\cong V\otimes Y\boxtimes e_\L$, which implies that $V_T\otimes Y\ast Y\cong V\otimes Y$. This means that $Y\in e\D_U(U)$ must be a quasi-idempotent. Now $e\D_U(U)$ is the bounded derived category of a modular category. Hence by Proposition \ref{p:qidbm} $Y$ must be of the form $W\otimes e$ for some $W\in \Dmon$ and hence $Y\boxtimes e_\L\cong W\otimes e\boxtimes e_\L$. Hence by Corollary \ref{c:altminqi} $e\boxtimes e_\L$ is a minimal quasi-idempotent in $e\D_U(G)\cong e\D_U(U)\boxtimes \D(T)$.

Now let $X\in e\D_U(U)\boxtimes \D(T)$ be any minimal quasi-idempotent. Using  Lemma \ref{l:mellin}(ii) and the fact that each $X\in e\D_U(U)\boxtimes \D(T)$ can be expressed as a direct sum of objects of the form $A\boxtimes B$, we see that there exists a multiplicative local system $\L$ on $T$ such that $(e\boxtimes e_\L)\ast X\neq 0$. Since both $X$ and $e\boxtimes e_\L$ are minimal quasi-idempotents, they must be equivalent by Definition \ref{d:minqi}. Moreover, since $e\boxtimes e_\L$ is an indecomposable, we must have $X\cong W\otimes (e\boxtimes e_\L)$ for some $W\in \Dmon$ (cf. Cor. \ref{c:allminqi} and its proof). 
\epf

\brk\label{r:monfun}
In the future we will denote the canonical equivalence $e\D_U(U)\boxtimes \D(T)\rar{\cong}e\D_U(G)$ by $X\mapsto \bar{X}$. Note that this equivalence also defines a canonical triangulated monoidal functor $\D(T)\rar{}e\D_U(G)$ (cf. \cite[\S8.5]{De2}) which we will also denote by $N\mapsto \bar{N}$. In particular, image of the minimal quasi-idempotent $e\boxtimes e_\L$ in $e\D_U(G)$ will be denoted by $\bar{e_\L}$. An object $M\boxtimes N\in e\D_U(U)\boxtimes \D(T)$ corresponds to the object $M\ast \bar{N} \in e\D_U(G)$.
\erk

\subsection{The possibly disconnected case}
Let $G$ be a neutrally solvable group and let $(H,\N)$ be a Heisenberg admissible pair for $G$ and let $e\in \DG$ be the corresponding Heisenberg idempotent. In the next few sections, we will classify all the minimal quasi-idempotents in the $\Qlcl$-linear  triangulated braided category $e\DG$. In this section, we study some of these minimal quasi-idempotents. 

We continue to use the notation from the introduction of \S\ref{s:thc}. In particular $U$ is the unipotent radical of $G$ and $G^\circ=TU$ where $T$ is some maximal torus. We denote the surjection $G^\circ \onto T$ by $h\mapsto \bar h$. We have an action of $\Pi_0$ on $T$ coming from the conjugation action of $\Gamma=G/U$ on $T$ and this gives us an action of $G$ on $T$. We denote this action by $g:t\mapsto{}^gt$. 

\subsubsection{The category $e\D_U(G)$}
We will first work with the category $e\D_U(G)$. We have $e\D_U(G)=\bigoplus\limits_{G^\circ g\in \Pi_0}e\D_U(G^\circ g)$. By \cite[Defn. 8.7, Prop. 8.8]{De2} for each $g\in G$ we have an  equivalence 
\beq
e\D_U(Ug)\boxtimes \D(T)\cong e\D_U(G^\circ g)\mbox{ defined by } M\boxtimes N\mapsto M\ast \bar{N},
\eeq
where $N\mapsto \bar{N}$ is the canonical monoidal functor $\D(T)\rar{}e\D_U(G^\circ)$ as defined in {\it loc. cit.} (see also Remark \ref{r:monfun}). Moreover, by \cite[Thm. 6.2]{De2} $e\D_U(Ug)\cong D^b(\tMg)$, where $\tMg$ is an invertible module category over the pointed modular category $\tM_{U,e}\cong \M(K_\N,\theta)$.

Note that for each $h\in G^\circ$ we have an equivalence of $e\D_U(U)$-module categories
\beq
e\D_U(Ug)\xto{(\cdot)\ast {\bar\delta_{{}^{g^{-1}}{\bar{h}}}}}e\D_U(Uhg).
\eeq
Now for each $a\in \Pi_0$, let us choose a lift $\ta\in G$ such that $\t{1}=1$. For $a,b\in \Pi_0$, let
\beq
\ta\tb=f(a,b)\t{ab}, \mbox{ where } f(a,b)\in G^\circ.
\eeq Then by \cite[Prop. 8.9]{De2}, the convolution $e\D_U(G^\circ\ta)\times e\D_U(G^\circ\tb)\rar{\ast}e\D_U(G^\circ\tilde{ab})$ (for $a,b\in \Pi_0$) can be identified with the composition
\beq
e\D_U(U\ta)\boxtimes \D(T)\times e\D_U(U\tb)\boxtimes \D(T)\rar{}e\D_U(U\ta\tb)\boxtimes \D(T)\rar{\cong}e\D_U(U\t{ab})\boxtimes \D(T)
\eeq defined by
\beq\label{e:convinedug}
(M_1\boxtimes N_1)\ast(M_2\boxtimes N_2)=(M_1\ast M_2\ast \bar{\delta}_{{}^{b^{-1}a^{-1}}(\bar{f(a,b)^{-1}})})\boxtimes (\delta_{{}^{b^{-1}a^{-1}}(\bar{f(a,b)})}\ast b^{-1}(N_1)\ast N_2).
\eeq

Let us define an auxiliary triangulated monoidal category $(\D,\ast')$ as follows:
\beq
\D=\bigoplus\limits_{a\in \Pi_0}e\D_U(U\ta)
\eeq
and we define $M_1\ast' M_2:=M_1\ast M_2\ast \bar\delta_{{}^{b^{-1}a^{-1}}(\bar{f(a,b)^{-1}})}$ for $M_1\in e\D_U(U\ta), M_2\in e\D_U(U\tb)$ and where the convolution on the right hand side comes from the category $e\D_U(G)$. Then $\D$ can be equipped with an associativity constraint using the isomorphisms $\bar\delta_s\ast \bar\delta_t\cong \bar\delta_{st}$ in $e\D_U(G)$ for $s,t\in T$. This associativity constraint is not quite canonical, however this does not matter in the proof of the next lemma. Moreover note that by \cite{De1} each $e\D_U(U\ta)$ is the bounded derived category of a finite semisimple abelian category and $\D$ is equivalent to the bounded derived category of a fusion category.

\blem\label{l:minqiedug}
Let $\L$ be a $\Pi_0$-equivariant multiplicative local system on $T$ and $e_\L\in \D(T)$ the corresponding minimal quasi-idempotent. Then $e\boxtimes e_\L\in e\D_U(U)\boxtimes \D(T)\cong e\D_U(G^\circ)\subset e\D_U(G)$ is a minimal weakly central quasi-idempotent in $e\D_U(G)$. (The image of $e\boxtimes e_\L$ in $e\D_U(G^\circ)\subset e\D_U(G)$ is $\beel$.) We have 
\beq
{}^{\beel}e\D_U(G)=\bigoplus\limits_{G^\circ g\in \Pi_0}{}^{\beel}e\D_U(G^\circ g)\cong\bigoplus\limits_{G^\circ g\in \Pi_0} e\D_U(Ug)\boxtimes {}^{e_\L}\D(T).
\eeq
In other words, the full subcategory ${}^{\beel}e\D_U(G^\circ g)\subset e\D_U(G^\circ g)\cong e\D_U(Ug)\boxtimes \D(T)$ consists of objects of the form $\bar{M\boxtimes e_\L}$ with $M\in e\D_U(Ug)$. Hence we also obtain (see Definition \ref{d:qi} for notation)
\beq
{}^{\beel}e\D_U(G^\circ g)^\Delta=\bigoplus\limits_{G^\circ g\in \Pi_0}{}^{\beel}e\D_U(G^\circ g)^\Delta\cong\bigoplus\limits_{G^\circ g\in \Pi_0} e\D_U(Ug)\boxtimes {}^{e_\L}\D(T)^\Delta.
\eeq
\elem
\bpf
Using the $\Pi_0$-equivariance of $\L$ and (\ref{e:convinedug}), it is easy to check that $\beel$ is a weakly  central quasi-idempotent in $e\D_U(G)$. Let $A\boxtimes B\in e\D_U(Ug)\boxtimes \D(T)$. Then we have $(e\boxtimes e_\L)\ast(A\boxtimes B)\cong A\boxtimes (g^{-1}(e_\L)\ast B)$. Since $\L$ is $\Pi_0$-equivariant, we obtain $(e\boxtimes e_\L)\ast(A\boxtimes B)\cong A\boxtimes (e_\L\ast B)\cong (\L\ast B)_1\otimes A\boxtimes e_\L$ using (\ref{e:elm}). Then as in the proof of Proposition \ref{p:mqics}, we conclude that each object in ${}^{(e\boxtimes e_\L)}e\D_U(Ug)\boxtimes \D(T)$ is of the form ${M\boxtimes e_\L}$. Now by Proposition \ref{p:qidbm} the nonzero quasi-idempotents in the category $(\D,\ast')$ defined above are all of the form $V\otimes e$ for some $V\in \Dmon$. Hence we see that all the nonzero quasi-idempotents in ${}^{\beel}e\D_U(G)$ are of the form $V\otimes \beel$ for some $V\in \Dmon$. Hence by Corollary \ref{c:altminqi} $\beel$ is a minimal weakly central quasi-idempotent in $e\D_U(G)$.
\epf

\subsubsection{The category $e\D_{G^\circ}(G)$}
Let us now study the category $e\D_{G^\circ}(G)$ which is a braided $\Pi_0$-crossed category:
\beq
e\D_{G^\circ}(G)=\bigoplus\limits_{G^\circ g\in \Pi_0}e\D_{G^\circ}(G^\circ g).
\eeq
By \cite[Prop. 8.10]{De2}, we have identifications
\beq
e\D_{G^\circ}(G^\circ g)\cong e\D_U(Ug)\boxtimes \D^g_T(T),
\eeq
where $\D^g_T(T)$ denotes the category $T$-equivariant $\Qlcl$-complexes for the $g$-twisted conjugation action of $T$ on itself. In terms of this identification, for $g_1,g_2\in G$, the convolution (cf. \cite[Prop. 8.9]{De2})
\beq
e\D_{G^\circ}(G^\circ g_1)\times e\D_{G^\circ}(G^\circ g_2)\rar{\ast}e\D_{G^\circ}(G^\circ g_1g_2) \mbox{ corresponds to the convolution }
\eeq
\beq
e\D_{U}(Ug_1)\boxtimes \D^{g_1}_T(T)\times e\D_{U}(Ug_2)\boxtimes \D^{g_2}_T(T)\rar{\ast}e\D_{U}(Ug_1g_2)\boxtimes \D^{g_1g_2}_T(T) \mbox{ defined by}
\eeq
\beq\label{e:convinedgg}
(M_1\boxtimes N_1)\ast (M_2\boxtimes N_2):=(M_1\ast M_2)\boxtimes ({}^{g_2^{-1}}(N_1)\ast N_2).
\eeq

Using the results of the previous section, we will prove
\bcor\label{c:beeledgg}
In the setup of Lemma \ref{l:minqiedug}, we may consider $\beel$ as an object of $e\D_{G^\circ}(G^\circ)\subset e\D_{G^\circ}(G)$. Then $\beel$ is a minimal weakly central quasi-idempotent in the braided $\Pi_0$-crossed category $e\D_{G^\circ}(G)$. We have
\beq\label{e:beeledgg}
{}^{\beel}e\D_{G^\circ}(G)=\bigoplus\limits_{G^\circ g\in \Pi_0}{}^{\beel}e\D_{G^\circ}(G^\circ g)\cong\bigoplus\limits_{G^\circ g\in \Pi_0} e\D_U(Ug)\boxtimes {}^{e_\L}\D^g_T(T).
\eeq
All objects of ${}^{\beel}e\D_{G^\circ}(G^\circ g)\cong e\D_U(Ug)\boxtimes {}^{e_\L}\D^g_T(T)$ are of the form $V\otimes (M\boxtimes e_\L)$ for some $M\in e\D_U(Ug)$ and $V\in \D_T(1)$. We also obtain (see Definition \ref{d:qi} for notation)
\beq\label{e:beeledggdelta}
{}^{\beel}e\D_{G^\circ}(G)^\Delta=\bigoplus\limits_{G^\circ g\in \Pi_0}{}^{\beel}e\D_{G^\circ}(G^\circ g)^\Delta\cong\bigoplus\limits_{G^\circ g\in \Pi_0} e\D_U(Ug)\boxtimes {}^{e_\L}\D^g_T(T)^\Delta.
\eeq
\ecor
\bpf
Note that we have a canonical monoidal functor $\D_T(T)\rar{}e\D_{G^\circ}(G^\circ)\subset e\D_{G^\circ}(G)$, which we also denote by $X\mapsto \bar{X}$. Then we see that $\beel$ is a $V_T$-quasi-idempotent in $e\D_{G^\circ}(G)$. We see that it is weakly central from (\ref{e:convinedgg}). Furthermore, we see that $\beel$ is a minimal weakly central quasi-idempotent in $e\D_{G^\circ}(G)$ using Lemma \ref{l:minqiedug} and the forgetful functor $e\D_{G^\circ}(G)\rar{}e\D_U(G)$. The proof of (\ref{e:beeledgg}) is similar to the proof of Lemma \ref{l:minqiedug} using (\ref{e:convinedgg}). To prove the next statement, we only need to further observe that all objects of ${}^{e_\L}\D_T^g(T)$ are of the form $V\otimes e_\L$ for some $V\in \D_T(\pt)$. First note that $e_\L$ is perverse up to a shift and is $\Pi_0$-equivariant. Using the isomorphism between the multiplicative local systems $\L$ and ${}^g\L$ on $T$, we obtain a canonical equivariant structure on $e_\L$ for the $g$-twisted action of $T$ on itself. Thus we see that $e_\L\in \D^g_T(T)$. Now using a similar argument as in the proof of Proposition \ref{p:eldt} we conclude that all objects of ${}^{e_\L}\D^g_T(T)$ are indeed of the form $V\otimes e_\L$ for some $V\in \D_T(1)$. Statement (\ref{e:beeledggdelta}) about the triangulated subcategories also follows.
\epf

\brk\label{r:edggdgg}
In the setting above, $\beel$ may also be considered as a minimal weakly central quasi-idempotent in $\D_{G^\circ}(G)$ and we see that ${}^{\beel}e\D_{G^\circ}(G)$ and ${}^{\beel}\D_{G^\circ}(G)$ are equal as full subcategories of $\D_{G^\circ}(G)$ and we may use either notation to denote it.
\erk

%\brk\label{r:spec}
%We can identify $\D_T(1)\subset \D_T(T)$ with the full triangulated subcategory of $e\D_{G^\circ}(G^\circ)\cong e\D_U(U)\boxtimes \D_T(T)\subset e\DGc$ generated by the unit object $e$. Hence we have an action of $\D_T(1)$ on $e\DGc$ which we denote by $\otimes$. By Corollary \ref{c:beeledgg}, the additive category ${}^{\beel}e\D_{G^\circ}(G)$ is closed under the $\otimes$-action of $\D_T(1)$ and splits as a direct sum of the additive categories closed under tensoring by $\D_T(1)$  generated by $\bar{M\boxtimes e_\L}$ where $M$ ranges over the simple objects of the heart $\tM_{Ug,e}\subset e\D_U(Ug)$ as $g$ ranges over representatives of the connected components of $G$. The objects $\bar{M\boxtimes e_\L}\in {}^{\beel} e\DGc$ and their shifts are characterized by this property. We call these objects the special indecomposable objects of ${}^{\beel}e\DGc$. Notice that these special indecomposable objects have been defined purely in terms of the triangulated monoidal structure of $e\DGc$ and the minimal idempotent $\beel$.
%\erk

\subsubsection{Truncated convolution and fusion categories}
Let $\L$ be a $\Pi_0$-equivariant multiplicative local system on $T$ and $\beel\in e\D_{G^\circ}(G)$ the associated minimal weakly central quasi-idempotent. We will now associate with $\beel$ a braided $\Pi_0$-crossed fusion category which we will denote by $\tMbeel$. 
\bprop\label{p:beeltrsubcat}
(i) Let $\L\in \C(T)(\Qlcl)^{\Pi_0}$. Then ${}^{\beel}e\D_{G^\circ}(G)^\Delta$ is a $\Qlcl$-linear triangulated braided $\Pi_0$-crossed semigroupal category. It has a unique non-degenerate bounded t-structure, denoted here by $(\D^{\leq 0}, \D^{\geq 0})$ with heart denoted by $\tMbeel^\Delta$, such that $\D^{\leq 0} \ast \D^{\leq 0} \subset \D^{\leq 0}$ and $\D^{\leq 0} \ast \D^{\leq 0}\not\subset \D^{\leq -1}$.\\
%It has a unique $t$-structure, denoted here by $(\D^{\leq 0}, \D^{\geq 0})$, such that $\D^{\leq 0}\ast \D^{\leq 0}\subset \D^{\leq 0}$ but $\D^{\leq 0}\ast \D^{\leq 0}\nsubseteq \D^{\leq -1}$. \\
(ii) For each $i\in \Z, X\in {}^{\beel}e\D_{G^\circ}(G)^\Delta$ let $X^i\in \tMbeel^\Delta$ denote the $i$-th cohomology of $X$ with respect to this $t$-structure. The truncated convolution
\beq
\underline{\ast}:\tM_{G,\beel}^\Delta\times \tM_{G,\beel}^\Delta\rar{}\tM_{G,\beel}^\Delta \hbox{, defined as } X\uast Y:=(X\ast Y)^0 \h{ for } X,Y\in \tMbeel^\Delta
\eeq 
equips $\tMbeel^\Delta$ with the structure of a braided $\Pi_0$-crossed semigroupal category.\\
(iii) Let $\tMbeel\subset \tMbeel^\Delta$ be the full subcategory formed by the semisimple objects. Then $\tMbeel$ is closed under truncated convolution,  $\beel$ is a unit object for truncated convolution on $\tMbeel$ and $(\tMbeel,\uast, \beel)$ is a braided $\Pi_0$-crossed fusion category. The simple objects of $\tMbeel$ are of the form $\bar{M\boxtimes e_\L}$ with $M$ a simple object in $\tMg$ for some $g\in G$ and $e_\L$ is considered as an object of $\D_T^g(T)$.\\
%direct sums of the indecomposable objects $\bar{M\boxtimes e_\L}$ (with simple $M\in \tMg$) described above.\\
(iv) The identity component of $\tMbeel$ is equivalent to the pointed modular category corresponding to the metric group $(K_\N,\theta)$ (see the introduction of \S\ref{s:thc}) and the Frobenius-Perron dimension of $\tMbeel$ (denoted by $\FPdim(\tMbeel)$) is equal to $|\Pi_0|\cdot |K_\N|\in\Z$.
\eprop
\bpf
To prove (i), let us first construct one $t$-structure satisfying the desired properties. We have described the categories ${}^{\beel}e\D_{G^\circ}(G)\subset {}^{\beel}e\D_{G^\circ}(G)^\Delta$ in Corollary \ref{c:beeledgg}. Now by \cite[Thm. 6.2]{De2}, we have $e\D_U(Ug)\cong D^b(\tMg)$ where $\tMg$ is an invertible $\tM_{U,e}$-module category. For $g\in G$, consider the  $t$-structure on $e\D_{U}(Ug)$ whose heart is $\tMg$. Also let us consider the $t$-structure on the category ${}^{e_\L}\D_T^g(T)^\Delta$  such that the object $e_\L\in {}^{e_\L}\D^g_T(T)\subset {}^{e_\L}\D^g_T(T)^\Delta$ lies in its heart (which we denote by $\tM_{T,g,e_\L}^\Delta$). In this way, we obtain a $t$-structure on $e\D_U(Ug)\boxtimes {}^{e_\L}\D^g_T(T)^\Delta$. We transfer this $t$-structure to ${}^{\beel}e\D_{G^\circ}(G^\circ g)^\Delta$ using Corollary \ref{c:beeledgg}. Let us denote the heart of this $t$-structure by $\tM_{G^\circ g, \beel}^\Delta$. Thus we obtain a $t$-structure on all of ${}^{\beel}e\D_{G^\circ}(G)^\Delta$. We obtain a grading on the heart of this $t$-structure
\beq\label{e:Delgrading}
\tMbeel^\Delta = \bigoplus\limits_{G^\circ g\in \Pi_0} {\tM_{G^\circ g,\beel}^\Delta}\cong\bigoplus\limits_{G^\circ g\in \Pi_0} \tMg\boxtimes \tM_{T,g,e_\L}^\Delta.
\eeq %and that $\tMbeel$ is a semisimple abelian category with finitely many simple objects. 
Using \cite[Thm. 6.2(iii)]{De2}, Corollary \ref{c:beeledgg} and (\ref{e:convinedgg}) we see that this $t$-structure on ${}^{\beel}e\D_{G^\circ}(G)$ satisfies the properties desired in (i). 

Let us now prove the uniqueness of the $t$-structure from statement (i). By Appendix \ref{s:sasha}, the hearts of all the possible different non-degenerate bounded $t$-structures on ${}^{\beel}e\D_{G^\circ}(G)^\Delta$ must have as their simple objects the various shifts of the simple objects of $\tMbeel^\Delta$. Hence we see that there is a unique such $t$-structure which also satisfies the desired properties from (i).

To prove (ii), using the same argument as in \cite[Thm. 3.1.3.1]{Li}, we can construct canonical associativity as well as crossed braiding isomorphisms for the truncated convolution by truncating the corresponding structures from $e\D_{G^\circ}(G)^\Delta$.

We now prove (iii) and (iv). It is clear from (\ref{e:convinedgg}) that $\tMbeel$ is closed under truncated convolution and hence is a braided $\Pi_0$-crossed semigroupal category. From (\ref{e:Delgrading}) we obtain the grading
\beq\label{e:tmbeelgrading}
\tMbeel = \bigoplus\limits_{G^\circ g\in \Pi_0} {\tM_{G^\circ g,\beel}}\cong\bigoplus\limits_{G^\circ g\in \Pi_0} \tMg\boxtimes \tM_{T,g,e_\L},
\eeq 
where $\Vec\cong\tM_{T,g,e_\L}\subset \tM_{T,g,e_\L}^\Delta$ is the full subcategory formed by (finite) direct sums of the simple object $e_\L$.
We will now prove that $(\tMbeel,\uast,\beel)$ is a fusion category. First note that we have a canonical morphism $\delta_1\rar{a}e_\L$ in $\D_T(T)$. The construction of this morphism is exactly the same as in the unipotent group case as described in \cite[\S8.1.2]{B1}. Moreover, using a similar argument to the one in \cite[\S8.3]{B1} we can prove that if we convolve this morphism with $e_\L$, we get a morphism $e_\L\rar{}e_\L\ast e_\L\cong V_T\otimes e_\L$ which is an isomorphism onto the direct summand $e_\L$ of $V_T\otimes e_\L$. Using this, we obtain a similar morphism $e\rar{a}\beel$ in $e\D_{G^\circ}(G)$. Now each simple object of $\tMg\boxtimes \tM_{T,g,e_\L}\cong \tM_{G^\circ g,\beel}$ is of the form $M\boxtimes e_\L$ for some simple $M\in \tMg$. Let $X\in\tM_{G^\circ g, \beel}$ be any object. Convolving the morphism $a$ with $X$, we obtain a morphism $X\rar{}\beel\ast X$. Then using the previous observations, we see that the induced morphism $X=X^0\rar{}(\beel\ast X)^0=\beel\uast X$ obtained by taking the $0$-th cohomology with respect to our chosen $t$-structure is a (functorial) isomorphism. This makes $\beel$ into a left unit in $\tMbeel$. Similarly, it is also a right unit. Thus we see that $(\tMbeel,\uast,\beel)$ has the canonical structure of a braided $\Pi_0$-crossed monoidal category with identity component $\tM_{G^\circ,\beel}\cong \tM_{U,e}$ (by (\ref{e:tmbeelgrading})), a pointed modular category which corresponds to the metric group $(K_\N,\theta)$ (cf. \cite{De1}). Now by \cite[\S6.1]{De1} we conclude that $\tMbeel$ is rigid and hence a fusion category as desired. The statement about the Frobenius-Perron dimension follows from the fact that $\FPdim(\M_{U,e})=|K_\N|$ and \cite{DGNO}.
\epf

\brk\label{r:spherical}
Under the identification $e\D_U(Ug)\boxtimes \D(T)\cong e\D_U(G^\circ g)$, the contravariant (weak) duality functor $(\cdot)^\vee=\f{D}^-(\cdot)[2\dim H](\dim H):e\D_U(G^\circ g)\rar{}e\D_U(G^\circ g^{-1})$ corresponds to 
$$M\boxtimes N\mapsto M^\vee\boxtimes \f{D}^-g(N)=\f{D}^-M[2\dim H](\dim H)\boxtimes \f{D}^-g(N).$$ Now consider the identification $\tM_{G^\circ g,\beel}\cong \tM_{Ug,e}\boxtimes \tM_{T,g,e_\L}$. Note that the (rigid) duality in $\tM_{T,1,e_\L}(\cong \Vec)$ is given by $\f{D}^-(\cdot)[2\dim T](\dim T)$. Then we see that the duality functor $(\cdot)^*:\tM_{G^\circ g,\beel}\rar{}\tM_{G^\circ g^{-1},\beel}$ corresponds to $M\boxtimes N\mapsto \f{D}^-M[2\dim H](\dim H)\boxtimes \f{D}^-g(N)[2\dim T](\dim T)$. Hence we see that the (rigid) duality in $\tM_{G,\beel}$ is given by $\f{D}^-(\cdot)[2\dim H+2\dim T](\dim H+\dim T)$. We see that we have a natural identification of the square of the duality functor with the identity functor, or in other words $\tM_{G,\beel}$ has a natural spherical structure.
\erk

\subsection{Minimal quasi-idempotents in $e\DG$}\label{s:minqiedgg}
In this section we consider the category $e\DG$ and classify all minimal quasi-idempotents in it. We first obtain the following corollary of our previous results:
\bcor\label{c:minqiedgg}
(i) Let $\L\in \C(T)(\Qlcl)^{\Pi_0}$ be as before. Then we can consider $\beel$ as an object of $e\D_G(G^\circ)\subset e\DG$. Then $\beel$ is a minimal quasi-idempotent in the braided triangulated category $e\DG$ as well as the category $\DG$. \\
(ii) Furthermore, ${}^{\beel}e\D_{G}(G)^{\Delta}$ is a $\Qlcl$-linear triangulated braided semigroupal category. It has a unique $t$-structure, denoted here by $(\D^{\leq 0}, \D^{\geq 0})$, such that $\D^{\leq 0}\ast \D^{\leq 0}\subset \D^{\leq 0}$ but $\D^{\leq 0}\ast \D^{\leq 0}\nsubseteq \D^{\leq -1}$. \\
(iii) Let $\Mbeel^\Delta\subset {}^{\beel}\D_{G}(G)^\Delta$ denote the heart of the above $t$-structure. For each $i\in \Z, X\in {}^{\beel}\D_{G}(G)^\Delta$ let $X^i\in \Mbeel^\Delta$ denote the $i$-th cohomology of $X$ with respect to this $t$-structure. Then $(\M_{G,\beel}^\Delta,\underline{\ast})$ has the structure of a  braided semigroupal category, where 
\beq
\underline{\ast}:\M_{G,\beel}^\Delta\times \M_{G,\beel}^\Delta\rar{}\M_{G,\beel}^\Delta
\eeq denotes the truncated convolution functor defined by $X\uast Y:=(X\ast Y)^0$ for $X,Y\in \Mbeel^\Delta$. \\
(iv) Let $\Mbeel\subset \Mbeel^\Delta$ denote the full subcategory formed by the semisimple objects. Then $\Mbeel$ is closed under truncated convolution and $(\Mbeel,\uast,\beel)$ is a non-degenerate braided fusion category. We have an identification $\Mbeel=(\tMbeel)^{\Pi_0}$ and $\FPdim(\M_{G,\beel})=|\Pi_0|^2\cdot |K_\N|.$
\ecor
\bpf
Note that from the canonical monoidal functor $\D(T)\rar{}e\D_U(G^\circ)$, we also obtain a monoidal functor $\D_\Gamma(T)\rar{}e\D_G(G^\circ)\subset e\DG$ which we also denote by $X\mapsto \bar X$. Now we can consider $e_\L$ as an object of $\D_T(T)$ and since $\L$ is $\Gamma$-equivariant, we may consider $e_\L$ as an object of $\D_\Gamma(T)$. Thus we can consider the object $\beel$ as an object of $e\D_G(G)$. Then we see that $\beel$ is a minimal quasi-idempotent by using Lemma \ref{l:minqiedug} and the forgetful functor $e\D_G(G)\rar{}e\D_U(G)$. Statements (ii), (iii) and (iv) also follow readily from Proposition \ref{p:beeltrsubcat} after $\Pi_0$-equivariantization. Note that the non-degeneracy of the braided fusion category $\Mbeel$ follows from \cite[Prop. 4.56]{DGNO}.
\epf

We will now classify all minimal quasi-idempotents in the $\Qlcl$-linear triangulated braided category $e\DG$. We will use the facts about induction functors from Appendix \ref{a:iqimc}, where we have studied the induction of quasi-idempotents satisfying the ``(geometric) Mackey condition''.

\bthm\label{t:minqiedgg}
(i) Let $\L\in \C(T)(\Qlcl)$ be any multiplicative local system on $T$. Let $\Pi_0'\subset \Pi_0$ be the stabilizer of $\L$. Let $\Gamma'\subset \Gamma$ (resp. $G'\subset G$) be the subgroup containing $T$ (resp. $G^\circ$) such that $\Gamma'/T= \Pi_0'$ (resp. $G'/G^\circ= \Pi_0'$). Then $\beel$ is a minimal quasi-idempotent in $e\DGp$ satisfying the Mackey criterion with respect to $G$ and $f_\L:=\indg\beel\in e\DG$ is a minimal quasi-idempotent in $e\DG$. The isomorphism class of the minimal quasi-idempotent $f_\L\in e\DG$ only depends on the $\Pi_0$-orbit of $\L$.\\
(ii) The strongly semigroupal braided functor from Proposition \ref{p:indmac} is an equivalence and induces an equivalence
\beq
\indg:{}^{\beel}e\DGp^\Delta\rar{\cong}{ }^{f_\L}e\DG^\Delta
\eeq
 of $\Qlcl$-linear triangulated braided semigroupal categories.\\
(iii) Let $f\in e\DG$ be any minimal quasi-idempotent. Then there exists a multiplicative local system $\L$ on $T$ such that  $f\sim f_\L$. Moreover, the mapping $\L\mapsto f_\L$ defines a bijection between the set of orbits $\C(T)(\Qlcl)/\Pi_0$ and the set of equivalence classes of minimal quasi-idempotents in $e\DG$.\\
(iii$'$) Each equivalence class of minimal quasi-idempotents in $e\DG$ contains a unique (up to isomorphism) indecomposable object. This indecomposable object is of the form $f_\L$ for some $\L\in \C(T)(\Qlcl)$. Any minimal quasi-idempotent in the corresponding equivalence class is of the form $W\otimes f_\L$ for some $W\in \Dmon$.
\ethm

We first prove the following:
\blem\label{l:multls}
Let $X\in e\D_U(G)$ be any nonzero object. Then there exists a multiplicative local system $\L$ on $T$ such that $\beel\ast X\neq 0$. For any $X\in e\DG$, there exists a multiplicative local system $\L$ on $T$ such that $f_\L\ast X\neq 0$.
\elem
\bpf
Without loss of generality (say by passing to an indecomposable component of $X$), we may suppose that $X\in e\D_U(G^\circ g)\cong e\D_U(Ug)\boxtimes \D(T)$ corresponds to an object of the form $A\boxtimes B$ for some $g\in G$, $A\in e\D_U(Ug)$ and $B\in \D(T)$. Then by (\ref{e:convinedug}) or by \cite[Defn. 8.7, Prop. 8.8]{De2} we have $(\eel)\ast(A\boxtimes B)\cong A\boxtimes({}^{g^{-1}}(e_\L)\ast B)$. Now using Lemma \ref{l:mellin}(ii), we complete the proof of the first part of the Lemma.

Now suppose that $X\in \DG$. Applying the forgetful functor, we may consider $X$ as an object of $\D_U(G)$. Then by the first part, there exists an $\L$ such that $\beel\ast X\neq 0$. Now by definition $f_\L=\indg \beel$ is isomorphic (as an object of $e\D_U(G)$) to a direct sum of conjugates of $\beel$. Hence we conclude that $f_\L\ast X\neq 0$.
\epf

\begin{proof}[Proof of Theorem \ref{t:minqiedgg}]
We see that $\beel\in e\DGp$ is a minimal quasi-idempotent by Corollary \ref{c:minqiedgg}. Moreover, if $x\in G-G'$, then the multiplicative local systems $\L$ and ${}^x\L$ are non-isomorphic. Hence we deduce that $\beel\ast{}^x(\beel)=0$. This means that $\beel\in e\DGp$ satisfies the Mackey condition with respect to $G$. Then by Appendix \ref{a:iqimc}, $f_\L=\indg\beel$ is a quasi-idempotent in $e\DG$. If we consider $f_\L$ as an object of $e\D_U(G)$ (or even of $e\D_{G^\circ}(G)$), it is isomorphic to a direct sum 
\beq
f_\L\cong \bigoplus\limits_{G' g\in G/G'}{}^g(\beel).
\eeq
Moreover it is also clear that for any $g\in G$, $f_{{}^g(\L)}\cong f_\L$. To complete the proof of (i), we must show that $f_\L$ is a minimal quasi-idempotent. We will first prove (ii). This would also complete the proof of (i) by Corollary \ref{c:altminqi} since we already know that $\beel\in e\DGp$ is a minimal quasi-idempotent. 

Now we have ${}^{\beel}\DGp\cong \left({}^{\beel}\D_{G^\circ}(G')\right)^{\Pi_0'}$ and ${}^{f_\L}\DG\cong{}^{f_\L}\left(e\D_{G^\circ}(G)^{\Pi_0}\right)$. Then we see that 
\beq
\indg: \left({}^{\beel}\D_{G^\circ}(G')\right)^{\Pi_0'}\rar{} { }^{f_\L}\left(e\D_{G^\circ}(G)^{\Pi_0}\right)
\eeq
is an equivalence, completing the proof of (i), (ii).

To prove (iii), note that by Lemma \ref{l:multls}, given a minimal quasi-idempotent $f\in e\DG$, there exists a multiplicative local system $\L$ on $T$ such that $f_\L\ast f\neq 0$. Since both $f$ and $f_\L$ are minimal quasi-idempotents, we conclude that $f\sim f_\L$. Moreover it is clear that $f_\L\sim f_{\L'}$ if and only if $\L$ and $\L'$ lie in the same $\Pi_0$-orbit in $\C(T)(\Qlcl)$. Statement (iii$'$) also follows using  the same argument as in the proof of Corollary \ref{c:allminqi}. 
\end{proof}
%Recall that by Corollary \ref{c:minqiedgg} we have a preferred $t$-structure on ${}^{\beel}e\DGp$ with heart $\M_{G',e,\beel}$. Also by Proposition \ref{p:indmac}(iii), for $N\in \DG$, we can consider $\beel\ast N$ as an object of ${}^{\beel}e\DGp$. Consider the $\Qlcl$-linear functor $R^0:\DG\rar{}\M_{G',e,\beel}$ defined by $N \mapsto (\beel\ast N)^0$. By Proposition \ref{p:indmac}(v), for $M\in \M_{G',e,\beel}$ we have a functorial isomorphism $\beel\ast \indg M\rar\cong \beel\ast M$, and hence we obtain functorial isomorphisms $(R^0\circ \indg) M\rar\cong \beel \uast M\rar\cong M$ after taking the $0$-th cohomology with respect to our chosen $t$-structure. This proves that $\indg:\M_{G',e,\beel}\rar{}{}^{f_\L}e\DG$ is fully faithful. Hence we conclude that $\indg:{}^{\beel}e\DGp\rar{}{}^{f_\L}e\DG$ is fully faithful.
%Let $N\in \DG$ be any object. Then by Proposition \ref{p:indmac}(iv) we have a functorial isomorphism $f_\L\ast N\rar\cong \indg(\beel\ast N_{G'})$. Recall using Remark \ref{r:vt} and Proposition \ref{p:indmac} that $f_\L$ is a $V_T$-quasi-idempotent in $\DG$. If $N\in {}^{f_\L}\DG$, then we obtain that $V_T\otimes N\cong \indg(\beel\ast N_{G'})$.

\brk\label{r:vtqi}
Note that since $e_\L\in \D_T(T)$ is a $V_T$-quasi-idempotent, each $f_\L\in e\D_G(G)$ is also a $V_T$-quasi-idempotent. In other words, each indecomposable minimal quasi-idempotent in $e\DG$ is a $V_T$-quasi-idempotent.
\erk

\subsection{The modular category and character sheaves associated with a minimal quasi-idempotent}\label{s:mccsamqi}
In this section we will define the set $CS_e(G)$ of all character sheaves associated with the Heisenberg idempotent $e\in \DG$.  As a corollary of the results of the previous section, we can associate a modular category with any minimal quasi-idempotent in $e\DG$:
\bcor\label{c:modcat}
(i) Let $f$ be any minimal  quasi-idempotent in $e\DG$. Without loss of generality  assume that $f$ is indecomposable. Then ${}^fe\DG^\Delta={}^f\DG^\Delta$ is a $\Qlcl$-linear triangulated braided semigroupal category.  It has a unique $t$-structure, denoted here by $(\D^{\leq 0}, \D^{\geq 0})$, such that $\D^{\leq 0}\ast \D^{\leq 0}\subset \D^{\leq 0}$ but $\D^{\leq 0}\ast \D^{\leq 0}\nsubseteq \D^{\leq -1}$. \\
(ii) Let $\M_{G,f}^\Delta\subset {}^{f}\D_{G}(G)^\Delta$ denote the heart of the above $t$-structure. For each $i\in \Z, X\in {}^{f}\D_{G}(G)^\Delta$ let $X^i\in \M_{G,f}^\Delta$ denote the $i$-th cohomology of $X$ with respect to this $t$-structure. Then $(\M_{G,f}^\Delta,\underline{\ast})$ has the structure of a braided semigroupal category, where 
\beq
\underline{\ast}:\M_{G,f}^\Delta\times \M^\Delta_{G,f}\rar{}\M^\Delta_{G,f}
\eeq denotes the truncated convolution functor defined by $X\uast Y:=(X\ast Y)^0$ for $X,Y\in \M^\Delta_{G,f}$. \\
(iii) Let $\M_{G,f}\subset \M_{G,f}^\Delta$ denote the full subcategory formed by the semisimple objects. Then $\M_{G,f}={}^f\DG\cap \M_{G,f}^\Delta$ is closed under truncated convolution and $(\M_{G,f},\uast,f)$ has the structure of a non-degenerate braided fusion category. There is a natural spherical structure on $\M_{G,f}$, thus giving it the structure of a modular category.\\
(iv) The Frobenius-Perron dimension $\FPdim(\M_{G,f})$ is the square of an integer.\\
(v) Suppose that the ground field $\k$ is $\Fqcl$. Then the natural spherical structure on $\M_{G,f}$ is positive integral, i.e. categorical dimensions of all objects of $\M_{G,f}$ are positive integers.
\ecor
\bpf
Statements (i), (ii), (iii), (iv) follow from Corollary \ref{c:minqiedgg}, Theorem \ref{t:minqiedgg}, Remark \ref{r:spherical} and the fact that $|K_\N|=p^{2k}$ for some $k\in \Z_{\geq 0}$ (cf. \cite{Da}). Statement (v) follows from \cite[\S7.2]{De4}.
\epf

\brk
The natural spherical structure on $\M_{G,f}$ can also be defined using the twist $\theta$ in the category $\DG$ (cf. \cite{BD}).
\erk

Let us now define character sheaves associated with the Heisenberg idempotent $e\in \DG$. First note that if $f,f'\in e\DG$ are minimal quasi-idempotents which are not equivalent, then $f\ast f'=0$ and hence ${}^f\DG\cap {}^{f'}\DG=0$. Let $\f{L}_e(G)$ denote the set of isomorphism classes of indecomposable minimal quasi-idempotents in $e\DG$, or equivalently, the set of equivalence classes of minimal quasi-idempotents in $e\DG$. By Theorem \ref{t:minqiedgg}(iii) we have an identification $\f{L}_e(G)\cong \C(T)(\Qlcl)/\Pi_0$. As we will define below, the set $\f{L}_e(G)$ parametrizes the $\f{L}$-packets of character sheaves in the category $e\DG$.
\bdefn\label{d:csfg}
(i) Let $f\in e\DG$ be an indecomposable minimal quasi-idempotent. Let $\CS_{f}(G)=\CS_{e,f}(G)$ denote the (finite) set of isomorphism classes of simple objects of the modular category $\M_{G,f}\subset {}^f\DG\subset {}^f\D_G(G)^\Delta\subset \DG$. Then we say that the finite set $\CS_f(G)$ is the $\f{L}$-packet of character sheaves on $G$ associated with the (indecomposable) minimal quasi-idempotent $f\in e\DG$. Note that in particular $f$ is itself a character sheaf.\\
(ii) The set $\CS_e(G)$ of character sheaves on $G$ associated with the Heisenberg idempotent $e\in \DG$ is defined as the (disjoint) union of all the $\f{L}$-packets of character sheaves in $e\DG$:
\beq
\CS_e(G):=\coprod\limits_{f\in \f{L}_e(G)} \CS_f(G)= \coprod\limits_{\<\L\>\in \C(T)(\Qlcl)/\Pi_0} \CS_{f_\L}(G).
\eeq
\edefn

\brk\label{r:csedgg}
Let us note that the set $CS_e(G)$ as well as its $\f{L}$-packet decomposition has been described purely in terms of the $\Qlcl$-linear triangulated braided monoidal structure of the category $e\DG$: We look at a minimal quasi-idempotent in this category and then we look at the associated full subcategory which has a canonical $t$-structure which is used to define the character sheaves. It is also clear that for each $C\in \CS_e(G)$, we have $\Hom_{\DG}(C,C)=\Qlcl$. Thus we have now completed the proof of the Heisenberg case version of Theorem \ref{t:main1}.
\erk

\subsection{Perversity of character sheaves in the Heisenberg case}\label{s:pcshc}
We will now see that all the character sheaves in the set $CS_e(G)$ (where $e$ is a Heisenberg idempotent) are perverse up to a shift. We must remark however that this will not necessarily be true for character sheaves associated with general minimal idempotents in $\DG$.

We continue using our previous notations and conventions. In particular $e=\N\boxtimes \can_H$ is the Heisenberg idempotent associated with a Heisenberg admissible pair $(H,\N)$. Note that for each $G^\circ g\in \Pi_0$, we have an equivalence $e\D_U(Ug)\boxtimes \D(T)\cong e\D_U(G^\circ g)$ which has been extensively used previously. Note that there is a perverse $t$-structure on both sides of the equivalence. 
\blem\label{l:perversity}
The equivalence $e\D_U(Ug)\boxtimes \D(T)\cong e\D_U(G^\circ g)$ respects the perverse $t$-structure on both the sides. As a consequence, the equivalence $e\D_U(Ug)\boxtimes \D^g_T(T)\cong e\D_{G^\circ}(G^\circ g)$ also  respects the perverse $t$-structure on both the sides.
\elem
\bpf
Let us first consider the canonical equivalence $e\D_U(U)\boxtimes \D(T)\cong e\D_U(G^\circ)$. This equivalence has been defined in \cite[\S8.4]{De2}. In the definition we use an extension of the Heisenberg admissible pair $(H,\N)$ to a central admissible pair $(L,\N')$. Recall that the equivalence is then defined using a certain induction (or averaging) functor which defines an equivalence
\beq
\av_{U/U_TL}:e_{\N'}\D_{U_TL}(U_TLT)\rar\cong e\D_U(G^\circ),
\eeq
where $H\subset U_T\subset U$ is such that $U_T/H=(U/H)^T$. We refer to \cite[\S8.4]{De2} for the details and notations. Now $e_{\N'}$ is a closed idempotent in $\D_{U_TL}(U_TLT)$ and $e\cong \av_{U/U_TL}(e_\N')$ is a closed idempotent in $\D_U(G^\circ)$. Hence by the argument from \cite[\S5.6]{BD}, the canonical arrow $\av_{U/U_TL} N\rar{}\h{Av}_{U/U_TL}N$ is an isomorphism for each $N\in e_{\N'}\D_{U_TL}(U_TLT)$. Then using the results and arguments from \cite[\S7.5]{BD} we see that the equivalence $\av_{U/U_TL}:e_{\N'}\D_{U_TL}(U_TLT)\rar\cong e\D_U(G^\circ)$ takes the perverse $t$-structure on $e_{\N'}\D_{U_TL}(U_TLT)$ to a suitably shifted perverse $t$-structure on $e\D_U(G)$ and also that the equivalence $e\D_U(U)\boxtimes \D(T)\cong e\D_U(G^\circ)$ preserves the perverse $t$-structures on both the sides.

Now let us consider a general connected component $G^\circ g$ and  the equivalence $e\D_U(Ug)\boxtimes \D(T)\cong e\D_U(G^\circ g)$ defined by $M\boxtimes N\mapsto M\ast \bar{N}$. Note that $e[-\dim H]$ is perverse and hence the previous argument shows that the functor $N\mapsto \bar N$ (which is defined as the composition $\D(T)\rar{}e\D_U(U)\boxtimes \D(T)\rar{}e\D_U(G^\circ)\subset e\D_U(G)$) takes a perverse sheaf in $\D(T)$ to a perverse sheaf shifted by $\dim H$. 

Note that we have the full subcategory $\tM_{Ug,e}\subset e\D_U(Ug)$ of perverse sheaves shifted by $\dim H$. By \cite[Thm. 6.2(i)]{De2} $\tM_{Ug,e}$ is a semisimple abelian category and $e\D_U(G)\cong D^b\tM_{Ug,e}$. Now the fact that the equivalence $e\D_U(Ug)\boxtimes \D(T)\cong e\D_U(G^\circ g)$ preserves the perverse $t$-structures follows from Lemma \ref{l:convperv} proved below.
\epf

\blem\label{l:convperv}
Let $M\in \tM_{Ug,e}$. Then the triangulated functor $M\ast(\cdot):e\D_U(G^\circ)\rar{}e\D_U(G^\circ g)$ is exact with respect to the perverse $t$-structures.
\elem
\bpf
By \cite[Lem. 2.1.3.1]{Li} it suffices to prove that $M\ast P$ is perverse whenever $P\in e\D_U(G^\circ)$ is perverse. Note the object $M$ has a rigid dual $M^\vee\in \tM_{Ug^{-1},e}$ (cf. \cite[Thm. 6.2(vi)]{De2}) and we have $M^\vee\ast M\cong e\oplus X\in \tM_{U,e}=\M_{U,e}$. In particular, if $X\in e\D_U(G)$ is non-zero, then $M\ast X$ must also be non-zero. Now by \cite[Prop. 5.2]{De1} we must have $M\ast P\in {}^p\D^{\geq 0}(G^\circ g)$. On the other hand, $M^\vee\ast M\ast P\in \Perv(G^\circ)$ by Lemma \ref{l:perversity} for the case of $e\D_U(G^\circ)$ which has already been established. Hence we conclude that $M\ast P$ must also be perverse. This completes the proof of the lemma.
\epf

\bthm\label{t:perversity}
(i) Let $\L\in \C(T)(\Qlcl)^{\Pi_0}$ and let $\beel\in e\D_{G^\circ}(G)$ be the corresponding minimal weakly central quasi-idempotent. Then the preferred $t$-structure on ${}^{\beel}\D_{G^\circ}(G)^\Delta$ (from Proposition \ref{p:beeltrsubcat}) is equal to the perverse $t$-structure shifted by $\dim H+\dim T$.\\
(ii) Let $\L\in \C(T)(\Qlcl)$ be any multiplicative local system on $T$ and let $f_\L\in e\D_{G}(G)$ be the corresponding minimal quasi-idempotent. Then the preferred $t$-structure on ${}^{f_\L}\D_{G}(G)^\Delta$ (cf. Corollary \ref{c:modcat}) is equal to the perverse $t$-structure shifted by $\dim H+\dim T$. In particular, all the character sheaves associated with the Heisenberg idempotent $e$ are perverse sheaves shifted by $\dim H+\dim T$.
\ethm
\bpf
To prove (i), note that by Lemma \ref{l:perversity}, the equivalence $e\D_U(Ug)\boxtimes {}^{e_\L}\D^g_T(T)^\Delta\cong {}^{\beel}\D_{G^\circ}(G^\circ g)^\Delta$ (cf. Corollary \ref{c:beeledgg}) respects the perverse $t$-structures. In the proof of Proposition \ref{p:beeltrsubcat} we have constructed the preferred $t$-structure on ${}^{\beel}\D_{G^\circ}(G^\circ g)$ using the perverse $t$-structure on $e\D_U(Ug)$ shifted by $\dim H$ and the perverse $t$-structure on ${}^{e_\L}\D^g_T(T)^\Delta$ shifted by $\dim T$. This precisely corresponds to the perverse $t$-structure on ${}^{\beel}\D_{G^\circ}(G^\circ g)$ shifted by $\dim H+\dim T$. Statement (i) now follows.

Statement (ii) follows from (i) using Corollary \ref{c:minqiedgg}, Theorem \ref{t:minqiedgg} and Corollary \ref{c:modcat}. Since the character sheaves in the $\f{\L}$-packet associated with the minimal quasi-idempotent $f_\L$ are defined to lie in the heart of this preferred $t$-structure (cf. Definition \ref{d:csfg}) we see that all character sheaves associated with $e$ are perverse sheaves shifted by $\dim H+\dim T$.
\epf

\section{Irreducible characters associated with Heisenberg idempotents}\label{s:icahi}
In this section, we assume that our base field $\k$ is equal to $\Fqcl$ and that our neutrally solvable group $G$ is equipped with an $\Fq$-Frobenius map $F:G\rar{}G$. As we have stated before, we should also consider all pure inner forms (which are parametrized by the finite set $H^1(F,G)=H^1(F,\Pi_0)$) of the Frobenius.

Now suppose that the Heisenberg idempotent $e\in \DG$ (coming from the Heisenberg admissible pair $(H,\N)$) is $F$-stable. Since $e$ is supported on $H$, $H$ must be an $F$-stable (connected normal) subgroup. Moreover, the multiplicative local system $\N\in H^*$ must also be $F$-stable. Hence the Heisenberg admissible pair $(H,\N)$ is $F$-stable and hence it comes by base change from a Heisenberg admissible pair $(H_0,\N_0)$ defined over $\Fq$. Let $e_0\in \DGn$ be the corresponding Heisenberg idempotent. By taking the sheaf-function correspondence (cf. \S\ref{s:i}) for the multiplicative local system $\N_0$, we obtain for each of the forms $H^t_0(\Fq)$, the corresponding character $\T^t_{\N_0}:H^t_0(\Fq)\rar{}\Qlcl^\times$ and we may combine these into the function $\T_{\N_0}\in \FunG$ by extension by zero outside $H$. Similarly, we have the function $\T_{e_0}=\frac{\T_{\N_0}}{q^{\dim H}}\in \FunG$. Moreover $\T_{e_0}$ is an idempotent in $\FunG$ and is supported only $H$.  Consider the Hecke subalgebra $\T_{e_0}\FunG$ with unit $\T_{e_0}$. Let $\Irrep_{e}(G,F)\subset \Irrep(G,F)$ denote the set of irreducible representations (of all pure inner forms of $G^F$) such that the idempotent $\bar{\T_{e_0}}$ (``complex conjugate'' of the idempotent $\T_{e_0}$) acts by the identity. Taking the characters of the irreducible representations, we can consider $\Irrep_{e}(G,F)$ as an orthonormal basis of $\T_{e_0}\FunG$. In this section our goal is to study the set $\Irrepe$ and its relationship with the set $CS_e(G)^F$ of $F$-stable character sheaves associated with $e$.

\brk\label{r:fqform}
The $\Fq$-form $G_0$ of $G$ also defines for us the $\Fq$-forms $U_0, G^\circ_0$ and $T_0:=G^\circ_0/U_0$ of $U, G^\circ$  and $T$ respectively.
\erk

Let us now interpret the set $\Irrepe$ in terms of the category $\D^F_G(G)\cong \D_G(GF)$ of $F$-twisted conjugation equivariant $\Qlcl$-complexes on $G$ that is studied in \cite[\S2.4]{De3}. Recall from {\it loc. cit.}  that $\D^F_G(G)$ is equivalent to the bounded derived category of the category of representations of all pure inner forms of $G^F$. In particular $\D^F_G(G)$ also happens to be a semisimple abelian category with simple objects parametrized by $\Irrep(G,F)\times \Z$. (Here $\Z$ corresponds to the shifts in degree.) Given an irreducible representation $W\in \Irrep(G,F)$, we have the associated local system $W_{loc}\in \Sh^F_G(G)\subset \D^F_G(G)$.  Moreover, recall that $\D^F_G(G)\cong \D_G(GF)$ is a $\DG$-module category.  The following result is proved in \cite[Prop. 6.5]{De3}:
\bprop\label{p:irrepe}
If $W\in \Irrepe$, then $W_{loc}$ lies in the full subcategory $e\D^F_G(G)\subset \D^F_G(G)$. The map $W\mapsto W_{loc}$ sets up a bijection between the set $\Irrepe$ and the set of (isomorphism classes of) simple objects in $e\Sh^F_G(G)$.
\eprop

Combining the results above, let us state the following equivalent characterizations of the set $\Irrepe$:
\bcor\label{c:irrepe}
Let $W\in \Irrep(G,F)$ be an irreducible representation of a pure inner form, say $G^t_0(\Fq)$. Then the following are equivalent:\\
(i) $W\in \Irrepe$.\\
(ii) The subgroup $H^t_0(\Fq)\subset G^t_0(\Fq)$ acts on $W$ by the character $\T_{\N_0}^t:H^t_0(\Fq)\rar{}\Qlcl^\times$.\\
(iii) The dual (i.e., `complex conjugate') idempotent $\bar{\T_{e_0}}\in\FunG$ acts on $W$ trivially.\\
(iv) The character $\chi_W$ lies in $\T_{e_0}\FunG\subset \FunG$.\\
(v) $W_{loc}\in e\D_G(GF)$.\\
The set $\{\chi_W|W\in \Irrepe\}$ forms an orthonormal basis of $\T_{e_0}\FunG$.
\ecor

\subsection{$F$-stable minimal quasi-idempotents}
We have defined $\fL$-packets of character sheaves associated with $e$ in \S\ref{s:mccsamqi}. These $\f{L}$-packets are parametrized by the set $\fL_e(G)$ of isomorphism classes of indecomposable minimal quasi-idempotents in $e\DG$. Recall that by Theorem \ref{t:minqiedgg}, we have a canonical identification $\fL_e(G)=\C(T)(\Qlcl)/\Pi_0$. 

Note that in our current situation we also have a braided triangulated auto-equivalence $F^*:e\DG\rar{\cong}e\DG$ with inverse $F:=F_*:e\DG\rar{\cong}e\DG$. In particular, this induces a permutation of the set $\fL_e(G)$ as well as the set $CS_e(G)$ of character sheaves. Using Theorem \ref{t:minqiedgg} (and the same notation), we obtain:
\blem\label{l:faction}
For each $\L\in \C(T)(\Qlcl)$ we have canonical isomorphisms $$F^*f_\L=F^*\indg(\beel)\cong \ind_{F^{-1}(G)}^G(F^*\beel)\cong \ind_{F^{-1}(G)}^G(\bar{e_{F^*\L}})=f_{F^*\L}.$$ Hence the permutation of $\fL_e(G)$ induced by $F^*$ matches with the permutation of $\C(T)(\Qlcl)/\Pi_0$ induced by $F^*$.
\elem 

We want to partition the set $\Irrepe$ into what we will call $\fL$-packets of irreducible characters associated with the $F$-stable Heisenberg idempotent $e$. We will see that these $\fL$-packets are parametrized by the set $\fL_e(G)^F$.

Now we have the canonical triangulated monoidal functor $\D(T)\rar{}e\D_U(G^\circ)$ which is fully faithful onto a direct summand of the category $e\D_U(G^\circ)$. We will now prove that there is a canonical monoidal functor $\D(T_0)\rar{}e_0\D_{U_0}(G_0^\circ)$ which is fully faithful onto a direct summand of the category $e\D_{U_0}(G_0^\circ)$ and which is compatible with the functor $\D(T)\rar{}e\D_U(G^\circ)$ and extension of scalars.

Now by the construction of \cite[\S8.4]{De2}, the functor $\D(T)\rar{}e\D_U(G^\circ)$ can be constructed using a certain central admissible pair $(L,\N')$ for $G^\circ$. Now this central admissible pair is defined over $\F_{q^n}$ for some positive integer $n$, i.e. it comes from a central admissible pair $(L_1,\N_1')$ defined over $\F_{q^n}$ for $G^\circ_1:=G_0\otimes_{\F_{q}}\F_{q^n}$. Let $e_1$ be obtained from $e_0$ by base change to $\F_{q^n}$. Hence using the construction of \cite[\S8.4]{De2} over the field $\F_{q^n}$, we obtain a canonical functor $\D(T_1)\rar{}e_1\D_{U_1}(G^\circ_1)$ which is compatible with the $\Fq$-Frobenius.

Now an object of $\D(T_0)$ can be considered as an object of $\D(T_1)$ equipped with some descent data for the finite \'etale cover $T_1\rar{}T_0$. Now such descent data defines for us similar descent data for $e_1\D_{U_1}(G^\circ_1)$, or equivalently an object of $e_0\D_{U_0}(G^\circ_0)$. Hence we have proved:
\bprop\label{p:descent}
We have a canonical triangulated monoidal functor $\D(T_0)\rar{}e_0\D_{U_0}(G_0^\circ)$ which is fully faithful onto a direct summand of the category $e_0\D_{U_0}(G_0^\circ)$ and the following diagram commutes up to a natural isomorphism
$$\xymatrixcolsep{5pc}\xymatrix{
\D(T_0)\ar[r]\ar[d] & e_0\D_{U_0}(G^\circ_0)\ar[d]\\
\D(T)\ar[r] & e\D_U(G^\circ).\\
}$$
Similarly, we have a canonical triangulated monoidal functor $\D_{T_0}(T_0)\rar{}e_0\D_{G^\circ_0}(G_0^\circ)$ which is fully faithful onto a direct summand of the category $e_0\D_{G^\circ_0}(G_0^\circ)$.
\eprop

By Remark \ref{r:vtqi}, each indecomposable minimal quasi-idempotent in $e\DG$ is a $V_T$-quasi-idempotent. Note that the object $V_T=H^*_c(T,\can_T)\in \D_G(\h{Spec}(\k))$ can be obtained by base change from the object $V_{T_0}=H^*_c(T_0,\can_{T_0})\in \D_{G_0}(\h{Spec} \Fq)$. Also note that since $e_0$ is perverse up to a shift, it is the unique weak idempotent in $\DGn$ (up to isomorphism) whose base change to $\DG$ is isomorphic to $e$ (cf. \cite[\S6.1]{B}). 
\bprop\label{p:fstminqi}
(i) Let $f\in e\DG$ be an indecomposable minimal quasi-idempotent such that $F^*f\cong f.$  Then there exists a unique (up to isomorphism) $f_0\in \DGn$  whose base change to $\DG$ is isomorphic to $f$ and such that $f_0\ast f_0\cong V_{T_0}\otimes f_0$.\\
(ii) Let $f$ be as above. Say $f\cong f_\L$ with $\L\in \C(T)(\Qlcl)$. Then since $f$ is $F$-stable, the $\Pi_0$-orbit of $\L$ in $\C(T)(\Qlcl)$ is $F$-stable, i.e. $F^*\L\cong g(\L)$ for some $g\in \Pi_0$, i.e. $(gF)^*\L\cong \L$. Then there exists a unique multiplicative local system $\L_0^g$ on the $\Fq$-form $T_0^g$ (corresponding to the Frobenius $gF:T\to T$) whose base change to $T$ is $\L$. Then the $f_0$ from (i) can be constructed as $f_0:=f_{\L_0^g}\in \D_{G_0^g}(G_0^g)\cong\DGn$.
\eprop
\bpf
We have seen in Theorem \ref{t:perversity} that each indecomposable minimal quasi-idempotent $f$ in $e\DG$ is a perverse sheaf shifted by $\dim H+\dim T$. Hence the uniqueness part of (i) follows from \cite[\S6.1]{B}. To show existence, it suffices to prove (ii). Note that by Lemma \ref{l:faction}, we have $F^*f_\L\cong f_{F^*\L}$. Hence if $f_\L$ is $F$-stable we must have $(gF)^*\L\cong \L$ for some $g\in \Pi_0$. Hence $\L$ comes by base change from a unique multiplicative local system $\L_0^g$ on $T_0^g$. Here $T_0^g$ is the $\Fq$-form of $T$ corresponding to the Frobenius $gF:T\to T$. Then using Proposition \ref{p:descent} we can construct the object $f_{\L_0^g}\in \D_{G_0^g}(G^g_0)$ such that we have $f_{\L_0^g}\ast f_{\L_0^g}\cong V_{T_0^g}\otimes f_{\L_0^g}$. It is clear that by base change of $f_{\L_0^g}$ to $\DG$ we obtain $f_\L$. Moreover by \cite[\S4.4]{B} we have an identifications $\D_{G_0^g}(G_0^g)\cong\DGn$ and $\D_{G_0^g}(\h{Spec}(\Fq))\cong\D_{G_0}(\h{Spec}(\Fq))$ under which $V_{T_0^g}$ maps to $V_{T_0}$. This completes the proof.
\epf

\brk\label{r:bijctql}
The previous result means that given an $F$-stable $\fL$-packet associated with $e$ (i.e. an element of the set $(\C(T)(\Qlcl)/\Pi_0)^F=\fL_e(G)^F$) we can define the corresponding object $f_0\in e_0\DGn$. 
\erk

\brk\label{r:trf}
Let $f\in \f{L}_e(G)$ be an $F$-stable indecomposable minimal quasi-idempotent. We have $V_{T_0}\in \D_{G_0}(\h{Spec} \Fq)\subset \DGn$ considered as a complex supported at $1\in G_0$. Then we have seen in Proposition \ref{p:fstminqi} that $f_0\ast f_0\cong {V_{T_0}}\ast f_0$ and hence 
\beq\label{e:trf}\T_{f_0}\ast \T_{f_0}=\T_{V_{T_0}}\ast \T_{f_0}.\eeq
 Note that the function $\T_{V_{T_0}}$ is only supported on $1\in G$ and corresponding to the inner form $G^t_0(\Fq)$ we have $\T^t_{V_{T_0}}(1)=\frac{|T^t_0(\Fq)|}{q^{\dim T}}$. Consider the function 
$${\t{\T}}_{f_0}^t:=\frac{q^{\dim T}}{|T^t_0(\Fq)|}\cdot \T_{f_0}^t\in \Fun(G^t_0(\Fq)/\sim).$$ 
Then from (\ref{e:trf}), we conclude that $\t{\T}_{f_0}^t$ is an idempotent. Combining these functions over all the pure inner forms we obtain the idempotent $\t{\T}_{f_0}\in \T_{e_0}\FunG\subset \FunG$.
\erk

\subsection{$\fL$-packets of irreducible characters in $\Irrep_{e}(G,F)$}\label{s:lpicirrepe}
In this section, we will partition the set $\Irrepe$ into what we will call $\fL$-packets of irreducible characters which will be parametrized by the set $\fL_e(G)^F$. Recall that the category $e\D_G(GF)\cong e\D^F_G(G)$ is the bounded derived category of a finite semisimple category and hence is itself a semisimple abelian triangulated category. The simple objects of $e\D_G(GF)$ are parametrized by $\Irrepe\times \Z$ (cf. Proposition \ref{p:irrepe}) with $\Z$ corresponding to the degree shift.

We have (cf. \cite[Prop. 8.10]{De2})
\beq\label{e:edggf}
e\D_{G^\circ}(GF)=\bigoplus\limits_{G^\circ g\in \Pi_0}e\D_{G^\circ}(G^\circ gF)\cong \bigoplus\limits_{G^\circ g\in \Pi_0}e\D_{U}(U gF)\boxtimes \D_T^{gF}(T).
\eeq
Hence we see that each simple object of $e\D_{G^\circ}(GF)$ is of the form $\bar{M\boxtimes N}$, where $M$ is a simple object of $e\D_{U}(U gF)$ and $N$ a simple object of $\D^{gF}_T(T)$ for some $g\in G$. 
\blem\label{l:edgn}
Consider a connected component $G^\circ g\subset G$ and the Frobenius $gF:G^\circ \to G^\circ$. Let $W\in \Irrep_{e}(G^\circ, gF)=\Irrep_{e}({G^\circ}^{gF}).$ Then there exists a unique (up to isomorphism) $gF$-stable indecomposable minimal quasi-idempotent $\beel$ (where $\L\in \C(T)(\Qlcl)^{gF}$) in $e\D_{G^\circ}(G^\circ)$ such that $W_{loc}\in {}^{\beel}\D_{G^\circ}(G^\circ gF)$.
\elem
\bpf
By Proposition \ref{p:irrepe}, $W_{loc}\in e\D_{G^\circ}(G^\circ gF)\cong e\D_{U}(U gF)\boxtimes \D_T^{gF}(T)$ is a simple object and hence corresponds to an object of the form $M\boxtimes N\in e\D_{U}(U gF)\boxtimes \D_T^{gF}(T)$, with $M,N$ simple. In particular, $N$ corresponds to an irreducible character of $T^{gF}$, which in turn corresponds to a unique $gF$-stable multiplicative local system $\L$ on $T$. In fact we can check that $N$ must be isomorphic to $\L$ considered as an object of $\D^{gF}_T(T)$ (cf. \cite[\S6.1]{De3}). Then it is clear that $W_{loc}=\bar{M\boxtimes \L}\in {}^{\beel}\D_{G^\circ}(G^\circ gF)$ (cf. Corollary \ref{c:beeledgg}). The uniqueness is clear.
\epf

\bprop\label{p:lpdec}
Let $W\in \Irrepe$. Then there exists a unique (up to isomorphism) $F$-stable indecomposable minimal quasi-idempotent $f\in e\DG$ such that $W_{loc}\in {}^f\D_G(GF)={}^f\D_G(GF)^\Delta$. 
\eprop
\bpf
Let $W$ be an irreducible representation of the inner form $G^{gF}$ for some $g\in G$. Let $\O\subset G$ denote the $F$-twisted conjugacy class of $g$. Then by Lang's theorem we must have $\O=\coprod G^\circ g_i$ where the set $\{G^\circ g_i\}\subset \Pi_0$ is the $F$-twisted conjugacy class of the element $G^\circ g$ in $\Pi_0$. Then we have 
\[W_{loc}\in e\D_G(\O F)\cong e\D_{G^\circ}(\O F)^{\Pi_0}\cong \left(\bigoplus e\D_{G^\circ}(G^\circ g_iF)\right)^{\Pi_0}\cong e\D_{G^\circ}(G^\circ gF)^{\Pi_0^{gF}}.\]

Now let $W'\in \Irrep_{e}(G^\circ, gF)$ be an irreducible representation appearing in the restriction of $W$ to the subgroup ${G^\circ}^{gF}$. Now let $\L$ be the unique $gF$-stable multiplicative local system on $T$ guaranteed by Lemma \ref{l:edgn}, such that $W'_{loc}\in {}^{\beel}\D_{G^\circ}(G^\circ gF)$. 

Now using the same notation as we have used before, let $\Pi_0'\subset \Pi_0$ denote the stabilizer of $\L\in \C(T)(\Qlcl)$ and let $G'\subset G$ be such that $G'/G^\circ=\Pi_0'$. Then from Theorem \ref{t:minqiedgg} we know that $\beel\in e\DGp$ is a minimal quasi-idempotent satisfying the Mackey criterion with respect to $G$. Since $\L$ is $gF$-stable, $G'$ is a $gF$-stable subgroup of $G$.

We will now construct an irreducible representation $W_\L$ of ${G'}^{gF}$. Let $W_\L\subset W$ be the direct sum of all the isotypic components for the action of ${G^\circ}^{gF}$ on $W$ which correspond to the local system $\L$ as above. Then by the Clifford theory of finite groups, we see that ${G'}^{gF}$ acts on the space $W_\L$, that $W\cong \ind_{{G'}^{gF}}^{G^{gF}}W_\L$ and that $W_\L$ is irreducible. Moreover by construction, we have ${W_\L}_{loc}\in {}^{\beel}\D_{G'}(G'gF)$.

As before, we set $f_\L=\indg \beel$. This is an indecomposable minimal quasi-idempotent in $e\DG$. Since $\L$ is $gF$-stable, $f_\L$ is $F$-stable. Moreover we have (cf. Appendix \ref{a:iqimc} and \cite[\S4.3]{De3}) $W_{loc}\cong \indg {W_{\L}}_{loc}\in {}^{f_\L}\D_{G}(GF)$ as desired. The uniqueness is clear, since the convolution of distinct indecomposable minimal quasi-idempotents is zero.
%so we will consider $W_{loc}$ as an object of the $\Pi_0$-equivariantization of $e\D_{G^\circ}(GF)$. Applying the forgetful functor, $e\D_G(GF)\to e\D_{G^\circ}(GF)$, consider $W_{loc}$ as an object of $e\D_{G^\circ}(GF)$ and let $W'_{loc}$ in some $e\D_{G^\circ}(G^\circ gF)$ be a simple direct summand with $W'\in \Irrep_{e_0}({G^\circ}{gF})$. Since $W_{loc}$ is a simple object of $e\D_G(GF)$ all simple direct summands of $W_{loc}$ (when considered as an object of $e\D_{G^\circ}(GF)$) are given by the orbit of $W'_{loc}$ under the action of $\Pi_0$.
\epf

We can now define the $\fL$-packet decomposition of $\Irrepe$.
\bdefn\label{d:deflpirre}
Let $f$ be an $F$-stable indecomposable minimal quasi-idempotent in $e\DG$. We define 
\[\Irrep_{e,f}(G,F):=\{W\in \Irrepe|W_{loc}\in {}^f\D_G(GF)\}\]
and call this set the $\fL$-packet of irreducible representations associated with $f$. By Proposition \ref{p:lpdec} we have
\[\Irrepe=\coprod\limits_{f\in \fL_e(G)^F}\Irrep_{e,f}(G,F).\]
\edefn

In \S\ref{s:mccsamqi} we defined character sheaves in $e\DG$ in terms of minimal quasi-idempotents $f\in e\DG$ and by considering the full subcategory (denoted by $\M_{G,f}$) of semisimple objects in the heart  of a certain canonical $t$-structure on the triangulated category ${}^f\DG^\Delta\subset e\DG$. Using the same arguments we obtain similar results for the category $e\D_G(GF)$.

\bprop\label{p:tstrmodcat}
(i) Let $f$ be an $F$-stable indecomposable minimal quasi-idempotent in $e\DG$. Then ${}^f\D_G(GF)={}^f\D_G(GF)^\Delta$ is a $\Qlcl$-linear triangulated (as well as semi-simple abelian) module category (under convolution with compact supports) over the triangulated braided semigroupal category ${}^f\DG^\Delta$. There is a unique $t$-structure on ${}^f{\D_G(GF)^\Delta}$ which is compatible under convolution (cf. Corollary \ref{c:modcat}(i)) with the canonical $t$-structure on ${}^f\DG$. This $t$-structure is equal to the perverse $t$-structure shifted by $\dim H+\dim T$.\\
(ii) Let $\M_{GF,f}\subset {}^f\D_G(GF)$ denote the heart of the above $t$-structure. Then as in Corollary \ref{c:modcat}(ii) we can define a truncated convolution
\beq
\uast:\M_{G,f}\times \M_{GF,f}\rar{}\M_{GF,f}
\eeq
which provides $\M_{GF,f}$ with the structure of an invertible $\M_{G,f}$-module category. Moreover, there is a natural $\M_{G,f}$-module trace $\tr_{F,f}$ on $\M_{GF,f}$.\\
(iii) The map $W\mapsto M_W:=W_{loc}[\dim G+\dim H+\dim T]$ defines a bijection between the set $\Irrep_{e,f}(G,F)$ and the set $\O_{\M_{GF,e}}$ of (isomorphism classes of) simple objects of $\M_{GF,f}$.
\eprop
\bpf
The proofs of statements (i) and (ii) are the same as that of Corollary \ref{c:modcat} using (\ref{e:edggf}). Note that the $\M_{G,f}$-module category which is inverse to $\M_{GF,f}$ can be constructed as the heart of a preferred  $t$-structure on the category ${}^f\D_G(GF^{-1})$ (see also \cite[\S6]{ENO2}). Also note that for each $n\in \Z$, we have a duality functor $\M_{GF^n,f}\rar{}\M_{GF^{-n},f}$ given by $\f{D}^{-}(\cdot)[2\dim H+2\dim T](\dim H+\dim T)$. Hence as in Remark \ref{r:spherical} we can define a natural spherical structure on the monoidal category $\bigoplus\limits_{n\in \Z}\M_{GF^n,f}$. This gives us the desired natural $\M_{G,f}$-module trace in the category $\M_{GF,f}$. Statement (iii) is clear from Proposition \ref{p:lpdec} and Definition \ref{d:deflpirre}.
\epf

\brk
In view of the definition of character sheaves and the result above, we may think of the irreducible characters $\Irrepe$ as character sheaves in the category $e\D_G(GF)$. Recall from \cite[\S2.4.7]{De3} that associated with each object $M\in \D_G(GF)$ is its character $\chi_M\in \FunG$. In particular for each $W\in \Irrepe$ we have the function $\chi_{M_W}=(-1)^{\dim G+\dim H+\dim T}\chi_W=(-1)^{2d_e}\chi_{W}\in \T_{e_0}\FunG$, where $d_e=\frac{\dim U-\dim H}{2}=\frac{\dim G-\dim H-\dim T}{2}$ is the functional dimension of the Heisenberg idempotent $e$.
\erk

\brk\label{r:+tr}
We know that the modular category $\M_{G,f}$ is positive integral by Corollary \ref{c:modcat}. Let $\tr^+_{F,f}$ denote the positive $\M_{G,f}$-module trace on $\M_{GF,f}$ normalized according to \cite[\S1.3]{De5}. The natural trace $\tr_{F,f}$ also satisfies the previous condition. Hence we deduce that $\tr^+_{F,f}=\pm \tr_{F,f}$.
\erk

\subsection{The relationship between irreducible characters and character sheaves}\label{s:rbiccs}
Let us consider the action of the Frobenius on the set $CS_e(G)$. If $C\in CS_e(G)$ is such that we have an isomorphism $\psi:F^*C\rar{\cong}C$, then we have the associated trace of Frobenius function $\T_{C,\psi}\in \T_{e_0}\FunG$. For each $C\in CS_e(G)^F$, let $\psi_C:F^*C\rar{\cong}C$ denote some choice of an isomorphism. In this section we will describe the relationship between the two sets $\Irrepe, \{\T_{C,\psi}|C\in CS_e(G)^F\}\subset \T_{e_0}\FunG$. We have partitions
\beq
CS_e(G)^F=\coprod\limits_{f\in \fL_e(G)^F}CS_{e,f}(G)^F \mbox{ and }
\eeq
\beq
\Irrepe=\coprod\limits_{f\in \fL_e(G)^F}\Irrep_{e,f}(G,F).
\eeq
We know that $\Irrepe$ forms an orthonormal basis of $\T_{e_0}\FunG$.

Now if $f,f'\in \fL_e(G)^F$ are distinct, and if $C\in {}^f\DG$ and $M\in {}^{f'}\D_G(GF)$, then $C\ast M=0.$ Hence by \cite[Thm. 2.14]{De3} we see that the sets $\{\T_{C,\psi_C}|C\in CS_{e,f}(G)^F\}$ and $\Irrep_{e,f'}(G,F)$ are orthogonal to each other.

In other words, we are reduced to finding the relationship between the sets $\{\T_{C,\psi_C}|C\in CS_{e,f}(G)^F\}$ and $\Irrep_{e,f}(G,F)$ which correspond to the same $\fL$-packet. Let $C\in CS_{e,f}(G)^F$ equipped with an isomorphism $\psi_C:F^*C\rar{\cong}C$ and let $W\in \Irrep_{e,f}(G,F)$. Then we have $C\in \M_{G,f}$ and $M_W\in \M_{GF,f}$. We have the associated character $\chi_{M_W}=(-1)^{2d_e}\chi_W\in \T_{e_0}\FunG$. The inner product between the functions $\T_{C,\psi_C}$ and $\chi_{M_W}$ is given by (cf. \cite[Thm. 2.14]{De3})
\beq
\<\T_{C,\psi_C},\chi_{M_W}\>=\tr_F(\g_{C,\psi_C,M_W}),
\eeq
where $\g_{C,\psi_C,M_W}$ is the following composition in $^f\D_G(GF)$:
\beq
C\ast M_W\xto{\beta_{M_W,F^*C}^{-1}} M_W\ast F^*C\xto{\beta^{-1}_{F^*C,M_W}} F^*C\ast M_W\xto{\psi_C\ast \id_{M_W}}C\ast M_W
\eeq
with $\beta_{\cdot,\cdot}$ denoting the crossed braiding isomorphisms (cf. \cite[\S2.4.4]{De3}). Hence the matrix relating the sets $\{\T_{C,\psi_C}|C\in CS_{e,f}(G)^F\}$ and $\{\chi_W|W\in\Irrep_{e,f}(G,F)\}$ is equal to the matrix $\t{S}(F,f)$ whose entries are defined by 
\beq\label{e:defts}
\t{S}(F,f)_{C,W}:=\<\T_{C,\psi_C},\chi_{W}\>=\<\T_{C,\psi_C},(-1)^{2d_e}\chi_{M_W}\>=(-1)^{2d_e}\tr_F(\g_{C,\psi_C,M_W}),
\eeq
for $C\in \CS_{e,f}(G)^F ,W\in \Irrep_{e,f}(G,F)$.

Our next goal is to express the above transition matrix in terms of a certain crossed S-matrix. Truncating the automorphism $\g_{C,\psi_C,M_W}$ with respect to the canonical $t$-structure on $^f\DG$, we obtain the following composition in $\M_{GF,f}$:
\beq
\g_{C,\psi_C,M_W}^0:C\uast M_W\xto{\underline{\beta}_{M_W,F^*C}^{-1}} M_W\uast F^*C\xto{\underline{\beta}^{-1}_{F^*C,M_W}} F^*C\uast M_W\xto{\psi_C\uast \id_{M_W}}C\uast M_W.
\eeq
Now we have seen in Proposition \ref{p:tstrmodcat} that the invertible $\M_{G,f}$-module category $\M_{GF,f}$ is equipped with a canonical module trace $\tr_{F,f}$. We also have on $\M_{GF,f}$ the positive module trace $\tr^+_{F,f}$. Hence we are now in the situation studied in \cite{De5}. Consider the crossed S-matrix (cf. \cite[\S2]{De5}) $S^+(F,f)$  whose entries are defined by 
\beq
S^+(F,f)_{C,W}:=\tr^+_{F,f}(\g^0_{C,\psi_C,M_W})
\eeq 
for $C\in \CS_{e,f}(G)^F ,W\in \Irrep_{e,f}(G,F)$.
Note that in this case the modular auto-equivalence of $\M_{G,f}$ corresponding to the module category $\M_{GF,f}$ is given by $F=F_*:\M_{G,f}\rar{}\M_{G,f}$. In addition, let us assume that 
\bit
\item[$\bigstar$] the isomorphisms $\psi_C:F^*C\rar{\cong}C$ are chosen to satisfy the further requirements as in {\it loc. cit.} 
\eit
Then by \cite[Thm. 2.9]{De5} the entries of $S^+(F,f)$ are cyclotomic integers and the matrix $S^+(F,f)$ is unitary up to a scaling, namely
\beq\label{e:defs+}
S^+(F,f)\cdot \bar{S^+(F,f)}^T=\dim \M_{G,f}.\cdot I = \bar{S^+(F,f)}^T\cdot S^+(F,f),
\eeq 
where $\dim \M_{G,f}$ denotes the categorical dimension of the modular category $\M_{G,f}$.

%Let us consider the pairs $(C\ast M_W, \g_{C,\psi_C,M_W})$ and $(C\uast M_W, \g^0_{C,\psi_C,M_W})$ as objects of the $\Z$-equivariantization $\D_G(GF)^\Z$. Recall also that we have the object $V_{T_0}:=H^*_c(T,\can_T)\in \D_{G_0}(\h{Spec}\Fq)$. We may consider it as an object of $\D^{Weil}_{G_0}({\h{Spec} \Fq})$. Note that we have a natural action of $\D_G(\pt)^\Z$ on $\D_G(GF)^\Z$ which we will denote by $\tensor$. 
Let us now describe the relationship between the two traces $\tr_F$ and $\tr^+_{F,f}$ encountered above as well as that between the matrices $\t{S}(F,f)$ and $S^+(F,f)$:
\bthm\label{t:relheis}
(i) Let $M\in \M_{GF,f}\subset {}^f\D_G(GF)$ be an object supported on a single $G$-orbit $\O F\subset GF$ and let $\alpha:M\rar{}M$ be an endomorphism. Then
\beq
\tr_F(\alpha)=(-1)^{2d_e}\frac{q^{d_e}}{q^{\dim G}\cdot \sqrt{\dim \M_{G,f}}}\cdot\frac{q^{\dim T}}{|T_0^t(\Fq)|}\cdot\tr^+_{F,f}{(\alpha)},
\eeq where $tF$ is a point in the orbit $\O F$. %We have $$(C\ast M_W,\g_{C,\psi_C,M_W})\cong V_{T_0} \otimes (C\uast M_W, \g^0_{C,\psi_C,M_W})$$ as objects of $\D_G(GF)^\Z$. 
\\
(ii) Let $W\in \Irrep_{e,f}(G,F)$ be an irreducible representation of the pure inner form $G^t_0(\Fq)$. Then we have 
\beq
\dim W = \frac{\dim^+(M_W)}{\sqrt{\dim \M_{G,f}}}\cdot |\Pi_0^{tF}|\cdot q^{d_e},
\eeq where $\dim^+(M_W):=\tr_{F,f}^+(\id_{M_W})=\FPdim(M_W)$.\\
(iii) Let $C\in \M_{G,f}$, $\psi:F^*(C)\rar{\cong}C$ and $M\in \M_{GF,f}$. Then
\beq\label{e:trftr+}
\tr_F(\g_{C,\psi,M})=(-1)^{2d_e}\frac{q^{d_e}}{q^{\dim G}\cdot \sqrt{\dim \M_{G,f}}}\cdot \tr_{F,f}^+(\g^0_{C,\psi,M}) \hbox{ and hence}
\eeq
\beq
\t{S}(F,f)=\frac{q^{d_e}}{q^{\dim G}\cdot \sqrt{\dim \M_{G,f}}}\cdot S^+(F,f).
\eeq
\ethm
\brk
Note that ${\dim \M_{G,f}}$ is a perfect square by Corollary \ref{c:modcat} and we use its positive integral square root in the theorem above. Also, $d_e\in \frac{1}{2}\Z$ and we may need to use the positive square root of $q$ according to our chosen isomorphism $\Qlcl\cong \f{C}$.
\erk
We will prove this result in \S\ref{s:pot:relheis}. Let us now derive some corollaries.
\bcor\label{c:relheis}
Let $W\in \Irrep(G,F)$ be an irreducible representation of a pure inner form and let $f\in \f{L}_e(G)^F$. Then the following are equivalent:\\
(i) The representation $W$ lies in the $\f{L}$-packet $\Irrep_{e,f}(G,F)$.\\
(ii) The local system $W_{loc}$ lies in ${}^f\D_G(GF)$. \\
(iii) The character $\chi_W$ lies in the subspace $\t{\T}_{f_0}\FunG\subset \FunG$.\\
(iv) The dual (i.e., `complex conjugate') idempotent $\bar{\t{\T}_{f_0}}$ acts on $W$ as identity.\\
The set $\{\chi_W|W\in \Irrep_{e,f}(G,F)\}$ forms an orthonormal basis of the space $\t{\T}_{f_0}\FunG$.
\ecor
\bpf
The equivalence (i)$\Leftrightarrow$(ii) follows from Definition \ref{d:deflpirre} and the fact that ${}^f\D_G(GF)\subset e\D_G(GF)$. The equivalence (iii)$\Leftrightarrow$(iv) follows from standard character theory of finite groups. Let us now prove (i)$\Rightarrow$(iii). So suppose that $W\in \Irrep_{e,f}(G,F)$. Let $M_W\in \M_{GF,f}$ be the corresponding simple object. Then by \cite[Rem. 2.2]{De5}, we see that the row of the crossed S-matrix $S^+(F,f)$ corresponding to the unit object $f\in \M_{G,F}$ contains only nonzero entries. In particular, by (\ref{e:defs+}), $\tr_{F,f}^+(\g^0_{f_0, M_W})\neq 0$. Hence by Theorem \ref{t:relheis}(iii) and \cite[Thm. 2.14]{De3}, $\tr_F(\g_{f_0,M_W})=\<\T_{f_0},\chi_{M_W}\>=\l(\T_{f_0}\ast \chi_{M_W})\neq 0$. Hence $\T_{f_0}\ast \chi_{W}\neq 0$. This implies (iii).

Let us now prove that (iii)$\Rightarrow$(i). Suppose that $\chi_W\in \t{\T}_{f_0}\FunG$, i.e. $\t{\T}_{f_0}\ast \chi_W=\chi_W$. Hence we can prove that $\<\T_{f_0},\chi_{W}\>=\tr_{F}(\g_{f_0,W_{\loc}})=\l(\T_{f_0}\ast \chi_{W})\neq 0$. Hence $f\ast W_{\loc}\neq 0$ and consequently $W\in \Irrep_{e,f}(G,F)$.

The last statement of the corollary follows since $\{\chi_W|W\in \Irrep(G,F)\}$ forms an orthonormal basis of $\FunG$.
\epf

\bcor\label{c:sumsq}
We have 
\beq
\sum\limits_{W\in \Irrep_{e,f}(G,F)}\frac{(\dim W)^2}{|\Pi_0^{tF}|^2}=q^{2d_e},
\eeq
where in the summation $\<t\>\in H^1(F,G)$ denotes the $F$-twisted conjugacy class in $G$ on which $W_{\loc}$ is supported. In particular if $G$ is connected solvable, then 
\beq
\sum\limits_{W\in \Irrep_{e,f}(G_0(\Fq))}{(\dim W)^2}=q^{2d_e}.
\eeq
\ecor
\bpf
This follows from Theorem \ref{t:relheis}(ii) and the fact that $\sum\limits_{W\in \Irrep_{e,f}(G,F)}\dim^+(M_W)^2=\dim \M_{G,f}$.
\epf

\brk\label{r:nominqi}
We know that the $F$-stable indecomposable minimal quasi-idempotents are parametrized by the set $\f{L}_e(G)^F=(\C(T)(\Qlcl)/\Pi_0)^F$. We can prove, for example using the above, that 
\beq
|\f{L}_e(G)^F|=|(\C(T)(\Qlcl)/\Pi_0)^F|=\sum\limits_{\<t\>\in H^1(F,\Pi_0)}\frac{|T^{tF}|}{|\Pi_0^{tF}|}=\frac{1}{\Pi_0}\sum\limits_{t\in \Pi_0}|T^{tF}|.
\eeq
\erk

\subsection{Proof of Theorem \ref{t:relheis}}\label{s:pot:relheis}
In this section we will complete the proof of Theorem \ref{t:relheis} stated above. Using the definition of $\tr_F$ it is straightforward to check that statements (i) and (ii) of the Theorem are equivalent. Hence we will prove below statements (ii) and (iii).

\subsubsection{The connected case}\label{s:tcc}
Let us first consider the case that $G$ is a connected solvable group.
\begin{proof}[Proof of Thm. \ref{t:relheis}(ii) in the connected case]
Let $W$ be an irreducible representation of $G_0(\Fq)$ lying in the $\f{L}$-packet $\Irrep_{e,f}(G,F)$ under consideration. We want to prove that (in the connected case)
\beq\label{e:tpt1}
\FPdim(M_W)^2=\frac{|K_\N|}{q^{2d_e}}\dim(W)^2.
\eeq
In the connected case, we have canonical equivalences $e\DG\cong e\D_U(U)\boxtimes \D_T(T)$ and $e\D_G(GF)\cong e\D_U(UF)\boxtimes \D_T^F(T)$. The indecomposable minimal quasi-idempotent $f$ in $e\DG$ must be of the form $\beel$ for some $\L\in \C(T)(\Qlcl)^F$. The simple object $M_W\in \M_{GF,f}$ corresponds to an object of the form $M\boxtimes e_\L\in e\D_U(UF)\boxtimes \D_T^F(T)$. Now $M_W$ is just a suitably shifted local system
supported on $GF$. The rank of this local system is equal to the dimension of the irreducible representation $W$ of $G_0(\Fq)$. Now we have noted above that $M_W=M\ast\beel$ for some $M\in \M_{UF,e}\subset e\D_U(UF)$. Now $M$ is a suitably shifted local system on $UF$ corresponding to a representation $V$ of $U_0(\Fq)$ corresponding to the minimal idempotent $e\in \D_U(U)$. We will now prove that $\dim W=\dim V$. For this we compute the rank of the stalk of $M_W$ at the point $F\in GF$. We have
$M_W(F)=(M\ast \beel)(F)=\int\limits_U M(uF)\otimes \beel(F^{-1}(u^{-1}))$. Now we can consider $\beel$ as an object of $e\D_U(G)$ and $\beel|_U$ as an object of $e\D_U(U)$. From \cite[\S8.5]{De2} we can see that $\beel|_U\cong e[2\dim T](\dim T)$ is supported only on $H$. Also we have $M\in e\D_U(UF)$, hence we see that $M_W(F)\cong\int\limits_H M(F)[2\dim H+2\dim T](\dim H+\dim T).$ Thus we see that rank of the stalk of $M_W$ at $F$ is equal to the rank of the stalk of $M$ at $F$. Hence we have $\dim W=\dim V$. Statement (\ref{e:tpt1}) now follows from \cite[\S6.2]{De3}.
\epf

\begin{proof}[Proof of Thm. \ref{t:relheis}(iii) in the connected case]
To prove (\ref{e:trftr+}) it suffices to consider the case where the objects $C,M$ are simple. Hence let us compare the two automorphisms 
$$\g_{C,\psi_C,M_W}:C\ast M_W \rar{} C\ast M_W \hbox{ and }\g^0_{C,\psi_C,M_W}:C\uast M_W \rar{} C\uast M_W,$$
where $C,M_W$ are (simple objects) as before where we now assume that $G$ is a connected solvable group. As before, we have the canonical equivalences $e\DG\cong e\D_U(U)\boxtimes \D_T(T)$ and $e\D_G(GF)\cong e\D_U(UF)\boxtimes \D_T^F(T)$ and we have $f\cong \beel$ for some $\L\in \C(T)(\Qlcl)^F$. Then $C$ corresponds to an object $C'\boxtimes e_\L\in  e\D_U(U)\boxtimes \D_T(T)$ with $C'\in \M_{U,e}$ and $M_W$ corresponds to an object $M\boxtimes e_\L\in e\D_U(UF)\boxtimes \D_T^F(T)$ with $M\in \M_{UF,e}$. Then $\g_{C,\psi_C,M_W}$ corresponds to the automorphism
\beq
\g^U_{C',\psi_{C'},M}\boxtimes \g^T_{e_\L,e_\L}:(C'\ast M)\boxtimes (F^*e_\L\ast e_\L)\rar{}(C'\ast M)\boxtimes (F^*e_\L\ast e_\L)
\eeq
and $\g^0_{C,\psi_C,M_W}$ corresponds to the automorphism
\beq
\g^U_{C',\psi_{C'},M}\boxtimes {\g^T_{e_\L,e_\L}}^0:(C'\ast M)\boxtimes (F^*e_\L\uast e_\L)\rar{}(C'\ast M)\boxtimes (F^*e_\L\uast e_\L).
\eeq
Hence we are essentially reduced to studying the case of the torus. Let us look at the stalks of the automorphisms $\g^T_{e_\L,e_\L}$ and ${\g^T_{e_\L,e_\L}}^0$ at the point $1\in T$. Then using the explicit description of this stalk in the spirit of \cite[\S3]{De3} we obtain $\g^T_{e_\L,e_\L}(1):\int\limits_Te_\L{(F(t))}\boxtimes e_\L(t^{-1})\rar{}\int\limits_Te_\L(t)\boxtimes e_\L(F^{-1}(t^{-1}))$ and we can deduce that $\g^T_{e_\L,e_\L}\cong V_{T_0}\otimes {\g^T_{e_\L,e_\L}}^0.$ Hence we conclude that $\tr(\g^T_{e_\L,e_\L}(1))=\frac{|T_0(\Fq)|}{q^{\dim T}}\cdot \tr({\g^T_{e_\L,e_\L}}^0(1))$. Then using the definition of $\tr_F$, we obtain 
\beq\label{e:trg0}
\tr_F(\g_{C,\psi_C,M_W})=\frac{|T_0(\Fq)|}{q^{\dim T}}\cdot \tr_F(\g^0_{C,\psi_C,M_W}).
\eeq 
But now using Theorem \ref{t:relheis}(i), we obtain that 
\beq\tr_F(\g_{C,\psi_C,M_W})=(-1)^{2d_e}\frac{q^{d_e}}{q^{\dim G}\cdot \sqrt{\dim \M_{G,f}}}\cdot \tr_{F,f}^+(\g^0_{C,\psi_C,M_W})
\eeq
which proves (\ref{e:trftr+}).
Finally comparing the above with (\ref{e:defts}) we complete the proof of Theorem \ref{t:relheis}(iii) in the connected case.
\epf

\subsubsection{An intermediate case}\label{s:aic}
We now consider a general neutrally solvable group $G$, but we only consider (indecomposable) minimal quasi-idempotents $f$ of the form $\beel$ with $\L\in (\C(T)(\Qlcl)^{\Pi_0})^F$.

\begin{proof}[Proof of Thm. \ref{t:relheis}(ii) in the intermediate case]
Let $W\in \Irrep_{e,\beel}(G,F)$. For simplicity of notation, let us assume without loss of generality that $W$ is an irreducible representation of the pure inner form $G_0(\Fq)$. The object $M_W$ is supported on the $G$-conjugation orbit $\O F$ of $F\in GF$. Recall from Proposition \ref{p:beeltrsubcat} that we have the braided $\Pi_0$-crossed fusion category $\tMbeel\subset {}^{\beel}\D_{G^\circ}(G)$ and we have $\M_{G,\beel}\cong (\tMbeel)^{\Pi_0}$. Similarly we can construct the $\tMbeel$-module category $\tM_{GF,\beel}\subset { }^{\beel}\D_{G^\circ}(GF)$ (cf. \cite[\S6.2]{De3}) and prove that $\M_{GF,\beel}\cong (\tM_{GF,\beel})^{\Pi_0}$. Let us consider the object of $\tM_{GF,\beel}$ underlying $M_W$ and also restrict $W$ to the connected component $G^\circ_0(\Fq)$. The same argument as in \cite[\S6.2, in particular (96)]{De3} can now be used to complete the proof of the statement (ii) in the intermediate case by reducing it to the connected case.
\epf

\begin{proof}[Proof of Thm. \ref{t:relheis}(iii) in the intermediate case]
Let $C$ be any object of $\M_{G,\beel}\cong (\tM_{G,\beel})^{\Pi_0}$ and let $M$ be any object of $\M_{GF,\beel}\cong (\tM_{GF,\beel})^{\Pi_0}$. We can represent $C$ by a pair $\left(\bigoplus\limits_{g\in \Pi_0}C_g, \phi_C\right)$, where for each $g\in\Pi_0$, $C_g\in \tMbeel$ lies in the graded component labelled by $g$ and where $\phi_C$ denotes the $\Pi_0$-equivariant structure. Similarly we can represent $M$ by a pair $\left(\bigoplus\limits_{h\in \Pi_0}M_{hF}, \phi_M\right)$ where $M_{hF}\in \tM_{GF,\beel}$ lies in the graded component labelled by $hF\in \Pi_0F$. We have to compare the two automorphisms $\g_{C,\psi,M}:C\ast M\rar{}C\ast M$ and $\g^0_{C,\psi,M}:C\uast M\rar{}C\uast M$. For each $tF\in \Pi_0F$, let us compare the respective components with grading $tF$, namely $\g_{C,\psi,M}(tF):(C\ast M)_{tF}\rar{}(C\ast M)_{tF}$ and $\g^0_{C,\psi,M}(tF):(C\uast M)_{tF}\rar{}(C\uast M)_{tF}$. Now we have $(C\ast M)_{tF}=\bigoplus\limits_{h_1h_2=t}C_{h_1}\ast M_{h_2F}$. Now we have $C_{h_1}\in \tM_{G^\circ h_1,\beel}\subset {}^{\beel}\D_{G^\circ}(G^\circ h_1)\cong e\D_U(Uh_1)\boxtimes {}^{e_\L}\D_T^{h_1}(T)$ and $M_{h_2F}\in \tM_{G^\circ h_2F,\beel}\subset {}^{\beel}\D_{G^\circ}(G^\circ h_2F)\cong e\D_U(Uh_2F)\boxtimes {}^{e_\L}\D_T^{h_2F}(T)$, where by a slight abuse of notation we have used the letter $h_i$ to denote both $h_i\in \Pi_0$ as well as its lift in $G$. Then using a similar argument as in the connected case, we deduce that $\g_{C,\psi,M}(tF)\cong V_{T_0^t}\otimes \g_{C,\psi,M}^0(tF)$ and hence that
\beq
\tr^{G^\circ}_{tF}(\g_{C,\psi,M}(tF))=\frac{|T_0^t(\Fq)|}{q^{\dim T}}\cdot \tr^{G^\circ}_{tF}(\g^0_{C,\psi,M}(tF)),
\eeq where we consider $\g_{C,\psi,M}(tF)$ and $\g_{C,\psi,M}(tF)$ as morphisms in $\D_{G^\circ}(G^\circ tF)$ and where $\tr^{G^\circ}_{tF}$ denotes the trace in this category. Then using Theorem \ref{t:relheis}(i) in the connected case, we obtain that
\beq
\tr^{G^\circ}_{tF}(\g_{C,\psi,M}(tF))= (-1)^{2d_e}\frac{q^{d_e}}{q^{\dim G}\cdot \sqrt{\dim \M_{G^\circ, \beel}}}\cdot \tr_{tF,f}^+(\g^0_{C,\psi,M}(tF)).
\eeq
Now from the definition of the trace $\tr_F=\tr^G_F$ on $\D_G(GF)$ (see also \cite[\S6.2]{De3}), we obtain that for any endomorphism $\alpha$ in $\D_G(GF)$ we have
$$\tr^G_{F}(\alpha)=\frac{1}{|\Pi_0|}\sum\limits_{t\in \Pi_0}\tr^{G^\circ}_{tF}(\alpha_{tF}).$$
Moreover, it is clear that for any endomorphism $\alpha$ in $\M_{GF,\beel}\cong (\tM_{GF,\beel})^{\Pi_0}$, we have
$$\tr^{G,+}_{F,\beel}(\alpha)=\sum\limits_{t\in \Pi_0}\tr^{G^\circ,+}_{tF,\beel}(\alpha_{tF}).$$ Hence this completes the proof of Theorem \ref{t:relheis}(iii) in the intermediate case since we have $\dim \Mbeel=|\Pi_0|^2\cdot\dim \M_{G^\circ,\beel}$.
\epf

\subsubsection{The general case}\label{s:tgc}
We now consider the general case, namely $G$ is any neutrally solvable algebraic group and $f\cong f_\L\in e\DG$ is a general $F$-stable indecomposable minimal quasi-idempotent where $\<\L\>\in (\C(T)(\Qlcl)/\Pi_0)^F$. As in \S\ref{s:minqiedgg}, let $\Pi_0'\subset \Pi_0$ denote the stabilizer of $\L$ and let $G^\circ\subset G'\subset G$ be the corresponding subgroup. Then $\beel\in e\DGp$ is a minimal quasi-idempotent for $G'$ and we can apply the results of \S\ref{s:aic} to this case. Note that since the $\Pi_0$-orbit of $\L$ is $F$-stable, we must have $(gF)^*\L\cong \L$ for some pure inner form $gF:G\rar{}G$. Hence without loss of generality we may assume that $F^*\L\cong \L$. With this assumption, the subgroup $G'$ is $F$-stable.

Recall from $\S\ref{s:minqiedgg}$ that we have equivalences
\beq\label{e:indg}
\indg:{}^{\beel}\DGp^\Delta\rar{\cong} {}^{f_\L}\DG^\Delta \hbox{ and hence } \indg:\M_{G',\beel}\rar{\cong}\M_{G,f}.
\eeq
Similarly, we can prove that (under the assumption that $F^*\L\cong \L$) we have equivalences
\beq\label{e:indgf}
\indg:{}^{\beel}\D_{G'}(G'F)\rar{\cong} {}^{f_\L}\D_G(GF) \hbox{ and } \indg:\M_{G'F,\beel}\rar{\cong}\M_{GF,f}.
\eeq
Concretely, let $W'\in \Irrep_{e,\beel}(G',F)$ be an irreducible representation of an inner form ${G'_0}^h(\Fq)$ for some $h\in G'$. Let $W=\ind^{G_0^h(\Fq)}_{{G'_0}^h(\Fq)}$. Then by (\ref{e:indgf}), $W\in \Irrep_{e,f_\L}(G,F)$ and we have $\M_{G'F,\beel}\ni M_{W'}\longmapsto M_W\in \M_{GF,f_\L}$. Thus we have a bijection $\Irrep_{e,\beel}(G',F)\cong \Irrep_{e,f_\L}(G,F)$.

\begin{proof}[Proof of Thm. \ref{t:relheis}(ii) in the general case]
Using Theorem \ref{t:relheis} for $W'$, we have 
\beq
\dim W' = \frac{\dim^+(M_{W'})}{\sqrt{\dim \M_{G',\beel}}}\cdot |{\Pi_0'}^{hF}|\cdot q^{d_e}.
\eeq
Now $\dim W=\frac{|G^h_0(\Fq)|}{|{G'_0}^h(\Fq)|}\cdot\dim W'=\frac{|\Pi^{hF}_0|}{|{\Pi'_0}^{hF}|}\cdot\dim W'$, $\dim^+(M_{W'})=\dim^+(M_W)$ and $\M_{G',\beel}\cong \M_{G,f_\L}$. Hence we obtain that
\beq
\dim W = \frac{\dim^+(M_{W})}{\sqrt{\dim \M_{G,f_\L}}}\cdot |{\Pi_0}^{hF}|\cdot q^{d_e}
\eeq as desired.
\epf

\begin{proof}[Proof of Thm. \ref{t:relheis}(iii) in the general case]
By \S\ref{s:aic}, statement (iii) holds for $C'\in \M_{G',\beel}$, $\psi':F^*C'{\cong} C'$ and $M'\in \M_{G'F,\beel}$. Let $(C,\psi)=\indg (C',\psi')$ and $M=\indg M'$. Then we see that $\indg(\g_{C',\psi',M'})=\g_{C,\psi,M}$ and that $\indg(\g^0_{C',\psi',M'})=\g^0_{C,\psi,M}$. The statement (iii) in the general case then follows using (\ref{e:indg}) and (\ref{e:indgf}). 
\epf

\subsection{Shintani descent in the Heisenberg case}\label{s:sdthc}
We will now study Shintani descent in the Heisenberg case and prove Theorem \ref{t:main3} in the Heisenberg case. Recall that for each positive integer $m$, we have defined the $m$-th Shintani descent map (well defined up to scaling by $m$-th roots of unity, cf. \cite{De4})
\beq
\Sh_m:\Irrep(G,F^m)^F\hookrightarrow\FunG.
\eeq
We will continue to study the Heisenberg case and we will use all our previous notations and conventions from \S\ref{s:icahi}. In particular we have a Heisenberg admissible pair $(H_0,\L_0)$ defined over $\Fq$ and the corresponding idempotent $e_0\in \DGn$. In this section, we are interested in the restriction of the Shintani descent map to the subset $\Irrep_{e}(G,F^m)^F\subset \Irrep(G,F^m)^F$. Using the same argument as in \cite[Prop. 6.4]{De4}, we have
\beq
\Sh_m(\Irrep_{e}(G,F^m)^F)\subset \T_{e_0}\FunG
\eeq
and that the image forms an orthonormal basis of $\T_{e_0}\FunG$ which we call the $m$-th Shintani basis of $\T_{e_0}\FunG$. Now the set of irreducible characters 
\beq\label{e:lpdecf}
\Irrepe=\coprod\limits_{f\in \f{L}_e(G)^F}\Irrep_{e,f}(G,F)
\eeq
associated with $e_0$ is also an orthonormal basis of $\T_{e_0}\FunG$. Similar to (\ref{e:lpdecf}), we also obtain the $\f{L}$-packet decompositions
\beq\label{e:lpdecfm}
\Irrep_{e}(G,F^m)=\coprod\limits_{f\in \f{L}_e(G)^{F^m}}\Irrep_{e,f}(G,F^m),
\eeq
\beq\label{e:lpdecfmf}
\Irrep_{e}(G,F^m)^F=\coprod\limits_{f\in \f{L}_e(G)^F}\Irrep_{e,f}(G,F^m)^F.
\eeq
We will prove that Shintani descent respects the $\f{L}$-packet decompositions above and will then prove the analogue of Theorem \ref{t:main3} in the Heisenberg case.

\subsubsection{Some notation and results related to Shintani descent}
Let us now recall some notation and results from \cite{De4}. For any $m_1,m_2\in \Z$, we let $\D_{G}(GF^{m_1})^{F^{m_2}}$ denote the category consisting of pairs $(M,\psi)$ where $M\in \D_G(GF^{m_1})$ and $\psi:{F^{m_2}}^*M\rar{\cong}M$. We have a functor $\eta_m$ defined as the composition (cf. \cite[Lem. 5.1]{De4}):
\beq
\eta_m:\D_{G}(GF)\rar{}\D_G(GF)^{\id}\rar{\cong} \D_G(GF)^{F^m}.
\eeq
For $L\in \D_G(GF^m)^F$ and $M\in \D_{G}(GF)^{F^m}$, we have the automorphism (cf. \cite[\S5.2, (70)]{De4})
\beq\label{e:defzeta}
\zeta_{L,M}:L\ast M\rar{}M\ast F^*L\rar{}F^*L\ast {F^m}^*M\rar{}L\ast M.
\eeq

For $W\in \Irrep(G,F^m)^F$ and $V\in \Irrep(G,F)$, let $(W_{loc},\psi_W)\in \D_G(GF^m)^F$ and $V_{loc}\in \D_G(GF)$ be the associated objects, where $\psi_W:F^*W_{loc}\rar{\cong}W_{loc}$ is chosen according to {\it op cit, \S2.6}. Then the inner product between the Shintani descent $\Sh_m(W)$ and the character $\chi_V$ is described by \cite[Cor. 5.5]{De4} as below:
\beq\label{e:innprodsd}
\<\Sh_m(W),\chi_V\>=\tr_{F^{m+1}}\left(\zeta_{W_{loc},\psi_W, \eta_m(V_{loc})}:W_{loc}\ast V_{loc}\rar{}W_{loc}\ast V_{loc}\right).
\eeq

In \cite[\S4]{De4} we have defined the twists $\theta^{F^m}$ in the categories $\D_G(GF^m)$ for each $m\in \Z$, i.e., we have natural isomorphisms (where $F$ is the inverse functor to $F^*$):
\beq
\theta^{F^m}_M:M\rar{\cong} F^m(M) \hbox{ for each } M\in \M_{GF^m}
\eeq
satisfying certain compatibility relations with the crossed braidings (cf. \cite[Lem. 4.2.]{De4}).

Also for each $M\in \D_{G}(GF)$ we have the automorphism 
\beq\label{e:defnu}
\nu_{\theta^F_M}:M\ast M \xto{\beta^{-1}_{M,F^*M}} M\ast F^*M\xto{\id_M\ast F^*(\theta^F_M)}M\ast M \hbox{ in } \D_{G}(GF^2)
\eeq
and we have (cf. \cite[Thm. 4.3.]{De4})
\beq\label{e:de4trnu}
\tr_F{(\id_M)}=\tr_{F^2}(\nu_{\theta^F_M}).
\eeq

\subsubsection{Shintani descent and $\f{L}$-packet decomposition in the Heisenberg case}\label{s:sdlpdthc}
Let $V\in \Irrep_{e,f}(G,F)$ for some $f\in \f{L}_e(G)^F$. Then by Definition \ref{d:deflpirre} we have $V_{loc}\in {}^f\D_G(GF)$. Let $W\in \Irrep_{e,f'}(G,F^m)^F$ for some $f'\in \f{L}_e(G)^F$. Similarly we have $W_{loc}\in {}^{f'}\D_G(GF^m)$. Hence we see that if $f,f'$ are non-isomorphic, then $W_{loc}\ast V_{loc}=0$. Hence using (\ref{e:innprodsd}) we see that in this case we must have $\<\Sh_m(W),\chi_V\>=0$. Hence we see that for each $f\in \f{L}_e(G)^F$, $\Sh_m(\Irrep_{e,f}(G,F^m)^F)$ lies in the subspace $\t{\T}_{f_0}\FunG$ of $\T_{e_0}\FunG$ spanned by the characters $\Irrep_{e,f}(G,F)$, i.e. Shintani descent respects the $\f{L}$-packet decomposition, at least in the Heisenberg case.

We will now describe the matrix $\t{\Sh}_m(F,f)$ relating the sets $\Sh_m(\Irrep_{e,f}(G,F^m)^F)$ and $\Irrep_{e,f}(G,F)$. The entries of this matrix are given by
\beq
\t{\Sh}_m(F,f)_{W,V}=\<\Sh_m(W),\chi_V\> = \tr_{F^{m+1}}\left(\zeta_{W_{loc},\psi_W, \eta_m(V_{loc})}\right)
\eeq
for $W\in \Irrep_{e,f}(G,F^m)^F$ and $V\in \Irrep_{e,f}(G,F)$.

Suppose that $M\in \M_{GF^{m+1},f}\subset {}^f\D_G(GF^{m+1})$ is an object supported on a single $G$-orbit $\<tF^{m+1}\>\subset GF^{m+1}$ and that $\alpha:M\to M$ is an endomorphism. By Theorem \ref{t:relheis}(i) we have
\beq
\tr_F(\alpha)=(-1)^{2d_e}\frac{q^{(m+1)d_e}}{q^{(m+1)\dim G}\cdot \sqrt{\dim \M_{G,f}}}\cdot \frac{q^{(m+1)\dim T}}{|T^{tF^{m+1}}|}\cdot \tr_{F^{m+1},f}^+(\alpha).
\eeq
Using the same argument as in the proof of Theorem \ref{t:relheis}(iii), we can prove the following more general result:
\blem\label{l:trzeta}
Let $(L,\psi_L)\in (\M_{GF^{m},f})^F\subset \D_G(GF^m)^F$ and let $(M,\psi_M)\in (\M_{GF,f})^{F^{m}}\subset \D_G(GF)^{F^m}$. Let $\zeta_{L,\psi_L,M,\psi_M}:L\ast M \rar{\cong} L\ast M$ be the automorphism as defined by (\ref{e:defzeta}). Then
\beq
 \tr_{F^{m+1}}(\zeta_{L,\psi_L,M,\psi_M})=\frac{(-1)^{2d_e}}{\sqrt{\dim \M_{G,f}}}\cdot\left(\frac{q^{d_e}}{q^{\dim G}}\right)^{m+1}\cdot\tr_{F^{m+1},f}^+(\zeta^0_{L,\psi_L,M,\psi_M}).
\eeq
\elem
We will use this result in \S\ref{s:acthc} describe the matrices $\t{{\Sh}}_m(F,f)$ in terms of the Shintani matrices as defined in \cite[\S4]{De5} in the setting of modular categories.

\subsubsection{Twists in the categories $\M_{GF^m,f}$}
Previously, we have described the twists $\theta^{F^m}$ in the categories $\D_G(GF^m)$. Now for $f\in \f{L}_e(G)^F$ we will define new twists $\theta^{F^m,f,+}$ in the invertible $\M_{G,f}$-module categories $\M_{GF^m,f}\subset \D_G(GF^m)$. Note that we have a braided $\Z$-crossed rigid monoidal category (under truncated convolution)
\beq
\M_{G\<F\>,f}:=\bigoplus\limits_{m\in \Z} \M_{GF^m,f}\subset \D_G(G\<F\>),
\eeq
where $G\<F\>$ denotes the semidirect product of $G$ and $\Z=\<F\>$. It can be equipped with a unique positive spherical structure and the positive traces $\tr_{F^m,f}^+$ on $\M_{GF^m,f}$ correspond to this spherical structure. This spherical structure also gives rise to twists $\theta^{F^m,f,+}$ on the categories $\M_{GF^m,f}$ as follows:\\
For each $M\in \M_{GF^m,f}$, define the isomorphism $\theta^{F^m,f,+}_M$ as the composition
\beq
\theta_{M}^{F^m,f,+}:M\xto{\id_M\uast \c_M}M\uast M\uast M^*\xto{\beta_{M,M}}F(M)\uast M \uast M^*\xto{\id_{F(M)}\uast \e_{M^*}}F(M),
\eeq
where the positive spherical structure has been used implicitly in the last map of the composition. 

\blem\label{l:cfe}
There exists a number $c_{F,f}^+\in \Qlcl^\times$ such that for each $m\in \Z$ and $M\in \M_{GF^m,f}$ we have
\beq\label{e:defcf}
\theta^{F^m}_M=(c_{F,f}^+)^m\cdot \theta^{F^m,f,+}_M.
\eeq
\elem
\bpf
The twists $\theta^{F^m,f,+}$ are defined using a spherical structure on $\M_{G\<F\>,f}$ and hence satisfy properties (i),(ii) from \cite[\S4.1]{De5}. We can also check that the twists $\theta^{F^m}$ also satisfy these properties. For example for $a,b\in \Z$, $L\in \M_{GF^a,f}\subset \D_G(GF^a)$ and $M\in \M_{GF^b,f}\subset \D_G(GF^b)$ we have from \cite[Lem. 4.2.]{De4} that $\theta^{F^{a+b}}_{L\ast M}$ equals the composition
\beq
L\ast M\xto{\beta_{L,M}} F^a(M)\ast L\xto{\beta_{F^a(M),L}} F^b(L)\ast F^a(M) \xto{\theta^{F^a}_{F^b(L)}\ast \theta^{F^b}_{F^a(M)}}F^{a+b}(L)\ast F^{a+b}(M)\to F^{a+b}(L\ast M).
\eeq
Truncating the above with respect to the canonical $t$-structure, we see that $(\theta^{F^{a+b}}_{L\ast M})^0=\theta^{F^{a+b}}_{L\uast M}$ equals the composition
\beq
L\uast M\xto{\beta^0_{L,M}} F^a(M)\uast L\xto{\beta^0_{F^a(M),L}} F^b(L)\uast F^a(M) \xto{\theta^{F^a}_{F^b(L)}\uast \theta^{F^b}_{F^a(M)}}F^{a+b}(L)\uast F^{a+b}(M)\to F^{a+b}(L\uast M).
\eeq
Hence the twists $\theta^{F^m}$ correspond to a pivotal structure on $\M_{G\<F\>,f}$ (cf. \cite[\S4.1]{De4}). The lemma now follows using the same argument as in \cite[Prop. 7.2.]{De4}.
\epf

Next, let us compute $c_{F,f}^+$. We will prove the following result in \S\ref{s:pl:trnu}:
\blem\label{l:trnu}
Let $M\in \M_{GF,f}$ be an object supported on the $G$-orbit $\<tF\>\subset GF$. Then
\beq\label{e:trnu}
\tr_{F^2}(\nu_{\theta^F_M})=\frac{(-1)^{2d_e}}{\sqrt{\dim \M_{G,f}}}\cdot\frac{q^{2d_e}}{q^{2\dim G}}\cdot\frac{q^{\dim T}}{|T^t_0(\Fq)|}\tr^+_{F^2,f}(\nu^0_{\theta_M^F})
\eeq
where $\nu_{\theta^F_M}:M\ast M\rar{}M\ast M$ is defined by (\ref{e:defnu}) and where $\nu_{\theta^F_M}^0:M\uast M\rar{} M\uast M$ denotes its truncation. We have 
\beq\label{e:cfval}
c_{F,f}^+=\frac{q^{\dim G}}{q^{d_e}}.
\eeq
\elem 

\subsubsection{Proof of Lemma \ref{l:trnu}}\label{s:pl:trnu}
First let us prove that equations (\ref{e:trnu}) and (\ref{e:cfval}) are in fact equivalent. For each $M\in \M_{GF,f}\subset \D_G(GF)$, we have the automorphism $\nu_{\theta^F_M}$  defined by (\ref{e:defnu}). Similarly, we can define the automorphism $\nu_{\theta^{F,f,+}_M}$ in ${}^f\D_G(GF)$ as the composition
\beq\label{e:defnuf}
\nu_{\theta^{F,f,+}_M}:M\ast M \xto{\beta^{-1}_{M,F^*M}} M\ast F^*M\xto{\id_M\ast F^*(\theta^{F,f,+}_M)}M\ast M
\eeq
and its truncation
\beq\label{e:defnuftr}
\nu^0_{\theta^{F,f,+}_M}:M\uast M \xto{\left(\beta^{-1}_{M,F^*M}\right)^0} M\uast F^*M\xto{\id_M\uast F^*(\theta^{F,f,+}_M)}M\uast M.
\eeq
Using the fact that the twists $\theta^{F^m,f,+}$ and the traces $\tr_{F^m,f}^+$ in $\M_{GF^m,f}$ come from a spherical structure, we see that (cf. \cite[\S7.1]{De4})
\beq\label{e:idmnuf+}
\tr_{F,f}^+(\id_M)=\tr_{F^2,f}^+(\nu^0_{\theta^{F,f,+}_M}).
\eeq
Then from (\ref{e:defnu}), (\ref{e:defnuf}), \ref{e:defnuftr} and (\ref{e:defcf}) we see that 
\beq\label{e:nucfnu}
\nu_{\theta^F_M} = c_{F,f}^+\cdot\nu_{\theta^{F,f,+}_M} \hbox{ and } \nu_{\theta^F_M}^0 = c_{F,f}^+\cdot\nu_{\theta^{F,f,+}_M}^0.
\eeq

 By (\ref{e:de4trnu}) we have 
\beq
\tr_{F^2}(\nu_{\theta^F_M})=\tr_F(\id_M)
\eeq
\beq
=\frac{(-1)^{2d_e}}{\sqrt{\dim \M_{G,f}}}\cdot \frac{q^{d_e}}{q^{\dim G}}\cdot \frac{q^{\dim T}}{|T^t_0(\Fq)|}\cdot\tr_{F,f}^+(\id_M) \hbox{ }\hbox{ }\cdots\hbox{by Thm. \ref{t:relheis}(i)}
\eeq
\beq
=\frac{(-1)^{2d_e}}{\sqrt{\dim \M_{G,f}}}\cdot \frac{q^{d_e}}{q^{\dim G}}\cdot \frac{q^{\dim T}}{|T^t_0(\Fq)|}\cdot\frac{\tr_{F^2,f}^+(\nu_{\theta_{M}^{F}}^0)}{c_{F,f}^+} \hbox{ }\hbox{ }\cdots\hbox{by (\ref{e:idmnuf+}) and (\ref{e:nucfnu})}.
\eeq
Moreover, we know that $\tr_{F,f}^+(\id_M)\neq 0$. Hence we can now see that (\ref{e:trnu}) and (\ref{e:cfval}) are in fact equivalent.

Now let us prove that the lemma holds for tori. Hence suppose that $G=T$, a torus.  In this case, the Heisenberg idempotent $e$ must be the unit $\delta_1\in \D_T(T)$. Let $\L\in \C(T)(\Qlcl)^F$ be such that $f\cong e_\L\in \D_T(T)$. Note that in our situation we have
\beq\label{e:torcase}
d_e=0 \hbox{ and } \M_{T,e_\L}\cong \Vec.
\eeq 
Since we know that  (\ref{e:trnu}) and (\ref{e:cfval}) are equivalent, let us prove (\ref{e:trnu}) for the  torus $T$. Towards this end, we prove:
\blem\label{l:trf2}
Suppose that $G=T$ is a torus and that $f\cong e_\L\in \D_T(T)$ is an $F$-stable minimal quasi-idempotent as before. Then for $M\in \M_{TF,e_\L}$ we have
\beq
\tr_{F^2}(\nu_{\theta_M^F})=\frac{|\{s\in T|sF(s)=1\}|}{q^{\dim T}}\cdot\tr_{F^2}(\nu^0_{\theta_{M}^F}).
\eeq
\elem
\bpf
We know that $\M_{TF,e_\L}\cong \Vec$ and hence let us suppose that $M$ is the (unique) simple object, namely $e_\L$ translated by $F$. The proof is similar to the proof of (\ref{e:trg0}). Proceeding in the same way, let us look at the stalk (where we identify $M$ with the object $e_\L\in \D^F_G(G)\cong \D_G(GF)$, cf. \cite[\S4.2]{De4})
\beq
\nu_{\theta^F_M}(1): \int\limits_{h_1F(h_2)=1}e_{\L}({h_1})\otimes e_\L({h_2}) \rar{\cong} \int\limits_{h_1F(h_2)=1}e_\L({h_2})\otimes e_\L(F^{-1}(h_2^{-1}h_1)h_2). 
\eeq
Then using the same argument as in the proof of (\ref{e:trg0}) and \cite[\S4.2]{De4} we complete the proof of the lemma.
\epf

Note that we have $\{s\in T|sF(s)=1\}\subset T^{F^2}$. In fact, we have the short exact sequence
\beq
0\to \{s\in T|sF(s)=1\}\subset T^{F^2}\xto{s\mapsto sF(s)} T^F\to 0
\eeq
corresponding to the `norm' map. Hence we have 
\beq\label{e:tf2}
|T^{F^2}|=|T^F|\cdot |\{s\in T|sF(s)=1\}|.
\eeq
Now by Theorem \ref{t:relheis}(i) applied to the Frobenius $F^2:T\rar{}T$ and (\ref{e:torcase}) we have
\beq\label{e:trf2+}
\tr_{F^2}(\nu^0_{\theta_{M}^F})=\frac{1}{|T^{F^2}|}\cdot{\tr_{F^2,e_\L}^+(\nu_{\theta_{M}^{F}}^0)}.
\eeq 
Hence using Lemma \ref{l:trf2}, (\ref{e:tf2}) and (\ref{e:trf2+}) we get
\beq
\tr_{F^2}(\nu_{\theta_M^F})=\frac{1}{q^{\dim T}\cdot |T^{F}|}\cdot{\tr_{F^2,e_\L}^+(\nu_{\theta_{M}^{F}}^0)}.
\eeq
This completes the proof of (\ref{e:trnu}) and hence of Lemma \ref{l:trnu} in the case of the torus. In particular for a torus we have $c^+_{T,F,e_\L}=q^{\dim T}$.

Next let us consider the case of a connected solvable group $G=TU$. Let $f\cong \beel$, where $\L\in \C(T)(\Qlcl)^F$. In this case we have
\beq\label{e:edgf}
{}^{\beel}\D_G(G)\cong e\D_U(U)\boxtimes {}^{e_\L}\D_T(T) \hbox{ and } {}^{\beel}\D_G(GF)\cong e\D_U(UF)\boxtimes {}^{e_\L}\D^F_T(T).
\eeq
Since we know that Lemma \ref{l:trnu} holds for unipotent groups by \cite{De4} and since we have proved it for tori above, it follows that the lemma holds for connected solvable $G$ and in particular we have $c^+_{G,F,\beel}=\frac{q^{\dim G}}{q^{d_e}}$.

Finally using the connected case, we can also prove the result for all neutrally solvable groups.

\subsubsection{Shintani matrices and almost characters in the Heisenberg case}\label{s:acthc}
Recall that we have a (positive integral) modular category $\M_{G,f}$ and an invertible $\M_{G,f}$-module category $\M_{GF,f}$ equipped with the (positive) $\M_{G,f}$-module trace $\tr_{F,f}^+$. In this abstract categorical setting, for each positive integer $m$, we have defined the $m$-th Shintani matrix in \cite[\S4]{De5}. Let us denote the matrix we obtain this way by $\Sh_m^+(F,f)$. On the other hand, we have the matrix $\t{\Sh}_m(F,f)$ from \S\ref{s:sdlpdthc} which relates the two sets $\Sh_m(\Irrep_{e,f}(G,F^m)^F)$ and $\Irrep_{e,f}(G,F)$ and whose entries are given by
\beq\label{e:deftsh}
\t{\Sh}_m(F,f)_{W,V}=\<\Sh_m(W),\chi_V\> = \tr_{F^{m+1}}\left(\zeta_{W_{loc},\psi_W, \eta_m(V_{loc})}\right)=\tr_{F^{m+1}}\left(\zeta_{M_W,\psi_{M_W}, M_V,\psi_{M_V}}\right)
\eeq
for $W\in \Irrep_{e,f}(G,F^m)^F$ and $V\in \Irrep_{e,f}(G,F)$, where 
$$\psi_{M_W}:=\psi_W[\dim G+\dim H+\dim T]:F^*M_W\rar\cong M_W$$
is chosen in such a way that the composition (cf. \cite[Rem. 4.1.]{De4})
\beq\label{e:psicondi}
{F^m}^*M_W\xto{{F^{(m-1)}}^*(\psi_{M_W})}{F^{(m-1)}}^*M_W\to\cdots F^*M_W\xto{\psi_{M_W}}M_W
\eeq equals ${F^m}^*\theta^{F^m}_{M_W}$ and where $\psi_{M_V}:{F^m}^*M_V\rar{\cong} M_V$ is defined using the twist $\theta^F$ in $\M_{GF,f}\subset \D_G(GF)$ (cf. \cite[\S5.1.]{De4}).

Let us recall the definition of the matrix ${\Sh}^+_m(F,f)$ from \cite[\S4.2.]{De5}. For this we must choose $$\psi'_{M_W}:F^*M_W\rar\cong M_W$$
 in such a way that the composition 
\beq
{F^m}^*M_W\xto{{F^{(m-1)}}^*(\psi'_{M_W})}{F^{(m-1)}}^*M_W\to\cdots F^*M_W\xto{\psi'_{M_W}}M_W
\eeq 
equals ${F^m}^*\theta^{F^m,f,+}_{M_W}=\frac{{F^m}^*\theta^{F^m}_{M_W}}{\left(c^+_{F,f}\right)^m}$ and we define $\psi'_{M_V}:{F^m}^*M_V\rar{\cong} M_V$ using the twist $\theta^{F,f,+}=\frac{\theta^F}{c^+_{F,f}}$ in $\M_{GF,f}$. In particular, we may (and will) choose $\psi'_{M_W}=\frac{\psi_{M_W}}{c_{F,f}^+}$ where $\psi_{M_W}$ is chosen to be as before and we have $\psi'_{M_V}=\frac{\psi_{M_V}}{\left(c_{F,f}^+\right)^m}$.

Then the entries of the matrix $\Sh^+_m(F,f)$ are defined by
\beq
{\Sh}^+_m(F,f)_{W,V}:=\tr_{F^{m+1},f}^+\left(\zeta^0_{M_W,\psi'_{M_W}, M_V,\psi'_{M_V}}\right)=\frac{1}{(c_{F,f}^+)^{m+1}}\tr_{F^{m+1},f}^+\left(\zeta^0_{M_W,\psi_{M_W}, M_V,\psi_{M_V}}\right).
\eeq
Finally, by (\ref{e:deftsh}), Lemma \ref{l:trzeta} and Lemma \ref{l:trnu} we obtain:
\bthm\label{t:tshsh}
For each $W\in \Irrep_{e,f}(G,F^m)^F$ choose $\psi_{M_W}$ satisfying (\ref{e:psicondi}) and set $\psi'_{M_W}=\frac{\psi_{M_W}}{c^+_{F,f}}$. With these choices, we have
\beq
\t{\Sh}_m(F,f)=\frac{(-1)^{2d_e}}{\sqrt{\dim \M_{G,f}}}\cdot\Sh^+_m(F,f).
\eeq
\ethm
As a corollary we can now complete the proof of the Heisenberg case of Theorem \ref{t:main3}.
\bcor\label{c:tsh}
Let $e$ be an $F$-stable Heisenberg idempotent in $\DG$ and let $f$ be an $F$-stable  indecomposable minimal quasi-idempotent in $e\DG$. Then\\
(i) For each positive integer $m$, the image $\Sh_m(\Irrep_{e,f}(G,F^m))^F\subset \t\T_{f_0}\FunG$ is an orthonormal basis known as the $m$-th Shintani basis of $\t\T_{f_0}\FunG$.\\
(ii) There exists a positive integer $m_0$ such that for any positive integer $m$, the $m$-th Shintani basis of $\t\T_{f_0}\FunG$ depends (up to scalings by roots of unity) only on the residue of $m$ modulo $m_0$.\\
(iii) Let $m$ be a positive multiple of $m_0$. Then the $m$-th Shintani basis of $\t\T_{f_0}\FunG$ agrees with the basis $\left\{\frac{q^{\dim G}}{q^{d_e}}\cdot \T_{C,\psi_C}|C\in \CS_{e,f}(G)^F\right\}$ up to scaling by roots of unity.
\ecor
\bpf
By Theorem \ref{t:tshsh} and \cite[Prop. 4.7.]{De5}, the matrix $\t\Sh_m(F,f)$ which relates the sets $\Sh_m(\Irrep_{e,f}(G,F^m))^F$ and $\Irrep_{e,f}(G,F)$ is unitary. Hence statement (i) follows.

Using \cite[Thm. 4.10(ii)]{De5}, we see that there exists an $m_0$ such that the matrices $\Sh_m^+(F,f)$ only depend (up to scaling by roots of unity) on the residue of $m$ modulo $m_0$. This completes the proof of (ii).

By Theorem \ref{t:relheis}, the transition matrix between the sets $\left\{\frac{q^{\dim G}}{q^{d_e}}\cdot \T_{C,\psi_C}|C\in \CS_{e,f}(G)^F\right\}$ and $\Irrep_{e,f}(G,F)$ is given by $\frac{1}{\sqrt{\dim \M_{G,f}}}\cdot S^+(F,f)$, where $S^+(F,f)$ is the crossed S-matrix. Statement (iii) now follows from \cite[Thm. 4.10(iii)]{De5}.
\epf

This completes the proofs of the Heisenberg case versions of all the main results of this paper (Theorems \ref{t:main1}, \ref{t:main2} and \ref{t:main3}).

\section{Proof of the main results in the general case}\label{s:pmr}
In this section, we will use the results from the Heisenberg case to complete the proofs of all our main results in the general case. In fact, we will also state more precise versions of our main results.

\subsection{Definition of character sheaves and $\f{L}$-packets}\label{s:dcs}
In this section we will define the set $\CS(G)$ of all character sheaves on a neutrally solvable algebraic group $G$ over any algebraically closed field $\k$ of characteristic $p>0$. We will also describe an $\f{L}$-packet decomposition of the set $\CS(G)$ as well as the modular categories attached to such $\f{L}$-packets. 

As in \S\ref{s:i}, for a group $G$ as above, let $\hG$ denote the set of (isomorphism classes of) minimal idempotents in the braided monoidal category $\DG$. For each minimal idempotent $e\in \hG$ we will first define the set $\CS_e(G)$ of character sheaves in the braided monoidal category $e\DG\subset \DG$. By \cite[Thm. 2.28]{De2} for such an $e$, there exists an admissible pair $(H,\N)$ for $G$ which gives rise to the minimal idempotent $e$. Let $G'\subset G$ be the normalizer of $(H,\N)$ and let $e'_\N:=\N\otimes \can_H\in \DGp$ be the corresponding Heisenberg idempotent. Then according to the results {\it op cit.}, we have an equivalence of triangulated braided monoidal categories
\beq\label{e:indgedgp}
\indg:e'_\N\DGp \rar{\cong} e\DG \hbox{ with } e'_\N \mapsto e.
\eeq
\brk\label{r:ribbon}
As we have noted before, the triangulated braided monoidal category $e'_\N\DGp$ has the structure of a ribbon $\r$-category with the duality functor defined by $\f{D}^-_{G'}(\cdot)[2\dim H](\dim H)$ (cf. \cite[\S2.3]{De1}, \cite[Appendix A]{BD}). Transferring this structure along the above equivalence, we see that $e\DG$ must have the structure of a ribbon $\r$-category. On the other hand, $\DG$ is also a ribbon $\r$-category with duality functor $\f{D}^-_G$ and $e\in \DG$ is a locally closed idempotent by \cite[Thm. 2.28]{De2}. Hence by using the argument from \cite[Appendix A.6]{BD}, we can prove that $e\DG$ is a Grothendieck-Verdier category with dualizing object $\f{D}^-_Ge$. Hence, both $e$ (the unit object of $e\DG$) and $\f{D}^-_Ge$ are dualizing objects of $e\DG$. Hence by \cite[Rem. A.3]{BD}, the object $\f{D}^-_Ge$ must be invertible in $e\DG$.
\erk

Let $e\in \hG$ be a minimal idempotent with $(H,\N)$ an admissible pair giving rise to $e.$ With all notation as before, suppose that $T'$ is a maximal torus of $G'$. By Remark \ref{r:vtqi}, each indecomposable minimal quasi-idempotent in $e'_\N\DGp$ is a $V_{T'}$-quasi-idempotent. Hence it follows that each indecomposable minimal quasi-idempotent in $e\DG\cong e'_\N\DGp$ must also be a $V_{T'}$-quasi-idempotent. In particular, we see that the object $V_{T'}\in \Dmon$ is independent of the choice of the admissible pair $(H,\N)$ that gives rise to $e$ and hence $\dim T'$ is also independent of the choice of admissible pair. In particular, to each $e\in \hG$ we can associate an object $V_e\in \Dmon$ such that every indecomposable minimal quasi-idempotent in $e\DG$ is a $V_e$-quasi-idempotent. The object $V_e$ is supported in cohomological degrees $\{0,-1,\cdots,-\dim T'\}$.

For $e\in \hG$ as above, consider its stalk $e_1\in \D_G(1)$. Then we have identifications $e_1\cong \avg ({e'_\N}_1)\cong \avg(\Qlcl)[2\dim H](\dim H)$. Apply the forgetful functor and consider $e_1$ as an object of $D^b\Vec$. As such its cohomology is nonzero in degree $2(\dim(G/G')-\dim H)$ and vanishes in all higher degrees. In particular, we see that the integer $n_e:=\dim H-\dim(G/G')$ only depends on $e$ and not on the choice of the admissible pair $(H,\N)$ giving rise to $e$.

In view of these comments, let us make the following definitions:
\bdefn\label{d:miniinv}
Let $e\in \hG$ be a minimal idempotent. Then we have the associated object $V_e\in \Dmon$ as described above. Let $\tau_e\in \Z_{\geq 0}$ be such that $V_e$ is supported in cohomological degrees $\{0,\cdots,-\tau_e\}$. Define the functional dimension as $d_e:=\frac{\dim G - n_e -\tau_e}{2}$.
\edefn

\brk
If $e\in \hG$ comes from an admissible pair $(H,\N)$, we have $n_{e'}=\dim H$, $V_{e'}\cong V_e\cong V_{T'}$, $\tau_{e'}=\tau_e$ and $d_e=d_{e'}+\dim G-\dim G'$. Moreover, the definition of $d_e$ above agrees with the definition from Remark \ref{r:fundim} in the case of Heisenberg idempotents.
\erk

\bthm\label{t:modcatgen}
Let $G$ be a neutrally solvable group. Let $e\in \hG$ and let $f$ be any indecomposable minimal  quasi-idempotent in $e\DG$.\\ 
(i) Then ${}^fe\DG^\Delta={}^f\DG^\Delta$ is a $\Qlcl$-linear triangulated braided semigroupal category.  It has a unique $t$-structure, denoted here by $(\D^{\leq 0}, \D^{\geq 0})$, such that $\D^{\leq 0}\ast \D^{\leq 0}\subset \D^{\leq 0}$ but $\D^{\leq 0}\ast \D^{\leq 0}\nsubseteq \D^{\leq -1}$. \\
(ii) Let $\M_{G,f}^\Delta\subset {}^{f}\D_{G}(G)^\Delta$ denote the heart of the above $t$-structure. For each $i\in \Z, X\in {}^{f}\D_{G}(G)^\Delta$ let $X^i\in \M_{G,f}^\Delta$ denote the $i$-th cohomology of $X$ with respect to this $t$-structure. Then $(\M_{G,f}^\Delta,\underline{\ast})$ has the structure of a  braided semigroupal category, where 
\beq
\underline{\ast}:\M_{G,f}^\Delta\times \M_{G,f}^\Delta\rar{}\M_{G,f}^\Delta
\eeq denotes the truncated convolution functor defined by $X\uast Y:=(X\ast Y)^0$ for $X,Y\in \M_{G,f}$. \\
(iii) Let $\M_{G,f}\subset \M_{G,f}^\Delta$ denote the full subcategory formed by the semisimple objects. Then $\M_{G,f}={}^f\DG\cap \M_{G,f}^\Delta$ is closed under truncated convolution and $(\M_{G,f},\uast,f)$ has the structure of a non-degenerate braided fusion category. There is a natural spherical structure on $\M_{G,f}$, thus giving it the structure of a modular category.\\
(iv) The Frobenius-Perron dimension $\FPdim(\M_{G,f})$ is the square of an integer.\\
(v) Suppose that the ground field $\k$ is $\Fqcl$. Then the natural spherical structure on $\M_{G,f}$ is positive integral, i.e. categorical dimensions of all objects of $\M_{G,f}$ are positive integers.
\ethm
\bpf
Let $(H,\N)$ be an admissible pair (with normalizer $G'$) which gives rise to the minimal idempotent $e$ as above and let $e'_\N\in \DGp$ be the corresponding Heisenberg idempotent.  The theorem is now clear using Corollary \ref{c:modcat} for the Heisenberg idempotent $e'_\N\in \DGp$ and the equivalence (\ref{e:indgedgp}).
\epf

\bdefn\label{d:lp}
Let $\f{L}(G)$ denote the set of all indecomposable minimal quasi-idempotents in the triangulated braided monoidal category $\DG$. For $e\in \hG$, let $\f{L}_e(G)$ denote the set of all indecomposable minimal quasi-idempotents in the triangulated braided monoidal category $e\DG$.
\edefn

\bprop\label{p:lleg}
We have 
\beq
\f{L}(G)=\coprod\limits_{e\in \hG}{\f{L}_e(G)}.
\eeq
Each equivalence class of minimal quasi-idempotents in $\DG$ contains a unique (up to isomorphism) indecomposable object.
\eprop
\bpf
Let $f$ be any minimal quasi-idempotent. Then by \cite[\S5.5]{De2}, there exists a minimal idempotent $e\in\DG$ such that $e\ast f\neq 0.$ Then using the facts that $f$ is a minimal quasi-idempotent, $e$ is an idempotent and Lemma \ref{l:wv'}, we conclude that $e\ast f \cong f$, i.e. $f\in e\DG$. Both the statements of the proposition now follow (using Theorem \ref{t:minqiedgg}(iii$'$) and (\ref{e:indgedgp}) to prove the second statement).
\epf

Let us now define character sheaves on neutrally solvable groups.
\bdefn\label{d:csg}
Let $G$ be a neutrally solvable group. \\
(i) Let $f\in \f{L}(G)$ be an indecomposable minimal quasi-idempotent with the associated modular category $\M_{G,f}\subset {}^f\DG\subset {}^f\DG^\Delta\subset \DG$. We define $\CS_f(G)$ to be the (finite) set of isomorphism classes of simple objects of $\M_{G,f}$. The set $\CS_f(G)$ is known as the $\f{L}$-packet of character sheaves associated with $f\in \f{L}(G)$.\\
(ii) Let $e\in \hG$ be a minimal idempotent. The set of character sheaves on $G$ associated with $e$ is defined as the union $\CS_e(G):=\coprod\limits_{f\in \f{L}_e(G)}\CS_f(G)$.\\
(iii) The set of all character sheaves on $G$ is defined as the union $\CS(G):=\coprod\limits_{f\in \f{L}(G)}\CS_f(G)=\coprod\limits_{e\in \hG}\CS_e(G)$.
\edefn

\brk
Here we have implicitly used the fact that if $f,f'\in \f{L}(G)$ (resp. $e,e'\in \hG$) are non-isomorphic, then $f\ast f'=0$ (resp. $e\ast e'=0$).
\erk

With the following result, we see that the set $\CS(G)$ of character sheaves satisfies all the properties as desired in Theorem \ref{t:main1}:
\bprop\label{p:main1iv}
Let $\Phi:\DG\rar\cong \DG$ be any $\Qlcl$-linear triangulated braided monoidal auto-equivalence. Then $\Phi$ preserves the set $\CS(G)$ of character sheaves in $\DG$ as well as their idempotent and $\f{L}$-packet decompositions, i.e. $\Phi$ induces a permutation of the sets $\hG$, $\f{L}(G)$ and we have $\Phi(\CS_e(G))=\CS_{\Phi(e)}(G)$ and $\Phi(\leG)=\lL_{\Phi(e)}(G)$ for each $e\in \hG$ and $\Phi(\CS_f(G))=\CS_{\Phi(f)}(G)$ for each $f\in \f{L}(G)$. For each $f\in \lG$, the restriction of $\Phi$ to $\M_{G,f}\subset \DG$ induces an equivalence of modular categories:
\beq
\Phi:\M_{G,f}\rar{\cong} \M_{G,\Phi(f)}.
\eeq
\eprop
\bpf
The sets $\hG, \lG, \CS(G)$ as well as the modular categories $\M_{G,f}$ are all defined purely in terms of the $\Qlcl$-linear triangulated braided monoidal structure of $\DG$. Hence the proposition is obvious.
\epf

\subsection{Irreducible representations and admissible pairs}\label{s:irap}
For the remainder of this paper, we take the ground field to be $\Fpcl$ and we let $G$ be a neutrally solvable group over $\Fpcl$ equipped with an $\Fq$-Frobenius $F:G\rar{}G$. In the remainder of the paper, our goal is to study the irreducible characters of all pure inner forms $G_0^g(\Fq)$ using the theory of character sheaves.

We will begin by studying the set $\Irrep(G,F)$ in terms of admissible pairs. Let $\P_{\h{adm}}(G)$ denote the set of all admissible pairs $(H,\N)$ (where the multiplicative local system $\N$ is considered as a point of $H^*$). Then $G$ acts on $\Padm$ by conjugation and we let $[\Padm]$ denote the quotient for this action. By \cite[Thm. 2.28]{De2}, we have a surjective map
\beq\label{e:padmhg}
[\Padm] \onto \hG \hbox{ defined by } (H,\N)\mapsto \indg e'_\N.
\eeq

Now the Frobenius $F$ acts on the set $[\Padm]$. Consider the set of fixed points $[\Padm]^F$. We see that the geometric conjugacy class $\<(H,\N)\>$ is $F$-stable if and only if there exists a pure inner form $gF$ of $F$ such that the admissible pair $(H,\N)$ is $gF$-stable. In this case, the normalizer $G'$ must also be $gF$-stable. Such  a pair gives rise to an object $e'_{\N_0^g}\in \D_{{G'}^g_0}({G'}^g_0)$ and the induced object $e^g_0:=\ind_{{G'}^g_0}^{G^g_0}e'_{\N_0^g}\in \D_{G^g_0}(G^g_0)$. Now by \cite[\S4.4]{B} we have an identification $\D_{G^g_0}(G^g_0)\cong \D_{G_0}(G_0)$ and hence we obtain the idempotent $e_0\in \D_{G_0}(G_0)$. We also obtain the function $\T_{e_0}\in \FunG$ which will also be an idempotent. In other words, if the geometric conjugacy class of an admissible pair $(H,\N)$ is $F$-stable, then the corresponding minimal idempotent $e\in \hG$ comes from an idempotent $e_0\in \DGn$ by extension of scalars. Moreover, in this case, it is easy to check that the idempotent $\T^g_{e_0}\in \Fun(G^g_0(\Fq)/\sim)$ is nonzero. Indeed, using \cite[Prop. 4.12]{B} we can prove that $\T^g_{e_0}(1)$ is a positive rational number (see also \cite[Prop. 2.20, \S A.6]{B}). Hence the idempotent $\T_{e_0}\in \FunG$ obtained from an $F$-stable geometric conjugacy class of an admissible pair as above is nonzero. Moreover it is clear that $\T_{e_0}$ only depends on the geometric conjugacy class $\<(H,\N)\>\in [\Padm]^F$. In fact we have
\blem\label{l:tre0e}
Let $\<(H_1,\N_1)\>, \<(H_2,\N_2)\>\in [\Padm]^F$ and let $\T_{{e_1}_0}, \T_{{e_2}_0} \in \FunG$. Then either $\T_{{e_1}_0}\ast\T_{{e_2}_0}=0$ or $\T_{{e_1}_0}=\T_{{e_2}_0}$.
\elem
\bpf
Suppose that $\T_{{e_1}_0}\ast\T_{{e_2}_0}\neq 0$. Hence we must have ${e_1}_0\ast {e_2}_0\neq 0$ and hence $e_1\ast e_2\neq 0$. But both $e_1,e_2\in \hG$ are minimal idempotents. Hence we must have $e_1\cong e_2$. Also, $e_i\in \DG$ are simple, i.e. $\Hom(e_i,e_i)=\Qlcl$. Hence it follows that the associated idempotents $\T_{{e_1}_0},\T_{{e_2}_0}$ must be equal.
\epf

The map (\ref{e:padmhg}) induces a map 
\beq\label{e:padmfhgf}
[\Padm]^F\rar{}\hG^F.
\eeq
Let us denote the image of the above map by $\hG^F_{\h{adm}}$. This is the set of $F$-stable minimal idempotents in $\DG$ that can be obtained from an $F$-stable geometric conjugacy class of admissible pairs. The following lemma is clear:
\blem\label{l:uniqweil}
Given $e\in \hG^F_{\h{adm}}$, it has a unique Weil structure $\psi_e:F^*e\rar\cong e$ such that the associated trace function $\T_{e,\psi_e}\in \FunG$ is an idempotent. By a slight abuse of notation, we may often denote it by $\T_e$. Thus associated with each $e\in \hG^F_{\h{adm}}$ is a nonzero idempotent $\T_e\in \FunG$.
\elem

\bdefn
Suppose that $\<(H,\N)\>\in [\Padm]^F$ as above and let $\T_{e_0}\in \FunG$ be the corresponding idempotent with its `complex conjugate' being denoted $\bar\T_{e_0}$. Let $e\in \hG^F_{\h{adm}}$ be the corresponding minimal idempotent. Then we define 
\beq
\Irrep_{H,\N}(G,F)=\Irrepe:=\{W\in \Irrep(G,F)|\bar\T_{e_0}=\bar\T_e \hbox{ acts in $W$ as the identity}\}.
\eeq
By Lemmas \ref{l:tre0e}, \ref{l:uniqweil} the sets $\Irrep_{H,\N}(G,F)$ are all non-empty and are either equal or disjoint as we vary the pair $(H,\N)$. Moreover, for two such admissible pairs $(H_1,\N_1), (H_2,\N_2)$ we have $\Irrep_{H_1,\N_1}(G,F)=\Irrep_{H_2,\N_2}(G,F)$ if and only if the associated minimal idempotents in $\hG^F_{\h{adm}}$ are isomorphic.
\edefn

Below we give an alternative characterization of the set $\Irrep_{H,\N}(G,F)$ in the spirit of \cite[Def. 2.13]{B1}.
\bprop\label{p:altlp}
Let $\C\in [\Padm]^F$ be an $F$-stable geometric conjugacy class of admissible pairs, let $e$ be its image in $\hG^F_{\h{adm}}$ and let $\T_e\in \FunG$ be the corresponding idempotent. \\
(i) For $g\in G$, the idempotent $\T_e^g\in \Fun(G^g_0(\Fq)/\sim)$ is nonzero if and only if there exists a $gF$-stable admissible pair $(H,\N)$ in the geometric conjugacy class $\C$.\\
(ii) For $g\in G$, define $\Irrep_{\C}(G^g_0(\Fq))=\Irrep_{e}(G^g_0(\Fq)):=\Irrep(G^g_0(\Fq))\cap\Irrepe$. Let $W\in \Irrep(G^g_0(\Fq))$. Then $W\in \Irrep_{\C}(G^g_0(\Fq))$ if and only if $W$ occurs inside $\h{ind}_{H^g_0(\Fq)}^{G^g_0(\Fq)}\T_{\N_0^g}$ for some $gF$-stable admissible pair $(H,\N)\in \C$, where $\T_{\N_0^g}:H^g_0(\Fq)\to \Qlcl^\times$ is the multiplicative character associated with the $gF$-stable multiplicative local system $\N$. 
\eprop
\bpf
Let us prove statement (i). If an admissible pair $(H,\N)$ in $\C$ is $gF$-stable, then we have already seen that $\T_e^g\in \Fun(G^g_0(\Fq))$ is nonzero.  

In the other direction, suppose that $g\in G$ is such that $\T^g_e\neq 0$. Let us fix an admissible pair $(H,\N)\in \C\in [\Padm]^F$. Now there exists some pure inner form $tF$ of $F$ such that our chosen admissible pair $(H,\N)$ is $tF$-stable. For simplicity, let us replace our original Frobenius $F$ with $tF$ everywhere. With this choice, $(H,\N)$ is $F$-stable and hence its normalizer $G'$ is an $F$-stable subgroup of $G$. Now we have the Heisenberg idempotent $e'_0\in \D_{G'_0}(G'_0)$ and the corresponding minimal idempotent $e_0\in \DGn$ obtained by induction. Now we have assumed that $\T_{e_0}^g\in \Fun(G^g_0(\Fq)/\sim)$ is nonzero. Hence by \cite[Prop. 4.12]{B} we see that the $F$-twisted conjugacy class $\<g\>$ must lie in the image of the map $H^1(G',F)\rar{}H^1(G,F)$. In other words, there must be an $h\in G$ such that $hgF(h^{-1})\in G'$, i.e., such that $({}^{hgF(h^{-1})}H,{}^{hgF(h^{-1})}\N)=(H,\N)$. But since $(H,\N)$ is $F$-stable, this means that the admissible pair ${}^{h^{-1}}(H,\N)\in \C$ must be $gF$-stable as desired.

To prove (ii), by passing to pure inner form of $F$ if necessary, we may assume that $g=1$. Suppose that $(H,\N)\in \C$ is $F$-stable. (If there is no such pair in $\C$, then by part (i), $\T_e^1=0$ in $\Fun(G_0(\Fq)/\sim)$, in which case statement (ii) is vacuous.) Let $e'_0\in \D_{G'_0}(G'_0)$ be the corresponding Heisenberg idempotent. By \cite[Prop. 4.12]{B} we have
\beq
\T_e^1=\sum\limits_{\substack{\<g'\>\in H^1(G',F)\\g'=hF(h^{-1}) \h{ for }\\\h{ some }h\in G}}\ind_{{G'}^{g'}_0(\Fq)}^{G_0(\Fq)}\T^{g'}_{e'_0}.
\eeq
On the other hand, for $h\in G$, the pair ${}^{h^{-1}}(H,\N)\in \C$ is $F$-stable if and only if $hF(h^{-1})\in G'$. Statement (ii) now follows. (See also Corollary \ref{c:irrepe}.)
\epf

\subsection{Geometric reduction process for algebraic groups}\label{s:grpag}
In this section we prove a result relating irreducible representations and admissible pairs for any algebraic group $G$ over $\Fqcl$ equipped with an $\Fq$-Frobenius $F:G\rar{}G$. The following result is a straightforward generalization of \cite[Thm. 7.1]{B1}:
\bthm\label{t:wadm} 
Let $G$ be any algebraic group over $\Fqcl$ and let $F:G\rar{}G$ be an $\Fq$-Frobenius. Let $W$ be an irreducible representation of a pure inner form $G^g_0(\Fq)$. Let $(A,\K)$ be a pair consisting of a $gF$-stable connected unipotent normal subgroup $A\normal G$ and a $G$-invariant multiplicative local system $\K\in (A^*)^{gF}$ such that the subgroup $A^{gF}\subset G^g_0(\Fq)$ acts on $W$ by the character $\T_{\K^g_0}:A^{gF}\rar{}\Qlcl^\times$. Then there exists a $gF$-stable admissible pair $(H,\N)\geq (A,\K)$ for $G$ such that $W$ restricted to the subgroup $H^{gF}=H^g_0(\Fq)\subset G^g_0(\Fq)$ contains the associated character $\T_{\N^g_0}:H^g_0(\Fq)\rar{}\Qlcl^\times$ as a direct summand, or equivalently, such that $W$ occurs as a direct summand in $\ind_{H^g_0(\Fq)}^{G^g_0(\Fq)}\T_{\N_0^g}$. 
\ethm
\brk
If we take the pair $(A,\K)=(1,\Qlcl)$, then the result above implies that to any $W\in \Irrep(G,F)$ we can associate some admissible pair compatible with $W$ in the above sense.
\erk
\begin{proof}[Proof of Thm. \ref{t:wadm}]
The proof of \cite[Thm. 7.1]{B1} readily extends to the general case. Namely, let us first consider the case $[G,\R_u(G)]\subset A$. Now the argument from \cite[\S7.5]{B1} completes the proof in this case.

Now suppose that $[G,\R_u(G)]$ is not contained in $A$, i.e. $\R_u(G)/A=\R_u(G/A)$ is not central in $G/A$. Let $A\subset Z\subset \R_u(G)$ be such that $\R_u(Z(G/A))=Z/A$. Then $Z$ must be a proper subgroup of $\R_u(G)$ by our assumption. Now consider the map
$\phi_\K:\R_u(G)/A\rar{}(Z/A)^*$ as defined in \cite[\S7.6]{B1}. Using the argument {\it loc. cit.}, we obtain a $gF$-stable subgroup $A\subsetneq B\subset Z$ with an extension $\K'\in (B^*)^{gF}$ of $\K\in A^*$ such that the pair $(B,\K')$ is compatible with $W$. Let $G_1\subset G$ be the normalizer of the pair $(B,\K')$. If $G=G_1$, we can replace our original pair with the strictly larger pair $(B,\K')$ and proceed. Else $G_1$ is strictly smaller and we can proceed by induction following the argument from \cite[\S7.6]{B1}.
\epf
 
Combining Theorem \ref{t:wadm} with the results from \S\ref{s:irap} we obtain:
\bcor\label{c:irreppart}
Assume once again that $G$ is neutrally solvable equipped with a Frobenius. Then we have a partition
\beq
\Irrep(G,F)=\coprod\limits_{e\in \hG^F_{\h{adm}}}\Irrepe.
\eeq
\ecor

\subsection{Irreducible representations and character sheaves}\label{s:ircs}
In this section we will describe the relationship between irreducible characters and character sheaves on a neutrally solvable group $G$. In view of Corollary \ref{c:irreppart}, let us fix a minimal idempotent $e\in \hG^F_{\h{adm}}$ and study the set $\Irrepe=\Irrep_\C(G,F)$, where $\C\in [\Padm]^F$ is some pre-image of $e$ under the map (\ref{e:padmfhgf}). By passing to a pure inner form of $F$ if needed, let us assume that the geometric conjugacy class $\C$ contains an $F$-stable admissible pair $(H,\N)$ with normalizer $G'$. Let $e'\in \DGp$ be the corresponding Heisenberg idempotent. Then we have $e=\indg e'$ and an equivalence
\beq\label{e:indedgp}
\indg:e'\DGp\rar{\cong}e\DG
\eeq
of $\Qlcl$-linear triangulated braided monoidal categories compatible with the actions of the Frobenius on each side. We also have a $\Qlcl$-linear triangulated equivalence (cf. \cite[\S4.2]{De3})
\beq\label{e:indirrepe}
\indg:e'\D_{G'}(G'F)\rar{}e\D_G(GF)
\eeq
of module categories. Also note that the functor 
\beq
\indg:\D_{G'}(G'F)\rar{}\D_G(GF)
\eeq
identifies with induction of representations of pure inner forms of $G'_0(\Fq)$ (cf. \cite[\S4.3]{De3}). Also we have the `induction of class functions' map (cf. \cite[\S4.3]{De3}):
\beq\label{e:indfun}
\indg:\Fun([G'],F)\rar{}\FunG.
\eeq
Given an object $M'\in \D_{G'}(G'F)$ we have its character $\chi_{M'}\in \Fun([G'],F)$. We have seen {\it loc. cit.} that 
\beq
\chi_{\indg M'}=\indg\chi_{M'}.
\eeq

Also we have the compatible functors
\beq
\ind_{G'_0}^{G_0}:\D_{G'_0}(G'_0) \rar{} \DGn,
\eeq
\beq
\ind_{G'}^{G}:\D_{G'_0}^{\h{Weil}}(G'_0) \rar{} \D_{G_0}^{\h{Weil}}(G_0).
\eeq
By \cite[Prop. 4.12]{B1}, for each $(C',\psi')\in \D_{G'_0}^{\h{Weil}}(G'_0)$ we have 
\beq
\indg(\T_{C',\psi'})=\T_{\indg(C',\psi')}.
\eeq

Since the pair $(H,\N)$ is assumed to be $F$-stable, we have the corresponding Heisenberg idempotent $e'_0\in \D_{G'_0}(G'_0)$ and the idempotent $e_0\in \DGn$. We have equivalences 
\beq
\ind_{G'_0}^{G_0}:e'_0\D_{G'_0}(G'_0) \rar{\cong} e_0\DGn,
\eeq
\beq
\ind_{G'}^{G}:e'\D_{G'_0}^{\h{Weil}}(G'_0) \rar{\cong} e\D_{G_0}^{\h{Weil}}(G_0).
\eeq
By \cite[Lem. 5.1, 5.2]{De3}, (\ref{e:indfun}) gives rise to an isomorphism of Frobenius algebras
\beq\label{e:indefun}
\indg: \T_{e'}\Fun([G'],F)\rar{\cong} \T_e\FunG,
\eeq
where $\T_{e'}\in \Fun([G'],F)$ and $\T_e\in \FunG$ are the corresponding idempotents.
Combining this with (\ref{e:indirrepe}) and Corollary \ref{c:irrepe}, we obtain
\bprop
(i) The equivalence (\ref{e:indirrepe}) induces a bijection $\indg:\Irrep_{e'}(G',F)\rar\cong \Irrep_e(G,F)$.\\
(ii) For a $W\in \Irrep(G,F)$, the following are equivalent: 
\bit
\item[(a)] $W$ lies in $\Irrepe$.
\item[(b)] $\chi_W$ lies in $\T_e\FunG$.
\item[(c)] $W_{\h{loc}}$ lies in $e\D_G(GF)$.
\eit
(iii) The set $\{\chi_W|W\in \Irrepe\}$ is an orthonormal basis of $\T_e\FunG$.
\eprop

Also using (\ref{e:indedgp}) we have an equivalence between the theory of character sheaves in $e'\DGp$ and $e\DG$ compatible with the Frobenius actions.

Now suppose that $f\in \lL_e(G)^F$ is an $F$-stable indecomposable minimal quasi-idempotent in $e\DG$. As in the Heisenberg case, we define the $\lL$-packet of irreducible representations associated with $f$ as the set
\beq\label{e:deflpirrepe}
\Irrep_f(G,F)=\Irrep_{e,f}(G,F):=\{W\in \Irrep(G,F)|W_{\loc}\in {}^f\D_G(GF)={}^f\D_G(GF)^\Delta\}\subset \Irrepe.
\eeq 
By Proposition \ref{p:lpdec} and the results above we have
\beq\label{e:irrepedec}
\Irrepe=\coprod\limits_{f\in \lL_e(G)^F}\Irrep_{e,f}(G,F).
\eeq

By (\ref{e:indedgp}), let $f'$ be the corresponding minimal quasi-idempotent in $e'\DGp$. In this setting we have:
\bprop
The equivalence (\ref{e:indirrepe}) induces an equivalence $\indg:{}^{f'}\D_{G'}(G'F)\rar\cong {}^f\D_G(GF)$ and hence Proposition \ref{p:tstrmodcat} holds for every $e\in \hG^F_{\h{adm}}$. In particular ${}^f\D_G(GF)$ has a distinguished $t$-structure (namely, the perverse $t$-structure shifted by $n_e+\tau_e$, cf. Definition \ref{d:miniinv}) with heart, denoted by $\M_{GF,f}$, which is an invertible $\M_{G,f}$-module category (under truncated convolution) equipped with a natural $\M_{G,f}$-module trace. This invertible modular category corresponds to the modular auto-equivalence $F:={F^*}^{-1}:\M_{G,f}\rar\cong \M_{G,f}$. The map 
\beq 
W\mapsto M_W:=W_{\loc}[\dim G+n_e+\tau_e]\in \M_{GF,f}
\eeq
defines a bijection between the sets $\Irrep_f(G,F)$ and $\O_{\M_{GF,f}}$. Let $\t\T_f:=\indg \t\T_{f'}\in \T_e\FunG$, where $\t\T_{f'}\in \T_{e'}\Fun([G'],F)$ is the idempotent function associated with the quasi-idempotent $f'\in \lL_{e'}(G')^F$ (cf. Remark \ref{r:trf}). Let $W\in \Irrep(G,F)$. Then $W\in \Irrep_f(G,F)$ if and only if $\chi_W\in \t\T_f\FunG$. The set $\{\chi_W|W\in \Irrep_f(G,F)\}$ is an orthonormal basis of the Frobenius algebra $\t\T_f\FunG$.
\eprop
\bpf
The proposition is straightforward using the equivalences (\ref{e:indedgp}) and (\ref{e:indirrepe}). One only needs to observe that the equivalence (\ref{e:indirrepe}) takes local systems in $e'\D_{G'}(G'F)$ to local systems in $e\D_G(GF)$ and then use Definition \ref{d:miniinv}.
\epf

In the setting of the proposition, we have equivalences
\beq
\indg: \M_{G',f'}\rar\cong \M_{G,f} \hbox{ and } \indg: \M_{G'F,f'}\rar\cong \M_{GF,f}
\eeq
of modular categories and their module categories respectively. The induction functor takes the positive trace $\tr_{G',F,f'}^+$ in $\M_{G'F,f'}$ to the positive trace $\tr_{G,F,f}^+$ in $\M_{GF,f}$. Moreover, it is easy to check that for any endomorphism $\alpha\in \D_{G'}(G'F)$ we have $\tr_{G,F}(\indg \alpha)=\tr_{G',F}(\alpha)$ for the natural traces in $\D_{G'}(G'F)$ and $\D_G(GF)$.

Our next goal is to study the relationship between the sets $\CS_f(G)^F$ and $\Irrep_f(G,F)$. As in the Heisenberg case, for each $C\in \CS_f(G)^F$ let us choose $\psi_C:F^*C\rar\cong C$ so as to satisfy $\bigstar$ from \S\ref{s:rbiccs} and let $\t{S}(F,f)$ denote the transition matrix, i.e. the matrix with entries
\beq
\t{S}(F,f)_{C,W}:=\<\T_{C,\psi_C},\chi_{W}\>=\<\T_{C,\psi_C},(-1)^{2d_e}\chi_{M_W}\>=(-1)^{2d_e}\tr_F(\g_{C,\psi_C,M_W})
\eeq
for $C\in \CS_f(G)^F$ and $W\in \Irrep_f(G,F)$. On the other hand, let $S^+(F,f)$ denote the crossed S-matrix associated with the $\M_{G,f}$-module category $\M_{GF,f}$ equipped with the positive module trace. Using the results of this section, in particular (\ref{e:indefun}) along with Theorem \ref{t:relheis} we can now prove
\bthm\label{t:maingen}
(i) We have $\t{S}(F,f)=\frac{q^{d_e}}{q^{\dim G}\cdot \sqrt{\dim \M_{G,f}}}\cdot S^+(F,f)$.\\
(ii) Let $W'\in \Irrep_{f'}(G',F)$ be an irreducible representation of the pure inner form ${G'}^{g'F}$ for some $\<g'\>\in H^1(G',F)$. Then $W:=\ind_{{G'}^{g'F}}^{G^{g'F}} W'\in \Irrep_f(G,F)$ and we have
\beq
\dim W=\frac{\dim^+M_W}{\sqrt{\dim \M_{G,f}}}\cdot\frac{|T^{g'F}|}{|{T'}^{g'F}|}\cdot|\Pi_0^{g'F}|\cdot q^{d_e+\tau_e-\dim T},
\eeq 
where $T$ denotes the torus $G^\circ/\R_u(G)$. Moreover each $W\in \Irrep_f(G,F)$ can be obtained from a unique $W'\in \Irrep_{f'}(G',F)$ by induction. 
\ethm
\bpf
The result follows from the preceding remarks using the fact that $d_{G,e}=d_{G',e'}+\dim G -\dim G'$.
\epf

\subsection{Shintani descent in the general case}\label{s:sdgc}
In this section, we will prove results related to Shintani descent in the general case. Suppose that $e\in \hG^F_{\h{adm}}$. We continue to use all notation as before. In particular, let us assume that $e$ comes from an $F$-stable (by replacing $F$ by a pure inner form if necessary) admissible pair $(H,\N)$ for $G$ with normalizer $G'$. As before, let $\T_e\in \FunG$ denote the associated idempotent.

We have already studied Shintani descent in the Heisenberg case in \S\ref{s:sdthc}. Moreover, we have seen in \S\ref{s:grpag} that we can reduce the study of the general case to the Heisenberg case for the subgroup $G'$ using the induction functor $\indg$. 

Hence we obtain the following result in general:
\bthm\label{t:sdgc}
Let $e\in \hG^F_{\h{adm}}$ as above and let $m$ be a positive integer. Let us continue to use all notations and conventions as above. Then\\
(i) We have an $\lL$-packet decomposition
\beq
\Irrep_e(G,F^m)^F=\coprod\limits_{f\in \lL_e(G)^F}\Irrep_f(G,F^m)^F.
\eeq
(ii) The $m$-th Shintani descent $\Sh_m:\Irrep(G,F^m)^F\hookrightarrow \FunG$ respects $\lL$-packet and idempotent decompositions, namely $\Sh_m(\Irrep_e(G,F^m)^F)\subset \T_e\FunG$ with the image being an orthonormal basis and for each $f\in \lL_e(G)^F$, $\Sh_m(\Irrep_f(G,F^m)^F)\subset \t\T_f\FunG$ with the image again being an orthonormal basis. \\
(iii) For $f\in\lL_e(G)^F$, let $\t{\Sh}_m(F,f)$ denote the transition matrix relating the $m$-th Shintani basis of $\t\T_f\FunG$ as described above and the basis $\Irrep_f(G,F)$ (see also (\ref{e:deftsh})). Let $\Sh^+_m(F,f)$ be the $m$-th Shintani matrix associated with the modular category $\M_{G,f}$ and the module category $\M_{GF,f}$ equipped with the positive module trace. Then 
\beq
\t{\Sh}_m(F,f)=\frac{(-1)^{2d_e}}{\sqrt{\dim \M_{G,f}}}{\Sh}^+_m(F,f).
\eeq
(iv) There exists a positive integer $m_0$ such that for each $f\in \lL_e(G)^F$ and each positive integer $m$, the $m$-th Shintani basis of $\t\T_f\FunG$ only depends (up to scaling by roots of unity) on the residue class of $m$ modulo $m_0$.\\
(v) Let $m$ be a positive multiple of $m_0$. Then the $m$-th Shintani basis of $\T_e\FunG$ agrees with the basis $\left\{\frac{q^{\dim G}}{q^{d_e}}\cdot \T_{C,\psi_C}|C\in \CS_e(G)^F\right\}$ up to scaling by roots of unity.
\ethm

\brk
Let us recall that we have the twist $\theta^{F^m}$ in the category $\D_G(GF^m)$. On the other hand for each $f\in \lL_e(G)^F$ we also have the twist $\theta^{F^m,f,+}$ in $\M_{GF^m,f}$ corresponding to the positive module trace. On restricting to $\M_{GF^m,f}\subset \D_G(GF^m)$ we have (cf. Lemmas \ref{l:cfe}, \ref{l:trnu})
\beq
\theta^{F^m,f,+}=\frac{q^{md_e}}{q^{m\dim G}}\cdot \theta^{F^m}.
\eeq
\erk

\subsection{Conclusion}\label{s:c}
Finally in this section we complete the proof of all our main results in the general case. We have seen that the $\lL$-packets of character sheaves are parametrized by the set $\lL(G)$ of minimal indecomposable quasi-idempotents in $\DG$. We will now show that the $\lL$-packets of $\Irrep(G,F)$ are parametrized by the set $\lL(G)^F$. 
\bthm\label{t:fstap}
Let $f\in \lL(G)^F$ and let $e\in \hG$ be the minimal idempotent (cf. Proposition \ref{p:lleg}) such that $f\in \lL_e(G)$. Then $e\in \hG^F_{\h{adm}}$ and hence the $\f{L}$-packet $\Irrep_{e,f}(G,F)$ is non-empty. We have the $\lL$-packet decomposition
\beq
\Irrep(G,F)=\coprod\limits_{f\in \lL(G)^F}\Irrep_f(G,F).
\eeq 
\ethm
\bpf
Since $f$ is $F$-stable, it is clear that $e\in \hG^F$. Hence we must prove that $e$ in fact can be obtained from an $F$-stable conjugacy class of admissible pairs for $G$. By \cite{De2} we know that $e$ must come from some admissible pair $(H,\N)$ defined over $\Fqcl$. This pair must be defined over some finite field $\F_{q^m}$. Hence $e\in \hG^{F^m}_{\h{adm}}$. In particular we have the $\lL$-packet $\Irrep_{e,f}(G,F^m)$ and the bijection $\Irrep_{e,f}(G,F^m)\cong \O_{\M_{GF^m,f}}$. Since $f$ is $F$-stable, we have an action of $F$ on the $\lL$-packet $\Irrep_{e,f}(G,F^m)$. This action agrees with the action of $F$ on $\O_{\M_{GF^m,f}}$ under the above bijection. We also have a modular action of $F$ on the modular category $\M_{G,f}$ and the $\M_{G,f}$-module category $\M_{GF^m,f}$ corresponds to the action of $F^m$. Hence by \cite[Corollary 3.10]{De5} (used for $a=0$ and $a=m$) implies that $|\O_{\M_{GF^m,f}}^F|=|\O_{\M_{G,f}}^F|$. Hence we see that $\Irrep_{e,f}(G,F^m)^F$ is non-empty. By Corollary \ref{c:irreppart} we have an orthonormal decomposition 
\beq
\FunG=\bigoplus\limits_{e'\in \hG^F_{\h{adm}}}\T_{e'}\FunG.
\eeq 
Now using the fact that $\Sh_m(\Irrep(G,F^m)^F)\subset \FunG$ is a basis and Theorem \ref{t:sdgc}(ii), we see that $e$ must necessarily lie in $\hG^F_{\h{adm}}$. Once we know this, using Corollary \ref{c:irreppart} and (\ref{e:irrepedec}) we obtain the desired $\lL$-packet decomposition.
\epf

Once we have this last result, we can use our previous results about $\hG^F_{\h{adm}}$ to study all $\lL$-packets. Thus in view of Theorems \ref{t:maingen}, \ref{t:sdgc} the proofs of our main results (Theorems \ref{t:main2}, \ref{t:main3}) are now complete.

\subsection{An example}\label{s:example}
Let us workout the character sheaves and characters in an example, namely let $B=\left\{\left( \begin{smallmatrix} t_1 & a & c \\ 0 & t_2 & b\\ 0 & 0 & t_3\end{smallmatrix} \right): t_1t_2t_3=1 \right\}$ be the Borel subgroup of $SL_3(\k)$.

Let $T$ be the diagonal maximal torus of $B$. Let us consider certain 1-parameter subgroups of $T$. Let $T_{12}=\left\{\left( \begin{smallmatrix} t & 0 & 0 \\ 0 & t & 0\\ 0 & 0 & t^{-2}\end{smallmatrix}\right)\right\}, T_{23}=\left\{\left( \begin{smallmatrix} t^{-2} & 0 & 0 \\ 0 & t & 0\\ 0 & 0 & t\end{smallmatrix}\right)\right\}$ and $T_{13}=\left\{\left( \begin{smallmatrix} t & 0 & 0 \\ 0 & t^{-2} & 0\\ 0 & 0 & t\end{smallmatrix}\right)\right\}$. Also, let $\mu_3=\left\{\left( \begin{smallmatrix} \omega & 0 & 0 \\ 0 & \omega & 0\\ 0 & 0 & \omega\end{smallmatrix}\right):\omega^3=1 \right\}\subset T$.

The unipotent radical of $B$ is $U:=\left\{\left( \begin{smallmatrix} 1 & a & c \\ 0 & 1 & b\\ 0 & 0 & 1\end{smallmatrix}\right)\right\}$. The center (as well as commutator subgroup) of $U$ is $Z:=\left\{\left( \begin{smallmatrix} 1 & 0 & c \\ 0 & 1 & 0\\ 0 & 0 & 1\end{smallmatrix}\right)\right\}$. Let $V=U/Z =U^{ab} \cong\f{G}_a^2$. Let us fix a non-trivial character $\psi:\Z/p\Z\hookrightarrow \Qlcl^\times$ and hence we obtain the Artin-Schreier local system $\L$ on $\f{G}_a$. We can then identify the Serre dual of $V$ (which parametrizes multiplicative local systems on $V$) with the `linear dual' $V^*$ via pullback of the Artin-Schreier local system along linear functionals. Hence given a row vector $(x,y)\in V^*$, we have the corresponding multiplicative local system $\L_{(x,y)}$ on $V$. We have a canonical identification of Serre duals $U^*\cong V^*$ and by a slight abuse of notation, we will also consider $\L_{(x,y)}$ as a multiplicative local system on $U$ via the homomorphism $U\onto V$.

Then we have 5 minimal idempotents in $\D_B(B)$ corresponding to the following admissible pairs:
\begin{enumerate}
	\item The Heisenberg admissible pair $(U,\Qlcl)$ normalized by all of $B$.
	\item The admissible pair $(U,\L_{(1,0)})$ with normalizer $T_{12}U$.
	\item The admissible pair $(U,\L_{(0,1)})$ with normalizer $T_{23}U$.
	\item The admissible pair $(U,\L_{(1,1)})$ with normalizer $\mu_3U$.
	\item The admissible pair $(Z,\L)$ with normalizer $T_{13}U$.
\end{enumerate}

\brk
If $\L'$ is any multiplicative local system on $U$, then $(U,\L')$ is an admissible pair for $B$. However, the associated minimal idempotent only depends on the $T$-orbit of $\L'\in V^*\cong U^*$ and hence is covered by one of the first four cases. Similarly, for any non-trivial $\L'\in Z^*$, $(Z,\L')$ is an admissible pair which gives rise to the same minimal idempotent as case 5 above.
\erk

Of the five minimal idempotents in $\D_B(B)$, the cases $4$ and $5$ are more interesting and we will analyse them further. The analysis of the first three cases is simpler.

{\it The admissible pair $(U,\L_{(1,1)})$}: Suppose first that $\h{char} \k \neq 3$. The normalizer $\mu_3U\cong \mu_3\times U$ of this pair is a neutrally unipotent group. The object $e'_{(1,1)}:=\L_{(1,1)}[6](3)$ is a minimal idempotent in $\D_{\mu_3U}(\mu_3U)$ and $e_{(1,1)}:=\ind_{\muU}^B e'_{(1,1)}$ is the minimal idempotent in $\D_B(B)$ corresponding to this admissible pair. 

The category $e'_{(1,1)}\D_{\mu_3U}(\mu_3U)$ is equivalent to $D^b((\Vec_{\mu_3})^{\mu_3})$ where $(\Vec_{\mu_3})^{\mu_3}$ is the Drinfeld double of $\mu_3$. By Appendix \ref{s:qidcfc}, $e'_{(1,1)}$ is the unique minimal indecomposable quasi-idempotent in $e'_{(1,1)}\Dmu$. So in this case we have a unique $\f{L}$-packet of character sheaves in $e'_{(1,1)}\Dmu$ and it has 9 character sheaves. The underlying objects of $\D(\muU)$ for these 9 character sheaves are: the object $e'_{(1,1)}=\L_{(1,1)}[6](3)$ supported on $U$ and its two translates supported on the other connected components $\omega U, \omega^2 U$ of $\muU$. The 9 character sheaves come from the 3 different $\mu_3$-equivariance structures on each of the above 3 sheaves.

Finally we induce to $e_{(1,1)}\D_B(B)$ and obtain analogous results here. The modular category $\M_{B,e_{(1,1)}}$ is equivalent to the Drinfeld double $(\Vec_{\mu_3})^{\mu_3}$. There are 9 character sheaves in this $\f{L}$-packet.

Now let $\Fq$ be a finite field and let $\k=\Fqcl$. Let us consider the finite group $B(\Fq)$. Let $F:B\to B$ be the Frobenius endomorphism. The admissible pair $(U,\L_{(1,1)})$ can also be defined over $\Fq$ and it is clear that the minimal quasi-idempotent $e_{(1,1)}\in \D_B(B)$ is $F$-stable. The modular autoequivalence of $\M_{B,e_{(1,1)}}\cong (\Vec_{\mu_3})^{\mu_3}$, $F:\M_{B,e_{(1,1)}}\rar{}\M_{B,e_{(1,1)}}$ is induced by the Frobenius $F:\mu_3\rar{}\mu_3$.

This leads to two cases:\\
{\it (i) $\Fq$ does not contain the cuberoots of unity:} In this case $F$ interchanges the cuberoots $\omega$ and $\omega^2$ and it can be checked that the unit object $e_{(1,1)}\in \M_{B,e_{(1,1)}}\cong (\Vec_{\mu_3})^{\mu_3}$ is the only $F$-stable character sheaf in the $\f{L}$-packet. Hence there is a unique irreducible representation of $B(\Fq)$ lying in the corresponding $\f{L}$-packet $\Irrep_{e_{(1,1)}}(B(\Fq))$. This is the $(q-1)^2$-dimensional irreducible representation of $B(\Fq)$ obtained by inducing the `generic' character of $U(\Fq)$ corresponding to the $F$-stable multiplicative local system $\L_{(1,1)}$ on $U$.

{\it (ii) $\Fq$ contains the cuberoots of unity:} In this case $F$ acts as the identity on $\mu_3$ and hence all the 9 character sheaves in $\M_{B,e_{(1,1)}}$ are $F$-stable. In this case the $\f{L}$-packet $\Irrep_{e_{(1,1)}}(B(\Fq))$ consists of nine $\frac{(q-1)^2}{3}$-dimensional irreducible representations of $B(\Fq)$. In this case, the action of $T(\Fq)$ on the set of 1-dimensional characters of $U(\Fq)$ has three `generic orbits'. The induction to $B(\Fq)$ of each generic character of $U(\Fq)$ breaks into three irreducible representations. This gives the nine irreducible characters in this $\f{L}$-packet.

Finally, let us also consider the case $\h{char} \k =3$. In fact in this case things are simpler and we have $\M_{B,e_{(1,1)}}\cong \Vec$ and there is a unique character sheaf as well as a unique irreducible representation of $B(\Fq)$ in the corresponding $\f{L}$-packet. This completes the analysis of the admissible pair $(U,\L_{(1,1)})$.

{\it The admissible pair $(Z,\L)$}: The normalizer of this admissible pair is $T_{13}U$. Let $e'=\L[2](1)$ be the corresponding minimal idempotent in $\D_{T_{13}U}(T_{13}U)$ and $e=\ind_{T_{13}U}^B e'$ the associated minimal idempotent in $\D_B(B)$. 

Consider $U\supset H=\left\{\left( \begin{smallmatrix} 1 & a & c \\ 0 & 1 & 0\\ 0 & 0 & 1\end{smallmatrix}\right)\right\}\supset Z$. Let $\L'$ be the multiplicative local system on $H$ obtained by pulling back $\L$ along the evident projection $H\onto Z$. Then $\L'$ is a $T_{13}$-stable extension of $\L$ from $Z$ to $H$. Then by \cite[\S8.4]{De2} $(H,\L')$ is a {\it central} admissible pair for $B$ giving rise to the same minimal idempotent $e\in \D_B(B)$. The normalizer of the pair $(H,\L')$ is $T_{13}H$. Let $e'':=\L'[4](2)$ be the associated minimal idempotent in $\D_{T_{13}H}(T_{13}H)$. We have an identification $e''\D_{T_{13}H}(T_{13}H)\cong e''\D_{H}(H)\boxtimes \D_{T_{13}}(T_{13})\cong \D_{T_{13}}(T_{13})$ since $e''\D_H(H)\cong D^b\Vec$. 

Hence in this case, for each multiplicative local system $\K\in \C(T_{13})$ on the torus $T_{13}$, we obtain the corresponding minimal quasi-idempotent $e_\K\in \D_{T_{13}}(T_{13})$ and consequently the minimal quasi-idempotent $e''\boxtimes e_\K\in e''\D_{HT_{13}}(HT_{13})$. In other words, the set $\f{L}_{e''}(HT_{13})$ of $\f{L}$-packets of character sheaves in $e''\D_{HT_{13}}(HT_{13})$ can be identified with $\C(T_{13})$. Each such $\f{L}$-packet contains a unique character sheaf, namely the sheaf $e''\boxtimes e_\K$. The modular category associated with each $\f{L}$-packet is simply $\Vec$.

Consequently, character sheaves (which coincide with indecomposable minimal quasi-idempotents in this case) in $e\D_B(B)$ are in bijection with $\C(T_{13})$. Namely, given $\K\in \C(T_{13})$, $\ind_{HT_{13}}^B (e''\boxtimes e_\K)$ is the corresponding character sheaf in $e\D_B(B)$.

Now suppose $\Fq$ is a finite field, $\k=\Fqcl$ and $F:B\to B$ the Frobenius. We want to classify the corresponding irreducible characters of $B(\Fq)$. In this case, the admissible pair $(H,\L')$ and the subgroup $T_{13}$ is defined over $\Fq$. Hence the $F$-stable character sheaves $C$ in $e\D_B(B)$ correspond to $\C(T_{13})^F$ which in turn are in bijection with  characters $\chi$ of $T_{13}(\Fq)$. Given an $F$-stable character sheaf $C$ as above, we have the corresponding irreducible representation of $B(\Fq)$ constructed as follows: Let $\chi:T_{13}(\Fq)\to \Qlcl^\times$ be the character corresponding to $C$ as above and let $\psi_{\L'}:H(\Fq)\to \Qlcl^\times$ be the character corresponding to the $F$-stable multiplicative local system $\L'$ on $H$. Then the irreducible representation corresponding to $C$ is $\ind_{HT_{13}(\Fq)}^{B(\Fq)}\psi_{\L'}\chi$. It is $q(q-1)$-dimensional. Thus we see that $\Irrep_e(B(\Fq))$ contains $q-1$ $q(q-1)$-dimensional irreducible representations and each $\f{L}$-packet in $\Irrep_e(B(\Fq))$ is singleton.

\appendix
\section{Appendix: Quasi-idempotents in linear triangulated categories}\label{a:qiltc}
In this appendix we introduce and study the notion of quasi-idempotents in $\Qlcl$-linear triangulated categories. This notion plays a key role in our approach to the theory of character sheaves on neutrally solvable groups developed in the main text. 

\subsection{Definition of quasi-idempotents}
Let $(\D,\ast,\un)$ be any $\Qlcl$-linear triangulated monoidal category (cf. \cite[\S2.3.2]{De2}). Then we have an action of $D^b\Vec$ on $\D$ (which we will denote by $\otimes$). Since $\D$ is monoidal, this action can also be thought of in terms of the canonical (central) monoidal functor $(D^b\Vec,\otimes, \Qlcl)\rar{}(\D,\ast,\un)$. We note that for $V,W\in D^b\Vec$ and $X,Y\in \D$ we have functorial isomorphisms $(V\otimes X)\ast(W\otimes Y)\rar\cong (V\otimes W)\otimes (X\ast Y)$ and $(V\otimes W)\otimes X\rar\cong V\otimes (W\otimes X).$ We will often skip the parenthesis from our notation. The triangulated category $D^b\Vec$ has the standard $t$-structure, and for $V\in D^b\Vec$, $i\in \Z$, we let $V^i$ denote the $i$-th cohomology of $V$.

Suppose the triangulated category $\D$ is such that each object has a unique direct sum decomposition into indecomposable objects. Furthermore suppose that all shifts of indecomposable objects of $\D$ are non-isomorphic to each other. Then the condition that $V\otimes X\cong W\otimes X$ for some $V,W\in D^b\Vec$, $X\in \D$ implies that either $V\cong W$ in $D^b\Vec$ or $X=0$. 

\brk\label{r:ksa}
From now on we will always assume that $\D$ satisfies these properties. In the main text we will apply our results from this appendix to triangulated categories of $\Qlcl$-complexes. These categories satisfy the above properties.
\erk

\bdefn\label{d:monic}
We let $D^{\leq 0}_{mon}\Vec\subset D^b\Vec$ denote the full subcategory formed by objects $V\in D^b\Vec$ such that $V^i=0$ for each $i>0$ and $V^0\cong\Qlcl$. We say that an object $V$ of $D^b\Vec$ is monic if $V\in \Dmon$.
\edefn

\bdefn\label{d:qi}
We say that an object $e\in \D$ is a quasi-idempotent if $e\ast e\cong V\otimes e$ for some $V\in \Dmon$. In this case we say that $e$ is a $V$-quasi-idempotent. For such a quasi-idempotent, we let $^e\D\subset \D$ be the full subcategory formed by objects $X\in \D$ such that $e\ast X\cong V\otimes X$. Note that $^e\D$ is a right ideal in $\D$ and in particular it is a semigroupal category. Also, let us denote by ${}^{e}\D^{\Delta}\subset \D$ the full triangulated subcategory generated by ${}^e\D$. Then ${}^e\D^\Delta$ is also a right ideal in $\D$ and hence a triangulated semigroupal category.
\edefn

\brk\label{r:ed}
Note that we have a functor $e\ast(\cdot):\D\rar{} {}^e\D$. In general this functor may not be essentially surjective.
\erk

By Remark \ref{r:ksa}, if $e\in \D$ is a nonzero quasi-idempotent, then there is a unique $V\in D^b\Vec$ up to isomorphism such that $e\ast e\cong V\otimes e$. Consequently this $V$ must be monic and $e$ must be a $V$-quasi-idempotent.

\brk\label{r:moreqi}
If $e\in \D$ is a quasi-idempotent, then we can construct more quasi-idempotents from it as follows: Let $W\in \Dmon$.  Then if $e\in \D$ is a $V$-quasi-idempotent (for some $V\in \Dmon$), then $W\otimes e\in \D$ is a $W\otimes V$-quasi-idempotent. 
\erk

\blem\label{l:wv'}
Let $e,e'\in \D$ be objects such that $e'$ is a $V'$-quasi-idempotent (for some $V'\in \Dmon$) and $e\ast e'\cong W\otimes e\neq 0$ for some $W\in D^b\Vec.$ Then $W\cong  V'$.
\elem
\bpf
Convolving the given isomorphism with $e'$, we obtain
\beq
W\otimes (e\ast e')\cong (W\otimes e)\ast e' \cong (e\ast e')\ast e' \cong e\ast (e'\ast e')\cong e\ast (V\otimes e')\cong V'\otimes (e\ast e')
\eeq
and we have assumed that $e\ast e'\neq 0$. Hence by Remark \ref{r:ksa} we must have $W\cong V'$.
\epf

We say that an object $X\in \D$ is weakly central if $X\ast Y\cong Y\ast X$ for each $Y\in \D$. We do not require the existence of coherent braiding isomorphisms. In this paper we will mostly be concerned with weakly central quasi-idempotents.

\bdefn\label{d:minqi}
Let $e$ (resp. $e'$) be a nonzero weakly central $V$(resp. $V'$)-quasi-idempotent for some monic $V$ (resp. $V'$). Then $e\ast e'$ is a $V\otimes V'$-quasi-idempotent. We say that $e\leq e'$ if $e\ast e'\cong W\otimes e$ for some nonzero $W\in D^b\Vec$. (By Lemma \ref{l:wv'} we must have $W\cong V'$.) We say that $e\sim e'$ or that $e$ and $e'$ are equivalent if $e\leq e'$ and $e'\leq e$. We say that a nonzero weakly central quasi-idempotent $e$ is minimal if $e\ast e'=0$ or $e\leq e'$ for each weakly central quasi-idempotent $e'\in \D$. 
\edefn

\blem\label{l:ede'd}
Let $e,e'$ be nonzero weakly central quasi-idempotents in $\D$ and let $V,V'\in \Dmon$ be the corresponding monic objects. Then $e\leq e'$ if and only if ${}^e\D\subset {}^{e'}\D$. Moreover suppose that $e\sim e'$. Then we must have $V'\otimes e\cong V\otimes e'$ and that $^e\D= { }^{e'}\D$ as full subcategories of $\D$. 
\elem
\bpf
Let $e\leq e'$, i.e. $e\ast e'\cong V'\otimes e$. Suppose $X\in {}^e\D$, i.e.  $e\ast X \cong V\otimes X$. From this we obtain 
$$e'\ast e\ast X\cong V\otimes e'\ast X \mbox{ and}$$
$$V'\otimes e\ast X\cong  V'\otimes V\otimes X.$$ Now using $e'\ast e\cong e\ast e'\cong V'\otimes e$ we obtain that $$V\otimes e'\ast X\cong V\otimes V'\otimes X.$$ 
Hence using Remark \ref{r:ksa} we can deduce that $e'\ast X\cong V'\otimes X$, i.e. $X\in {}^{e'}\D$. Hence we must have $^e\D\subset { }^{e'}\D$ as full subcategories of $\D$.

On the other hand, suppose that ${}^e\D\subset {}^{e'}\D$. In particular we must have $e\in {}^{e'}\D$. Hence $e'\ast e\cong V'\otimes e$, i.e. $e\leq e'$.

Now suppose that $e\sim e'$. Then by Definition \ref{d:minqi}, we must have $V'\otimes e\cong e\ast e'\cong e'\ast e\cong V\otimes e'$.
\epf

As a consequence we obtain
\bcor\label{c:altminqi}
Let $e\in \D$ be a weakly central quasi-idempotent. Then $e$ is minimal if and only if any nonzero weakly central quasi-idempotent $e'\in \D$ that lies in the full subcategory ${}^e\D$ is equivalent to $e$.
\ecor

\subsection{Quasi-idempotents in derived categories of fusion categories}\label{s:qidcfc}
Let $(\M,\ast,\un)$ be a ($\Qlcl$-linear) fusion category.   In this section we will describe all the quasi-idempotents in the category $D^b\M$. First we prove that there are no non-trivial idempotents in $D^b\M$. 

Note that the category $D^b\M$ has a natural $t$-structure whose heart is $\M\subset D^b\M$. For an object $M\in D^b\M, i\in \Z$, we let $M^i\in \M$ denote its $i$-th cohomology with respect to this $t$-structure. Since $\M$ is semisimple,  $D^b\M$ is also in fact a semisimple abelian category whose simple objects are the shifts of the simple objects of $\M$. We have $M\cong \bigoplus\limits_{i\in \Z}M^i[-i]$ for all objects $M\in D^b\M$.

\blem\label{l:idindbm}
Let $\M$ be a fusion category as above. There are no idempotents in $D^b\M$ other than $0$ and $\un$.
\elem
\bpf
Let us extend the notion of Frobenius-Perron dimension from $\M$ to all of $D^b\M$ such that the shift preserves the Frobenius-Perron dimensions of objects. Let $e\in D^b\M$ be a nonzero idempotent, i.e. $e\ast e\cong e\neq 0$. Hence we must have $\FPdim(e)^2=\FPdim(e)$. Now by \cite[\S8]{ENO} $\FPdim(e)\geq 1$ since $e$ is a nonzero object. Hence we conclude that $\FPdim(e)=1$ and furthermore that $e$ must in fact be a simple object. From the isomorphism $e\ast e\cong e$ we deduce that $e$ must in fact be a simple object of $\M\subset D^b\M$. Moreover we have $\FPdim(e^*)=1$ and $\FPdim(e^*\ast e)=1$ and $\un$ is a direct summand of $e^*\ast e$. Hence we must have that $e^*\ast e\cong \un$, i.e. $e$ is invertible. Hence we deduce that $e\cong \un$. 
\epf

\bprop\label{p:qidbm}
The quasi-idempotents in the category $D^b\M$ are all of the form $V\otimes \un$ for some $V\in\Dmon$.
\eprop
\bpf
Let $V\in \Dmon$ and let $e\in D^b\M$ be a nonzero $V$-quasi-idempotent. We have
\beq\label{e:eeve}
e\ast e \cong V\otimes e.
\eeq
By comparing the highest cohomological degree in which both the sides are nonzero, we conclude that $e^i=0$ for each $i>0$ and that $e^0$ is nonzero. Now we have
\beq
e\cong\bigoplus\limits_{i\geq 0}e^{-i}[i] \mbox{ and }V\cong\Qlcl\oplus\bigoplus\limits_{i\geq i}V^{-i}[i].
\eeq
Expanding (\ref{e:eeve}) and looking at the 0-th cohomological degree we conclude that $e^0\ast e^0\cong e^0$ in $\M$.  By Lemma \ref{l:idindbm} we see that  $e^0\cong \un$. Now by successively looking at each of the negative cohomological degrees of (\ref{e:eeve}) and proceeding by induction, we deduce that $e^{-i}\cong V^{-i}\otimes e^0\cong V^{-i}\otimes \un$ for each $i\in \Z$. Hence $e\cong V\otimes \un$ as desired.
\epf

\bcor
The unit object $\un\in D^b\M$ is a minimal quasi-idempotent and there is a unique equivalence class of minimal quasi-idempotents in $D^b\M$.
\ecor

\subsection{Induction of quasi-idempotents satisfying the Mackey condition}\label{a:iqimc}
Let $G$ be any algebraic group over $\k$. Let $G'\subset G$ be a closed subgroup. In this section we consider quasi-idempotents $e$ in $\DGp$ satisfying the (geometric) Mackey criterion with respect to $G$ and the induced quasi-idempotents $\indg e$ in $\DG$. We will use the ($\Qlcl$-linear triangulated weakly semigroupal) induction functor $\indg:\DGp\rar{}\DG$ and its properties described in \cite{De2}, \cite{BD}.

The following result is a straightforward generalization of some results from \cite[\S4]{De2} and \cite[\S5]{BD}.
\bprop\label{p:indmac} (cf. \cite[Prop. 4.3]{De2})
Let $e\in \DGp$ be a $V$-quasi-idempotent such that $e\ast \delta_x\ast e=0$ for each $x\in G-G'$, where $\delta_x$ is the delta sheaf supported at $x$ and where we consider $e$ as an object of $\D(G)$ by extension by zero. Then we have:\\
(i) The object $f:=\indg e\in \DG$ is a $V$-quasi-idempotent.\\
(ii) If $M\in {}^e\DGp$, then $\indg M\in {}^f\DG$. The $\Qlcl$-linear functor
\beq
\indg: {}^{e}\DGp\rar{}{}^{f}\DG
\eeq
is (strongly) semigroupal and braided.\\
(iii) For each $N\in \DG$, we have a functorial isomorphism $e\ast N\rar\cong e\ast N_{G'}$ where we consider all the $\Qlcl$-complexes involved as objects in $\D_{G'}(G)$, where $N_{G'}\in \DGp$ denotes the restriction of $N$ to $G'$. In particular this means that $e\ast N$ has support contained inside $G'$ and we may consider $e\ast N\cong e\ast N_{G'}$ as an object of ${}^e\DGp$.\\
(iv) For each $N\in \DG$ we have a functorial isomorphism $f\ast N\rar\cong \indg(e\ast N_{G'})$.\\
(v) For each $M\in {}^e\DGp$ we have functorial isomorphisms 
\beq
e\ast (\indg M)|_{G'}\rar{\cong}e\ast \indg M\rar{\cong}e\ast M.
\eeq
\eprop
\bpf
For any $x\in G-G'$ and $M\in {}^e\DGp$, we have $e\ast \delta_x\ast e=0$ and hence $e\ast \delta_x\ast e\ast M\cong V\otimes e\ast \delta_x\ast M=0$. Hence $ e\ast \delta_x\ast M=0$. Similarly if $M'\in {}^e\DGp$ we can prove that $M'\ast \delta_x\ast M=0$. Now statements (i) and (ii) follow from \cite[Prop. 3.15]{De2}. The proofs of statements (iii), (iv) and (v) are similar to the proofs of \cite[Lemma 5.16 and Prop. 5.14]{B1}.
\epf

\subsection{$t$-structures on certain triangulated categories}\label{s:sasha}
We will need the following lemma about $t$-structures on certain special triangulated categories. The proof below was communicated by A. Beilinson.
\blem\label{l:sasha}
Let $\D$ be a triangulated category equipped with a non-degenerate bounded $t$-structure with heart $\A$. Suppose further that for any non-zero objects $A,B\in \A$, we have $\Hom(A,B)\neq 0$. Then any other non-degenerate bounded $t$-structure on $\D$ is a translate of the original one. 
\elem
\bpf
For $i\in\f{Z}$ let $H^i:\D\to\A$ be the $i$-th cohomology functor with respect to the given $t$-structure. Let $\B$ denote the heart of another non-degenerate bounded $t$-structure on $\D$. It is enough to prove that $\B$ lies in a translate of $\A$. Let $X$ be any non-zero object of $\B$. Suppose that $H^M(X), H^m(X)\in \A$ are the top and bottom (and hence $M\geq m$) non-vanishing cohomologies of $X$ with respect to the original $t$-structure. By our assumption, we have a non-zero morphism $f:H^M(X)\rar{}H^m(X)$ in $\A$. Look at the composition 
\beq
c:X\rar{} H^M(X)[-M]\xto{f[-M]} H^m(X)[-M]\rar{}X[m-M]
\eeq
where the first and the last map of this composition are the canonical truncation morphisms. Note that $H^M(c)=f$ is non-zero by construction and hence $c:X\rar{}X[m-M]$ is non-zero. But $X$ lies in the heart $\B$ of the other $t$-structure. The axioms of a $t$-structure then force that $m-M\geq 0$. Hence we conclude that $M=m$ and hence that $X$ lies in a translate of $\A$. This completes the proof. 
\epf

\brk\label{r:sasha}
If the heart $\A\subset \D$ above is such that it has only one simple object, and if every object of $\A$ has finite length, then $\A$ satisfies the hypothesis of the lemma. Hence the lemma will be applicable in the various cases that arise in the main text.
\erk

\end{document}